%% file: main_arxiv.tex
\documentclass[lettersize,journal]{IEEEtran}

\usepackage[table,xcdraw]{xcolor}

\usepackage{amsmath,amsfonts}
\usepackage{array}
\usepackage{textcomp}
\usepackage{stfloats}
\usepackage{url}
\usepackage{verbatim}
\usepackage{graphicx}
\usepackage{cite}
\usepackage{siunitx}
\usepackage{amssymb}
\usepackage{booktabs}
\usepackage{colortbl}
\usepackage{color}
\usepackage{changepage}
\usepackage{algcompatible}
\usepackage{algpseudocode}
\usepackage{tikz}
\usetikzlibrary{positioning}
\usepackage{subfigure}
\usepackage[section]{placeins}
\usepackage[T1]{fontenc}
\usepackage{multirow}
\usepackage[pagebackref,breaklinks,colorlinks]{hyperref}

\usepackage[nohead,margin=1.0in]{geometry}
\usetikzlibrary{plotmarks}
\usepackage{pgfplots}
\pgfplotsset{compat=1.16} 
\usepackage{mathtools}
\usepackage{commath}
\usepackage[normalem]{ulem}
\usepackage{algorithm2e, setspace}
\usepackage{algpseudocode}
\usepackage{dsfont}
\usepackage{etoolbox}
\usepackage{svg}

\makeatletter
\patchcmd{\@algocf@start}%
  {-1.5em}%
  {0pt}%
  {}{}%
\makeatother

\usepackage[capitalize]{cleveref}
\crefname{section}{Sec.}{Secs.}
\Crefname{section}{Section}{Sections}
\Crefname{table}{Table}{Tables}
\crefname{table}{Tab.}{Tabs.}

\numberwithin{equation}{section}

\newcommand{\R}{\mathbb{R}}

\newcommand\A{\mathcal{A}}
\newcommand\X{\mathcal{X}}
\newcommand\Y{\mathcal{Y}}
\newcommand\K{\mathcal{K}}
\renewcommand\L{\mathcal{L}}
\newcommand\M{\mathcal{M}}
\renewcommand\H{\mathcal{H}}

\newif\ifIncludeSupplementary%
\IncludeSupplementarytrue%

\newif\ifArxiv%
\Arxivfalse

\newcommand{\argmax}{\mathop{\rm argmax}}

\definecolor{hlbest}{rgb}{0.9,0.9,0.9}

\newcommand{\email}[1]{\nolinkurl{#1}}

\begin{document}

\title{Learning-Based Approaches for Reconstructions with Inexact Operators in nanoCT Applications}

\author{Tom Lütjen, Fabian Schönfeld, Alice Oberacker, Johannes Leuschner, Maximilian Schmidt, Anne Wald and Tobias Kluth%
\thanks{T.~Lütjen, F.~Schönfeld, J.~Leuschner, M.~Schmidt, and T.~Kluth are with the Center for Industrial Mathematics, University of Bremen, Germany (e-mail: \email{{jleuschn,maximilian.schmidt,tkluth}@uni-bremen.de}).}%
\thanks{A.~Oberacker is with Saarland University, Germany, (e-mail: \email{alice.oberacker@num.uni-sb.de}).}%
\thanks{A.~Wald is with the Institute for Numerical and Applied Mathematics, University of G\"ottingen (e-mail: \email{a.wald@math.uni-goettingen.de}).}
}

\markboth{}%
{}

\maketitle

\begin{abstract}
Imaging problems such as the one in nanoCT require the solution of an inverse problem, where it is often taken for granted that the forward operator, i.e., the underlying physical model, is properly known. 
In the present work we address the problem where the forward model is inexact due to stochastic or deterministic deviations during the measurement process. 
We particularly investigate the performance of non-learned iterative reconstruction methods dealing with inexactness and learned reconstruction schemes, which are based on U-Nets and conditional invertible neural networks. The latter also provide the opportunity for uncertainty quantification. 
A synthetic large data set in line with a typical nanoCT setting is provided and extensive numerical experiments are conducted evaluating the proposed methods.

\end{abstract}


\section{Introduction}

\IEEEPARstart{T}{omographic} X-ray imaging on small scales such as on the nano scale - in short: nanoCT - is an important imaging technique that allows a visualization of the inner structure of small objects in the micro- or nanometer range with a suitable resolution. In contrast to computerized X-ray tomography on larger scales, for example in medical imaging, small vibrations that are transferred from the environment as well as manufacturing tolerances make a significant impact on the data: the resulting relative motion between object and scanner lead to motion artifacts when classical reconstruction techniques such as the filtered backprojection are applied. In this sense, nanoCT is an intrinsically dynamic inverse problem, i.e., an inverse problem that is heavily underdetermined since only one projection per state of the object is available. Hence, the unknown dynamic, resp. inaccuracy in the forward map, needs to be taken into account.

Solving inverse problems typically relies on a proper formulation of the forward operator.
But often and as mentioned before the forward map/model is not known with sufficient accuracy. 
Typical examples can be found in several imaging applications such as magnetic particle imaging \cite{Gleich2005} and also in CT applications \cite{rtkms17,youssef_2013} as addressed in the present manuscript.
Different directions taking into account the model's inexactness explicitly and/or implicitly have been considered in the past. 
An explicit treatment is, for example, the total least squares approach \cite{golub1999tikhonov} explicitly taking into account a deviation on the operator when solving the original linear problem. 
This has also been extended and investigated with respect to regularization properties for general bilinear operators \cite{Bleyer:2013cw, KluthBathkeJiangMaass2020} and also for a combination of model- and data-based prior information \cite{KluthBathkeJiangMaass2020}.
The original total least squares approach can also be equivalently formulated by a problem taking the inexactness implicitly into account.
While this approach relies on the minimization of a certain type of functional, there also exists a class of iterative regularization methods taking into account the inexactness implicitly. 
In \cite{Blanke_2020}, an iterative reconstruction technique is proposed that relies on sequentially projecting iterates onto suitably defined subsets of the solution space, which are designed in such a way that  the solution set is included in each of these sets and measurement noise as well as modeling errors are reflected in the design of the subsets.  %
In the case of moving objects in CT, the unknown motion is interpreted as an inexactness with respect to a simpler model such as the (standard) Radon transform. This approach has been shown to yield promising results, at least for simulated data. \\
More recently, a learning-based approach has been proposed to explicitly correct an inexact operator for the purpose of parameter reconstruction \cite{lunz_learned_2021}, which still relies on the knowledge of the exact forward operator.
In general, several supervised learning-based approaches such as, for example, unrolled iteration schemes \cite{gregor_learning_2010,adler2017solving, adler_learned_2018} can be interpreted as methods taking into account operator inexactness implicitly. Typically, the modeled and potentially inexact operator information is somehow incorporated into the structure of the trainable reconstruction method and the training is then performed with data pairs obtained from the true operator.

In the present work we address the problem of image reconstruction with an inexact operator in the framework of nanoCT applications.
On the algorithmic side we consider sequential algorithms which take the inexactness of the operator into account and which still suffer from artifacts due to the inexactness. These methods are then extended by a learned post-processing scheme to obtain improved reconstructions less prone to the operator inexactness. In this context we use classical post-processing schemes via a UNet and we further develop new approaches based on conditional invertible neural networks, which allow for additional feature extraction for the purpose of uncertainty quantification. Furthermore, in order to evaluate the proposed methods a large synthetic dataset adapted to the nanoCT problem is generated. 

The work is structured as follows: 
In \cref{sec:problem} the general setting and the problem is specified. 
\cref{MaterialsMethods} then provides a detailed description of the non-learned sequential algorithms taking into account the operator inexactness, the learned post-processing schemes via neural networks with a particular focus on conditional invertible neural networks, and the data set generation. 
The methodological part is followed by extensive numerical experiments in \cref{sec:experiments} and we conclude with a discussion in \cref{sec:discussion}.

\section{Problem description and setting} 
\label{sec:problem}

\noindent The small scale of the object in the low micrometer or upper nanometer range means that even small disturbances during the measuring process have a relatively large impact on the scanning geometry and can lead to inconsistent data. This means that the position of the object undergoes changes during data acquisition due to the resulting unavoidable relative motion between object and tomograph. For example, the X-radiation used in the scanner is generated by the interaction of an electron beam with a tungsten needle. Both the electron beam and the position of the needle may vary slightly. However, the largest impact on the relative motion is caused by the environment, when vibrations are transmitted to the tomograph, or by manufacturing tolerances in the tomograph itself, such as a thermal drift or an imprecise motion of the rotating stage on which the object is placed \cite{Dremel_2018, Wang_2019}.
Additionally, the rotation axis may be tilted and the detector is skewed.  %
In order to prevent reconstructions of being impaired by strong artifacts due to the relative motion and geometry deviations, this inconsistency has to be compensated in the reconstruction process.

We use the following notation to outline our methods and results:
The unknown density function of the tested object (i.e., the ground truth) is denoted by $x \in \mathcal{X}$, the measurable data by $y \in \mathcal{Y}$, where $\mathcal{X}$ and $\mathcal{Y}$ are suitable vector spaces, e.g., scalar functions $\mathcal{X}\subset \{x:\Omega\subset \R^2 \to \R \}$ or $\mathcal{X}=\R^{n\times n}$ representing an image. The model describing the relation between density and data is denoted by $\mathcal{A}$, and we refer to it as the forward operator. Hence, we consider the underlying inverse problem as an operator equation
\begin{equation} \label{eq:ip}
 \mathcal{A}x=y, \quad \mathcal{A}: \mathcal{X} \to \mathcal{Y}.
\end{equation}
Noisy data $y^{\delta} := y +\xi$ with an additive noise vector $\xi$ is supposed to have a noise level $\delta$, such that $\left\lVert y^{\delta} - y \right\rVert_{\mathcal{Y}} \leq \delta.$ 
If we use an inexact forward operator $\mathcal{A}^{\eta}$, we assume that the respective error in the model is given by $\eta \in [0,\infty)$ with
\begin{displaymath}
    \left\lVert \mathcal{A}^{\eta} - \mathcal{A} \right\rVert_{\mathcal{X}\to\mathcal{Y}} \leq \eta,
\end{displaymath}
i.e., $\mathcal{A}^0=\mathcal{A}$.
The main goal is now to reconstruct an approximation $\Tilde{x}$ of $x$ from \eqref{eq:ip} for unknown $\mathcal{A}$ and $y$ but for given inexact $\mathcal{A}^\eta$ and potentially noisy $y^\delta$, i.e., from the perspective of the given $\mathcal{A}^\eta$ we have access to the corrupted measurement
\begin{equation}
y^{\delta} = \underset{=:y^\eta}{\underbrace{\mathcal{A}^\eta x}} + (\mathcal{A} - \mathcal{A}^\eta) x + \xi
\end{equation}
only, where $y^\eta$ would be a measurement obtained from $\mathcal{A}^\eta$.
Those desired general reconstruction methods for the inverse problem \eqref{eq:ip} are denoted by $\mathcal{T}: \mathcal{Y} \to \mathcal{X}$ in the remainder. Here, classical methods are combined with learning-based components in terms of neural network schemes $F_{\theta}:\mathcal{X} \to \mathcal{X}$ with network parameters~$\theta$.

\section{Materials and methods}
\label{MaterialsMethods}
\noindent Classical reconstruction methods for CT include analytical inversion formulas, such as the standard filtered back-projection (FBP), as well as iterative reconstruction algorithms \cite{natterer_mathematical_2001}.
However, these methods do not take into account potential inexactness of the forward operator as occurring for nanoCT.
One classical iterative scheme is given by the Kaczmarz algorithm \cite{kaczmarz1937angenaherte}, which in the context of CT is known as the algebraic reconstruction technique (ART) \cite{Gordon1970art}.
The Kaczmarz algorithm forms the basis for both non learning-based algorithms described in the following \cref{sec:classic_iterative}, which extend it in different ways in order to handle operator inexactness.
Then we consider the learning-based extension in \cref{sec:NN_postporcessing}.

\subsection{Sequential algorithms explicitly taking into account inexactness}
\label{sec:classic_iterative}
\subsubsection{Dremel}
\label{sec:dremel}
An approach for explicit geometry correction during the Kaczmarz algorithm was described in \cite{Dremel_2018}, which we will refer to as the \textit{Dremel} method.
It interleaves the iterations of Kaczmarz with correction steps that adapt the previously assumed geometry by maximizing the cross-correlation between the measured data and the forward-projection of the intermediate reconstruction over the plane of possible detector shifts.
We focus on the case in which shifts are identified for each angle independently.
While the correction step in the Dremel method naturally determines shifts in the detector plane, it can also be used to estimate object shifts (or source shifts) as they can be approximately translated into each other as described in \cite{Dremel_2018}.
We use the Dremel method to estimate the object shifts while reconstructing from the perturbed data as we consider data sets featuring random object shifts (\cref{sec:dataset}).
The full algorithm for the Dremel method is specified in the Supplementary Material~\ref{app:dremel} (\cref{alg:dremel}).

\subsubsection{RESESOP-Kaczmarz}
\label{sec:resesop}
The combination of regularizing sequential subspace optimization (RESESOP) as introduced in \cite{narkiss2005sequential}, and Kaczmarz' method, in short RESESOP-Kaczmarz, was first introduced in \cite{Blanke_2020} for linear semi-discrete inverse problems.
For a mathematical analysis of this method we refer to \cite{Blanke_2020}.
Within the present manuscript we write \textit{RESESOP} to refer to RESESOP-Kaczmarz. 
The general idea of the method is outlined in the following.
The respective algorithm in the version used for this work can be found in the Supplementary Material~\ref{app:resesop} (\cref{alg:resesop}).

We need to preface a few necessary concepts. 
Let us first define a set of linear mappings $\A_{k,l}: \X \rightarrow \Y_{k,l}$, %
where $k \in \K \coloneqq \{0, \ldots, K-1\}$ and $l \in \L \coloneqq \{0, \ldots, L-1\}$ are the index sets of scanner angles and detector points of a CT image.
$\A_{k,l}$ with $k \in \K$ and $l \in \L$ describes an X-ray of one scanner angle $k$ and one sensor $l$ on the detector plane. 
By $\A^*$ we denote the adjoint of an operator $\A$.
The semi-discrete inverse problem reads
\begin{equation}
   \A_{k,l}x = y_{k,l}, k \in \K, l \in \L,
\end{equation}
and the solution set can be defined as 
\begin{equation}
    \M_{\A, y} \coloneqq \{ x \in \X: \A_{k,l}x = y_{k,l} \text{ for all } k \in \K, l \in \L \}.
\end{equation}
Let $u \in \X \setminus \{ 0 \}$ and $\alpha, \xi \in \R$. A hyperplane is given by
\begin{equation}
    \H(u,\alpha) \coloneqq \{ x \in \X: \langle u, x \rangle = \alpha \}
\end{equation}
and a stripe is defined by 
\begin{equation}
    \H(u,\alpha, \xi) \coloneqq \{ x \in \X: |\langle u, x \rangle - \alpha| \leqslant \xi \}.
\end{equation}
The algorithm then includes multiple concepts:

\underline{SESOP component}
Sequential subspace optimization (SESOP) is an iterative method which can apply several search directions in each iteration step. %

The full iteration step for SESOP is given by
\begin{equation}
    x_{n+1} = x_n - \sum_{i \in I_n} \tilde{t}_{n,i} \A^\ast w_{n,i},
\end{equation}
with a chosen finite index set $I_n$, $w_{n,i} \in \Y$, search directions $\A^\ast w_{n,i}$, and $\tilde{t}_{n,i}$, $i\in I_n$, minimizing the function 
\begin{equation}
    h_n(t) \coloneqq \frac{1}{2} \norm{x - \sum_{i\in I_n} t_i u_i}^2 + \sum_{i\in I_n} t_i \alpha_i.
\end{equation}
This procedure is equivalent to calculating the metric projection $\mathcal{P}_{\H}(x_n)=: x_{n+1}$ of the current iterate $x_n$ onto the intersection of hyperplanes $\H := \bigcap_{i\in I_n} \H(u_i,\alpha_i)$ with $u_i=\A^\ast w_{n,i}$ and $\alpha_i=\langle w_{n,i}, y \rangle$. 

\underline{Kaczmarz component}
SESOP alone does not take into account the nature of a semi-discrete problem $\A_{q}x = y_{q}$, $q=0,...,Q-1$, (e.g., here $q=(l,k)$, $Q=K\cdot L$) since it is not defined for different realizations $\A_{q}$ of the forward mapping $\A$.
Kaczmarz' method, however, requires such a discretized setting, and it works by iteratively projecting onto the solution set of one realization $\A_{q}$ of $\A$. %
The iteration step can then be written as
\begin{equation}
    x_{n+1} = \mathcal{P}_{\M_{\A_{[n]}, y_{[n]}}}(x_n)%
\end{equation}
where $[i] = i$ mod $Q$. %
Kaczmarz' method thus makes it possible to merge subproblems in RESESOP-Kaczmarz.

\underline{Regularization}
The advantage of combining the SESOP algorithm with Kaczmarz' method is that the iteration can be adapted easily to regularizing semi-discrete inverse problems with modeling errors and noise levels depending on the subproblem. %
For inverse problems with inexact forward operator $\A^\eta$ and noisy data $y^\delta$, it is important to introduce a form of regularization.
As has been described in \cite{Blanke_2020}, SESOP-Kaczmarz can be regularized by replacing the hyperplanes $\H$ by stripes with a width defined by $\eta$ and $\delta$. Such a stripe is defined by
\begin{equation}
\begin{split}
    \H^{\delta, \eta, \rho} &\coloneqq \H\left(( \A^\eta ) ^* w,\langle w, y^\delta \rangle, (\delta + \eta \rho) \norm{w} \right) \\ 
    &= \{ x \in \X: \norm{\langle (\A^\eta)^* w, x \rangle - \langle w, y^\delta \rangle } \\
    &\leqslant (\delta + \eta \rho) \norm{w} \},
\end{split}
\end{equation}
with $\rho > 0$ and for a direction $w\in \mathcal{Y}$.
In the semi-discretized setting the method then terminates if a discrete version of the discrepancy principle is satisfied, i.e., for all $k,l$
\begin{equation}
    \norm{\A^\eta_{k,l}x_n - y_{k,l}} \leqslant \tau_{k,l} (\delta_{k,l} + \eta_{k,l} \rho ),
\end{equation}
with constants $\tau_{k,l} > 1$ and noise and inexactness levels $\delta_{k,l}$, $\eta_{k,l}$ given for any semi-discrete subproblem.
The full algorithm for the RESESOP method is specified in the Supplementary Material~\ref{app:resesop} (\cref{alg:resesop}).

\subsection{Post-processing via neural networks}
\label{sec:NN_postporcessing}
\noindent Neural networks $F_\theta$ are finding increasing applications in the field of CT \cite{Wang2020dl_tomo_reco, VillarragaGomez2022microscope_ct_dl}. These include both detection and reconstruction tasks \cite{AnayaIsaza2021overview_dl_medical, Rafiei2023cv_industrial_security, Szczykutowicz2022review_ct_reco}. A common approach for the latter is to post-process CT reconstructions of established methods $\mathcal{T}$, such as filtered backprojection, using neural networks \cite{jin_deep_2017, yang_low-dose_2018, yuan_sipid_2018, liu_interpreting_2020,  you_ct_2020}, e.g.,
\begin{equation}
    x_\mathrm{reco} = F_{\theta}\circ \mathcal{T}_{\text{FBP}}(y^\delta)
\end{equation}
where $F_\theta : \mathcal{X} \to \mathcal{X}$.
The idea is to correct artifacts and noise in the image using learned techniques. The post-processing networks are trained in a supervised way with (simulated) ground truth data from idealized scanning conditions, e.g., with a higher dose, more scan angles, or less motion of the scanned object. These scan settings result in superior image quality but may not be feasible in practical use due to time or security constraints.

The post-processing approaches have the advantage of working mainly in the space of images $\mathcal{X}$. Thus, convolutional neural networks, like the U-Net \cite{ronneberger_u-net_2015} architecture are suitable for this application. In addition, both the classic and the post-processed reconstruction are available to the user. This can increase the acceptance of data-driven approaches among users who already have years of experience evaluating images from a particular method. 
Another network type, so called invertible neural networks, also allows for image post-processing and it additionally allows for an uncertainty quantification. Compared to the U-Net architecture, the concept of conditional invertible neural networks, which are considered in this work, is less standard such that we highlight the important differences in architecture design, training, and reconstruction in the following.  
\subsubsection{Conditional invertible neural networks} 
\label{sec:CINN}
For our application in nanoCT, the uncertainties due to the inaccuracy of the forward model must be taken into account. To this end, it is helpful to switch to the statistical view of inverse problems \cite{kaipio_statistical_2005, tarantola_inverse_1982}. Here, instead of creating a single reconstruction, we are interested in recovering the whole conditional distribution $p_{\mathbf{x} \vert \mathbf{y}}(x|y^\delta)$, respectively $p_{\mathbf{x} \vert \mathbf{\Tilde{x}}}(x|\Tilde{x})$ with $\Tilde{x}=\mathcal{T}(y^\delta)$ in the post-processing setting of the present work. The idea is to approximate this distribution via a generative method, called conditional invertible neural network (iNN) \cite{Ardizzone2019cinn, Winkler2019cnf, anantha_padmanabha_solving_2021} which have been applied to CT reconstruction \cite{denker_conditional_2020, denker_conditional_2021} from exact operators so far.

\underline{Architectures:} Conditional invertible network architectures are typically given as mappings $G_\theta: \mathcal{X} \times \mathcal{X} \rightarrow \mathcal{Z}$ with latent space $\mathcal{Z}$. The conditioning is encoded in the second input variable and invertibility is then given with respect to the first variable, i.e., $\Tilde{G}_{\theta,\Tilde{x}}:= G_\theta(\cdot,\Tilde{x}):\mathcal{X} \to \mathcal{Z}$ is an invertible mapping for any $\Tilde{x} \in \mathcal{X}$. 
In the remainder we consider two specific architectures:
\begin{itemize}
    \item CiNN: A typical choice is a multi-scale architecture, based on NICE \cite{dinh_nice_2015} and RealNVP \cite{dinh_density_2017}. We use additive coupling blocks, the learned invertible downsampling proposed in \cite{etmann_iunets_2020}, and a ResNet conditioning network processing the conditioning input before feeding it into the coupling blocks.
    The conditioning ResNet does not need to be invertible.
    \item CiUNet: We also consider a conditional variant of the invertible U-Net proposed in \cite{etmann_iunets_2020}, which has been studied in \cite{denker_conditional_2021}. We again use additive coupling blocks, and the conditioning network is implemented as a (non-invertible) U-Net, which is connected to the coupling blocks of the invertible U-Net at each respective scale.
\end{itemize}
Further details and illustrations of the architectures are provided in the Supplementary Material~\ref{app:network_training}.

\underline{Training:}
The goal of the conditional iNN during training is to minimize the Kullback-Leibler (KL) divergence to the conditional distribution. Equivalently, for any $\Tilde{x}$ this can be done using the negative log-likelihood
\begin{align*}
        &\min_\theta \mathcal{D}_{\text{KL}}[p(x|\Tilde{x})|| q_\theta(x|\Tilde{x})] 
        \Leftrightarrow \min_\theta  - \mathbb{E}_{\mathbf{x} \vert \mathbf{\Tilde{x}}}[\log q_\theta(x|\Tilde{x})] 
\end{align*}
where $q_\theta $ is parameterized via the conditional iNN $G_\theta$ and an assumed distribution $p_{\mathbf{z}}$, which allows for easy sampling, for a variable $z$ defined on $\mathcal{Z}$. Those are linked by the relation $ x|\Tilde{x} = \Tilde{G}_{\theta,\Tilde{x}}^{-1}(z)|\Tilde{x}$, i.e,
\begin{equation}
q_\theta(x | \Tilde{x})= p_{\mathbf{z}}( \Tilde{G}_{\theta,\Tilde{x}}(x)) | \det J_{\Tilde{G}_{\theta,\Tilde{x}}}(x) |
\end{equation}
where $z|\Tilde{x}\sim p_{\mathbf{z}}$ for any $\Tilde{x}$. Joint training for all $\Tilde{x}$ in terms of expectation minimization and applying 
change of variables thus results in the considered training loss 
\begin{align*}
\ell(\theta) &= - \mathbb{E}_{\mathbf{\Tilde{x}}} \mathbb{E}_{\mathbf{x} \vert \mathbf{\Tilde{x}}}[\log q_\theta(x|\Tilde{x})] = - \mathbb{E}_{(\mathbf{x},\mathbf{\Tilde{x}})} [\log q_\theta(x|\Tilde{x})] \\
&= - \mathbb{E}_{(\mathbf{x}, \mathbf{\Tilde{x}})} \left[\log p_\mathbf{z}(\Tilde{G}_{\theta,\Tilde{x}}(x)) + \log | \det J_{\Tilde{G}_{\theta,\Tilde{x}}}(x) |\right] \\
        & \approx  - \frac{1}{N} \sum_{i=1}^N [ \log p_\mathbf{z}(G_\theta(x^{(i)},\Tilde{x}^{(i)})) \\ & + \log | \det \frac{\partial}{\partial x} G_\theta(x^{(i)}, \Tilde{x}^{(i)}) | ]     \end{align*}
for a given data set of tuples $(x^{(i)},\Tilde{x}^{(i)})$ including ground truth and reconstruction from a predefined method $\mathcal{T}$.
In the remainder we distinguish two variants of training for the conditional iNNs.
First, we consider ground truth images $x$ as the input variable or alternatively we consider the residual $\Delta x = \Tilde{x}-x$ as input variable for the invertible path of the network. The latter case is denoted by the suffix \textit{Res} in the present work (i.e., CiNNRes, CiUNetRes) as for given latent variable $z$ the reconstruction can be interpreted as a ResNet with respect to $\Tilde{x}$, i.e.,
\begin{equation}
    x_\mathrm{reco}= \Tilde{x} - \Tilde{G}_{\theta,\Tilde{x}}^{-1}(z).
\end{equation}

\underline{Reconstruction and uncertainty quantification:}
After training the conditional iNNs we can derive reconstructions and standard deviations from samples drawn from $p_\mathbf{z}$, i.e., given $M$ samples $z^{(m)} \in \mathcal{Z}$ and a $\Tilde{x}=\mathcal{T}(y^\delta)$, we obtain the final reconstruction $x_\mathrm{reco}$ via 
\begin{equation}
x_\mathrm{reco} = \frac{1}{M}\sum_{m=1}^M \Tilde{G}_{\theta,\Tilde{x}}^{-1}(z^{(m)}) =: F_\theta(\Tilde{x})
\end{equation}
which defines our learned post-processing routine $F_\theta$. 
In  addition this approach allows for an uncertainty quantification, e.g., in terms of pixel-wise standard deviation, i.e.,
\begin{equation}
\sigma_\mathrm{reco} = \left( \frac{1}{M}\sum_{m=1}^M (\Tilde{G}_{\theta,\Tilde{x}}^{-1}(z^{(m)}) - x_\mathrm{reco})^2 \right)^{\frac{1}{2}}.
\end{equation}
The Res case (here, $\Tilde{G}_{\theta,\Tilde{x}}^{-1}(z^{(m)})=\Delta x^{(k)}$) works analogously with minor adaptation, i.e.,
\begin{align}
x_\mathrm{reco} &= \tilde{x} - \frac{1}{M}\sum_{m=1}^M \Tilde{G}_{\theta,\Tilde{x}}^{-1}(z^{(m)}) =: F_\theta(\Tilde{x})\\
\sigma_\mathrm{reco} &= \left( \frac{1}{M}\sum_{m=1}^M (\Tilde{G}_{\theta,\Tilde{x}}^{-1}(z^{(m)}) - (\Tilde{x}-x_\mathrm{reco}))^2 \right)^{\frac{1}{2}}.
\end{align}

\subsection{Data set generation}
\label{sec:dataset}
\noindent In order to train the methods described in \cref{MaterialsMethods}, it is necessary to simulate tomographic imaging data including inexactness in the operator.
The generated dataset is composed of 32\,095 samples split into 30\,490 (95\%) for training, 1\,284 (4\%) for validation and 321 (1\%) for testing.

\begin{figure}
\centering
     \vspace{0.2em}
     \includegraphics[width=\columnwidth]{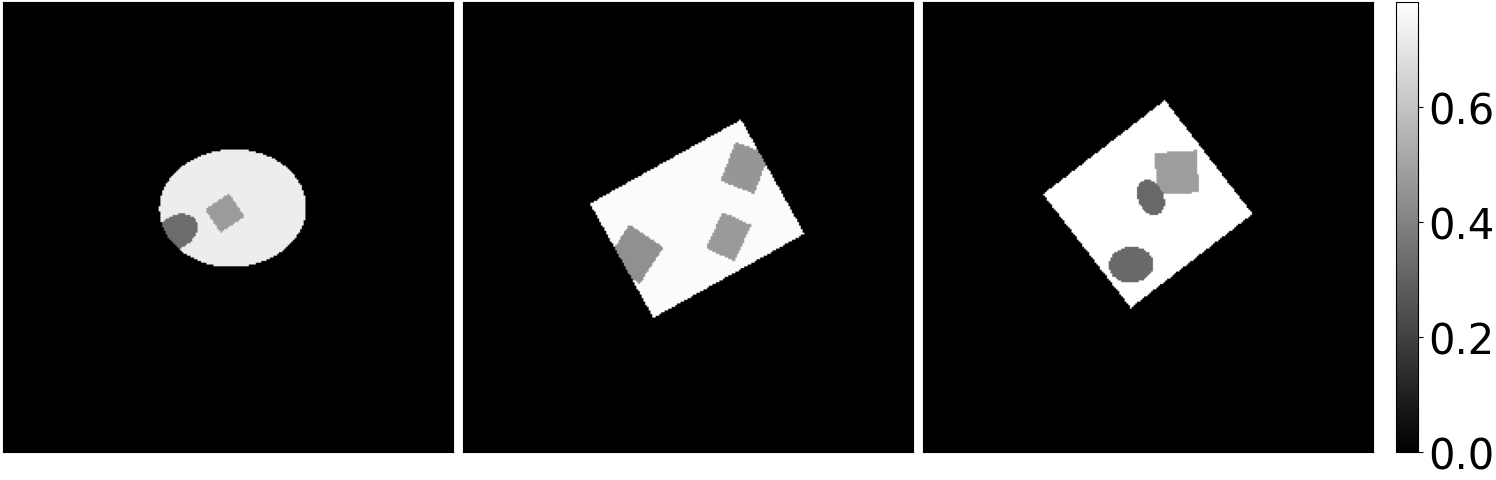}
      \caption{Three sample phantoms with $255\times 255$ pixels included in the generated test data set.}
      \label{phantoms}
\end{figure}

Phantoms $x^{(i)}$ were generated for a field of view of 255 $\times$ 255 pixels corresponding to spatial positions $r\in \R^2$ and are constructed from randomly generated rectangles and ellipses which overlay each other and have different levels of density ranging from 0 to 1, 0 being no density. One phantom has a main shape and up to 3 subshapes within (see \cref{phantoms}), which is inspired by common experimental settings in nanoCT \cite{Dremel_2018, Blaek2019, Cressa2022fib_sem_correlative}.

\begin{table}
\centering
\caption{Parallel and fan beam geometry. $S$ denotes the number of reconstruction pixels, $K$ defines the number of scanning angles (time steps), $L$ defines the number of detector pixels, and $D$ specifies the source radius in reconstruction pixel size unit. The detector extent is such that the projection beams cover the whole image.}
\include{tables/geometry.tex}
\label{geometry_tab}
\end{table}

\begin{figure}
\centering

      \includegraphics[width=\columnwidth]{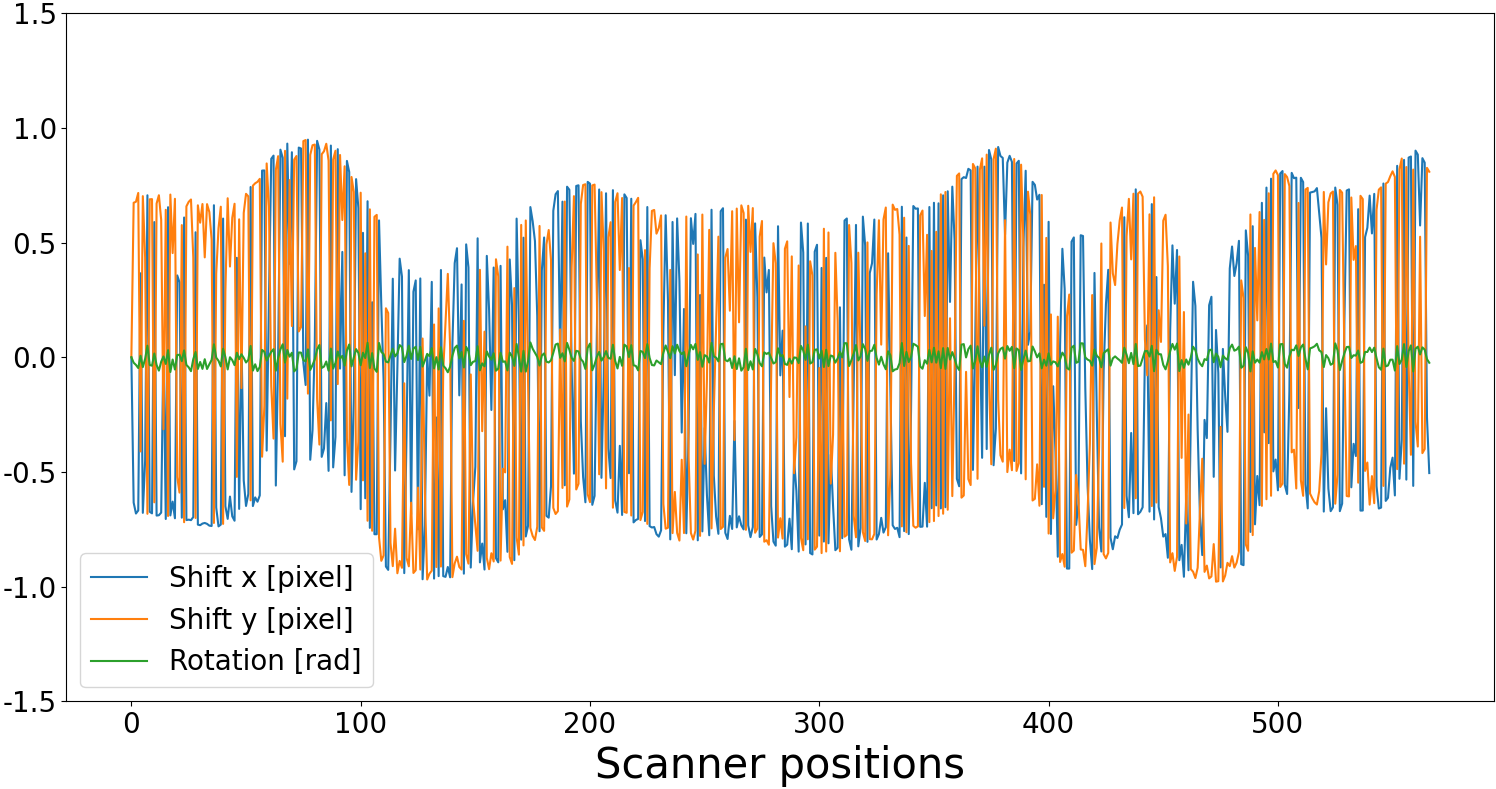}
      \caption{Randomly generated perturbations for the data set generation with respect to scanner positions $\phi$.}
      \label{perturbances}
\end{figure}

The objective in the considered problem setting is to simulate distortions in the measurements caused by small vibrations near the CT scanner during the scanning process. As a result, this is incorporated as movements of the entire phantom with respect to the scanner position/scanning angle $\phi$, where an example of randomly generated movements is illustrated in \cref{perturbances}. The vibrations have been mimicked by overlapping dampened sinus waves (38 for parallel and 9 for fan beam) with different starting points during the scanning process $\epsilon(\phi)$. The starting points and parameters for the sinus waves are randomly generated for each sample, but limited to a maximum $r_1$- and $r_2$-direction shift distance from the center. 

In addition to the sinus noise $\epsilon(\phi)$, which is a deterministic function with respect to $\phi$ randomly differing for each phantom $x$, there is a small random noise $\xi_{r-\mathrm{dir}}:=(\xi_{r_1-\mathrm{dir}}$, $\xi_{r_2-\mathrm{dir}})^t$ and $\xi_{\phi-\mathrm{dir}}$ taken normal distribution with randomly drawn standard deviation for each angle and applied to $r_1$-, $r_2$-direction shifts and rotations, each drawn for any angle $\phi$.
Here, the standard deviations are drawn from $\mathcal{N}(0.127,0.0254^2)$.
In summary, this can be formalized as considering the randomly moved phantom $x_\mathrm{non-perturbed}(r)$ for the particular scan angle $\phi$ as ($r\in \R^2$)
\begin{align}
&x_\mathrm{perturbed}(\phi)(r) =x_\mathrm{non-perturbed}(d(\phi,r)) \notag \\
&\text{with distortion} \notag \\
    &d(\phi,r)=R( \xi_{\phi-\mathrm{dir}}(\phi)) r + \xi_{r-\mathrm{dir}}(\phi)+ \epsilon(\phi)
\end{align}
where $R$ is the corresponding rotation matrix in $\R^2$.

\begin{figure}
\centering

      \includegraphics[width=\columnwidth]{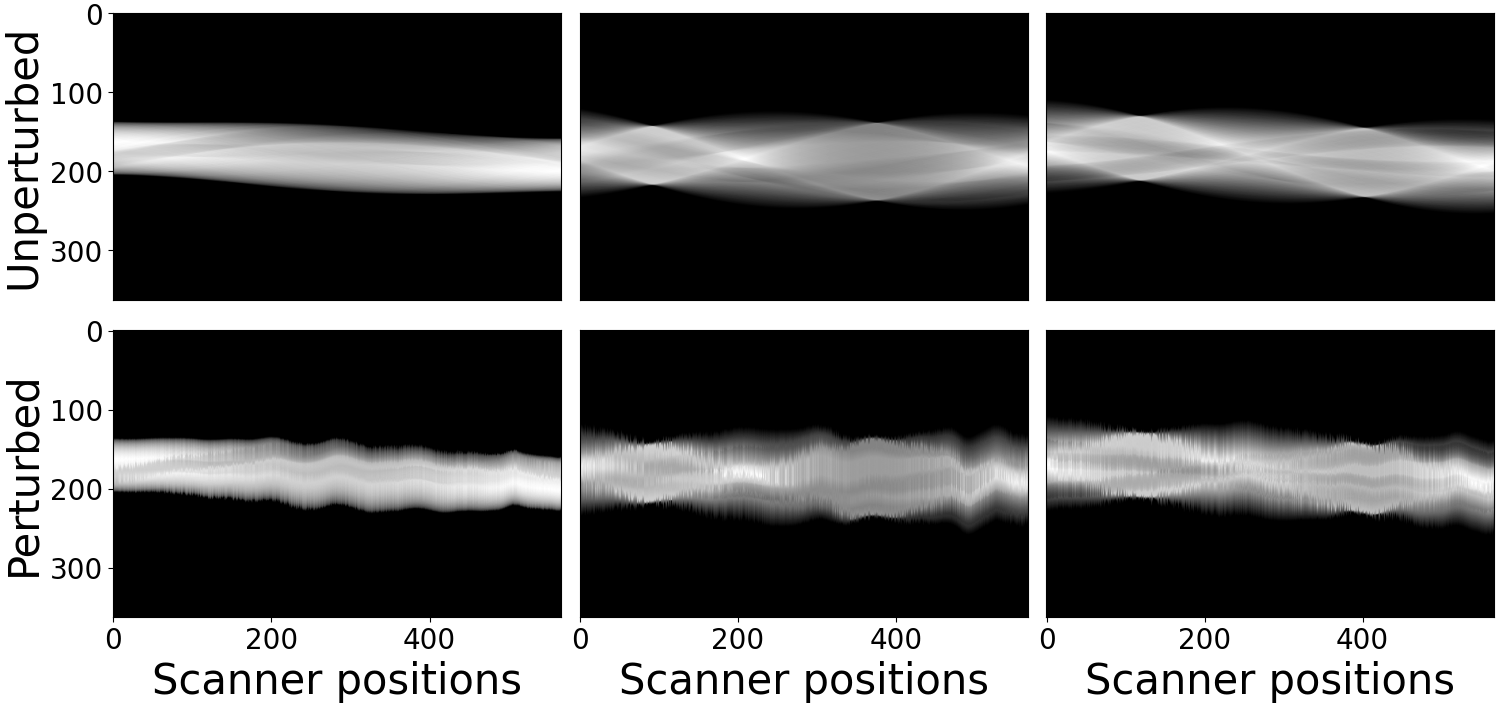}
      \caption{Non-perturbed (top) and perturbed (bottom) sinograms of the phantoms in \cref{phantoms} resulting from the perturbations like illustrated in \cref{perturbances}. Scanner positions $\phi$ are on the horizontal axis and detector positions are on the vertical axis.}
      \label{sinograms}
\end{figure}

The generated phantoms and perturbations were then used to calculate the sinograms $y^\eta$ and $y^\delta$ for the non-perturbed and the perturbed phantoms (see \cref{sinograms}) for a parallel- and a fan-beam geometry.
The sinograms were created using Operator Discretization Library (ODL) \cite{adler_operator_2018} and ASTRA Toolbox \cite{aarle_astra_2015} and the respective geometry parameters from \cref{geometry_tab}.

Note that the perturbed phantom $x_\mathrm{perturbed}(\phi)$, which is represented as a sequence of images with respect to $\phi$, is assigned to the exact but unknown operator $\mathcal{A}$ case and that the non-perturbed phantom $x_\mathrm{non-perturbed}$ is assigned to the inexact but known operator $\mathcal{A}^\eta$ case. Here we use $\mathcal{A}^\eta = \mathcal{A}^{\dagger}$ which is chosen as the operator without phantom perturbation in any case, i.e., measurement from exact unknown operator are $y$ with $y(\phi,s)=(\mathcal{A}^{\dagger}x_\mathrm{perturbed}(\phi))(\phi,s)$ ($s$ being the detector position) and hypothetical measurement from inexact known operator $y^{\eta}=\mathcal{A}^{\dagger}x_\mathrm{non-perturbed}$. 
In the present work emphasis is given on the operator inexactness such that the additive noise $\xi=0$ is chosen in any case, i.e., $y^\delta=y$ and $\delta=0$.
The data set includes the tuples $(x^{(i)},y^{(i)})$ where $x^{(i)}=x_\mathrm{non-perturbed}^{(i)}$ and $y^{(i)}=y^{\delta,(i)}$. In addition $y^{\eta,(i)}$ is also available in the data set.
The dataset is made available under \url{https://doi.org/10.5281/zenodo.8123498}.

\section{Numerical experiments and results}
\label{sec:experiments}

\begin{table}
\caption{Quantitative evaluation on the parallel beam test data for measurements from unperturbed phantoms.}
\label{tab:parallel_unperturbed_quantitative}
\include{tables/unperturbed_parallel_beam_test_data.tex}
\end{table}

\noindent We conducted numerical experiments on the simulated phantom data set described in \cref{sec:dataset}. 
For this we proceed as follows:

 First, we evaluate each non-trained reconstruction method $\mathcal{T}$ on the respective data sets $\{(x^{(i)},y^{(i)})\}_i$ for parallel and fan-beam geometry to obtain the reconstructions $\Tilde{x}_\mathcal{T}^{(i)}$ from reconstruction method $\mathcal{T}$ using $\mathcal{A}^\eta = \mathcal{A}^\dagger$ with $\mathcal{A}^\dagger$ being either a standard parallel or fan beam setting. Here, for $\mathcal{T}$ we used FBP, Dremel (see \cref{sec:dremel}), and RESESOP (see \cref{sec:resesop}). For the RESESOP method an estimate of $\eta$ needs to be obtained prior the reconstruction. Here, we assume that a conservative estimate of $\eta_{k,l}$ for any semi-discrete subproblem is given in terms of the maximum deviation over all detector positions for each angle, i.e, $\eta_{k,l} = \max_{\ell \in \mathcal{L}} (\|y_{k,\ell} - y^\eta_{k,\ell}\| )$ for any $l \in \mathcal{L}$. In the numerical experiments we also include the cases of $\pm 20 \%$ over- and underestimation.
 The Dremel method includes some algorithmic options, which were selected by evaluation on validation samples (see Supplementary Materials~\ref{app:dremel}).

Second, we train a learned post-processing routine $F_\theta$ using the data tuples $(x^{(i)}, \Tilde{x}^{(i)}_\mathcal{T})$. For the UNet architecture $F_\theta$, the loss is given by the Euclidean norm, i.e., minimize  $\sum_i \|F_\theta(\Tilde{x}_\mathcal{T}^{(i)}) - x^{(i)}\|^2$ with respect to network parameters $\theta$. The invertible network architectures are trained according to the descriptions in \cref{sec:CINN} using a standard Gaussian density $p_\mathbf{z}$. Here, we use a CiNN and a CiUNet as outlined in \cref{sec:NN_postporcessing}. All trainings are performed using the Adam optimizer, and the checkpoint with minimum validation loss is selected.

Third (for iNNs only), the iNNs serve as an uncertainty estimator in terms of an estimate for pixel-wise standard deviation.   

Finally, the resulting reconstruction schemes are then evaluated quantitatively (PSNR, SSIM) and qualitatively on the test set of the respective data set for parallel- and fan-beam geometry. 

Further algorithmic details are provided in the Supplementary Materials~\ref{app:dremel}, \ref{app:resesop}, and~\ref{app:network_training}. The code for the learned post-processing schemes via UNet and conditional iNNs is made available under \url{https://gitlab.informatik.uni-bremen.de/inn4ip/cond-inn4nanoct}.

\subsection{Architecture specifications}
\noindent The UNet architecture used in the present work is an adapted standard U-Net architecture with 5 scales and 32 to 128 channels using strided convolution downsampling, bilinear upsampling, skip connections and a sigmoid output activation. A detailed illustration of the used UNet can be found in the Supplementary Material~\ref{app:network_training} (\cref{fig:unet_scheme}).

The CiNN architecture has 6 scales and uses learned invertible downsampling, splitting, additive conditional coupling blocks, and a ResNet conditioning network.
In each downsampling the number of channels naturally is multiplied by 4, after which half of the channels are split off to form a part of the latent output $z$ (except after the first downsampling).
The other half of the channels is further processed by a conditional coupling block, followed by the next downsampling, until the coarsest scale with 128 channels is reached.
A small part of this signal is finally processed by random permutation and a conditional coupling block, before concatenating all parts of the latent output $z$.
See the Supplementary Material~\ref{app:network_training} (\cref{fig:cinn_scheme}) for an illustration.

The CiUNet architecture is based on the iUNet architecture \cite{etmann_iunets_2020} using 5 scales, learned invertible down- and upsampling, additive conditional coupling blocks, and skip connections forwarding half of the channels at each scale from the encoder to the decoder.
A non-invertible U-Net is used for the conditioning, connected to the conditional coupling blocks at the respective scales in reversed order, i.e., U-Net's decoder activations are connected to CiUNet's encoder coupling blocks.
After the downsampling in the encoder to the coarsest resolution, at which 32 channels are used, the features are upsampled back to the original resolution by the decoder, while concatenating with the forwarded channels from the skip connections.
The decoder output then forms the latent output $z$ (via reshaping).
A sequence of four additional layers, namely activation normalization, downsampling, conditional coupling and upsampling, is inserted at the beginning before the encoder.
An illustration and further details can be found in the Supplementary Material~\ref{app:network_training} (\cref{fig:iunet_scheme}).

Reconstructions and standard deviations are computed with $M=100$ samples from the latent space.

\subsection{Non-perturbed phantom case}
\noindent As a sanity check and also as a reference, we also trained and evaluated selected methods on unperturbed measurements, i.e., those measurements generated by unperturbed phantoms and $\mathcal{A}^\dagger=\mathcal{A}^\eta$ which is illustrated in \cref{tab:parallel_unperturbed_quantitative}.
All reconstruction methods $\mathcal{T}$ without post-processing deliver accurate reconstructions as expected. 
The post-processing via UNet and the residual iNNs (CiNNRes, CiUNetRes) result in the largest performance improvements in terms of PSNR as well as in SSIM. Here, the UNet provides the largest improvement in terms of PSNR. 
The post-processing via non-residual invertible networks results in at least similar or slightly improved performance in PSNR but the SSIM is decreased compared to the respective non-post-processed case. 
The image reconstructions of these methods are illustrated in Supplementary Material~\ref{app:figs_non_perturbed_parallel} (\cref{fig:reco_unperturbed_parallel_std_algo,fig:reco_unperturbed_parallel_U-Net,fig:reco_unperturbed_parallel_CINN,fig:reco_unperturbed_parallel_CI-U-Net}) for the phantoms in \cref{phantoms}.

\subsection{Perturbed phantom case }

\begin{table}
\caption{Quantitative evaluation on the parallel beam test data for measurements from perturbed phantoms. Overall best performance is highlighted in \textbf{bold} for any $\mathcal{T}$, best performance under the iNN post-processing networks is highlighted in \textcolor{teal}{color}, and overall best performance is \underline{underlined}.}
\include{tables/perturbed_parallel_beam_test_data.tex}
\label{tab:parallel_perturbed_quantitative}
\end{table}

\begin{table}
\caption{Quantitative evaluation on the fan beam test data for measurements from perturbed phantoms. Overall best performance is highlighted in \textbf{bold} for any $\mathcal{T}$, best performance under the iNN post-processing networks is highlighted in \textcolor{teal}{color}, and overall best performance is \underline{underlined}.}
\include{tables/perturbed_fan_beam_test_data.tex}
\label{tab:fan_perturbed_quantitative}
\end{table}

\noindent \underline{Quantitative results:} The quantitative results on the perturbed phantom test data set for parallel beam are summarized in \cref{tab:parallel_perturbed_quantitative}. 
First, we can immediately observe that the non-learned methods Dremel and RESESOP perform superior to FBP in terms of PSNR as well as in SSIM which is due to the design of those sequential reconstruction methods taking into account the inexactness.
Second, in the learned post-processing schemes the UNet approach delivers the largest improvement when compared to all other conditional iNN approaches. The largest improvement can be observed for RESESOP when combined with the UNet. Large improvements can even be observed for over and underestimated $\eta$. But also FBP and Dremel can highly benefit from the learned UNet post-processing particularly in terms of SSIM. 
Third, the conditional iNNs show a more diverse behavior. Dremel and RESESOP perform better when combined with the residual approaches CiNNRes and CiUNetRes. Interestingly, the performance becomes worse for the non-residual approaches CiNN and CiUNet (except for RESESOP and CiNN with respect to PSNR), and the non-residual CiUNet for RESESOP unexpectedly performs worse with the correct per-angle maximum error level $\eta$ than with over- and under-estimated $\eta$, which might be a rather random effect caused by apparent training difficulties. For FBP the situation is slightly different. Here, all iNNs are able to improve performance in terms of SSIM but for PSNR this holds true for CiNN and CiNNRes only. In summary, the UNet post-processing and in particular when combined with RESESOP provides the best performance and among all conditional iNNs the CiUNetRes provides the largest improvements.
The learned post-processing approaches are also evaluated on the fan beam data set for Dremel and RESESOP resulting in \cref{tab:fan_perturbed_quantitative}. 
Here, similar observations as in the parallel beam case can be made with minor differences in the PSNR results.
Finally, we investigated the robustness of the trained post-processing considering different reconstruction methods $\mathcal{T}$ for training and testing. The results for the UNet are illustrated in \cref{tab:parallel_perturbed_quantitative_mixing}. Networks trained on Dremel and FBP show similar or slightly worse behavior for reconstructions from different algorithms. In contrast, architectures trained with RESESOP reconstructions perform worse with other reconstruction algorithms, particularly in terms of PSNR. 
The RESESOP reconstructions might have an inherent characteristic which is advantageous for the UNet architecture.
Further quantitative results are provided in the Supplementary Materials~\ref{app:quant_perturbed_parallel}. 

\begin{figure}
    \centering
         Reconstructions \\[2mm]
    \includegraphics[width=\columnwidth]{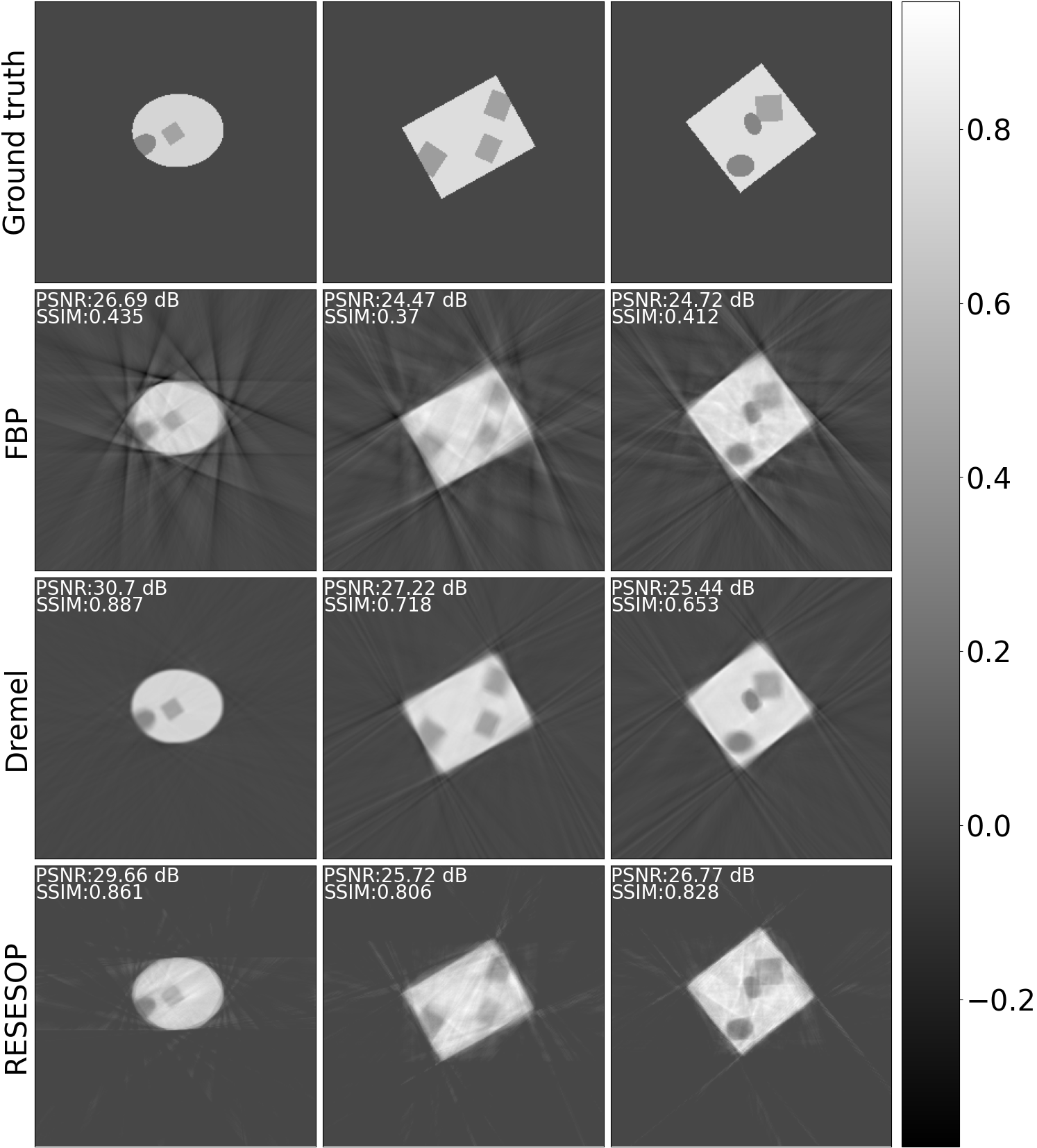}
    \caption{Image reconstructions using the non-trained reconstruction methods $\mathcal{T}$ on perturbed parallel beam data for the phantoms in \cref{phantoms}.}
    \label{fig:reco_perturbed}
\end{figure}

\begin{figure}[t]
    \centering
    Reconstructions \\[2mm]
    \includegraphics[width=\columnwidth]{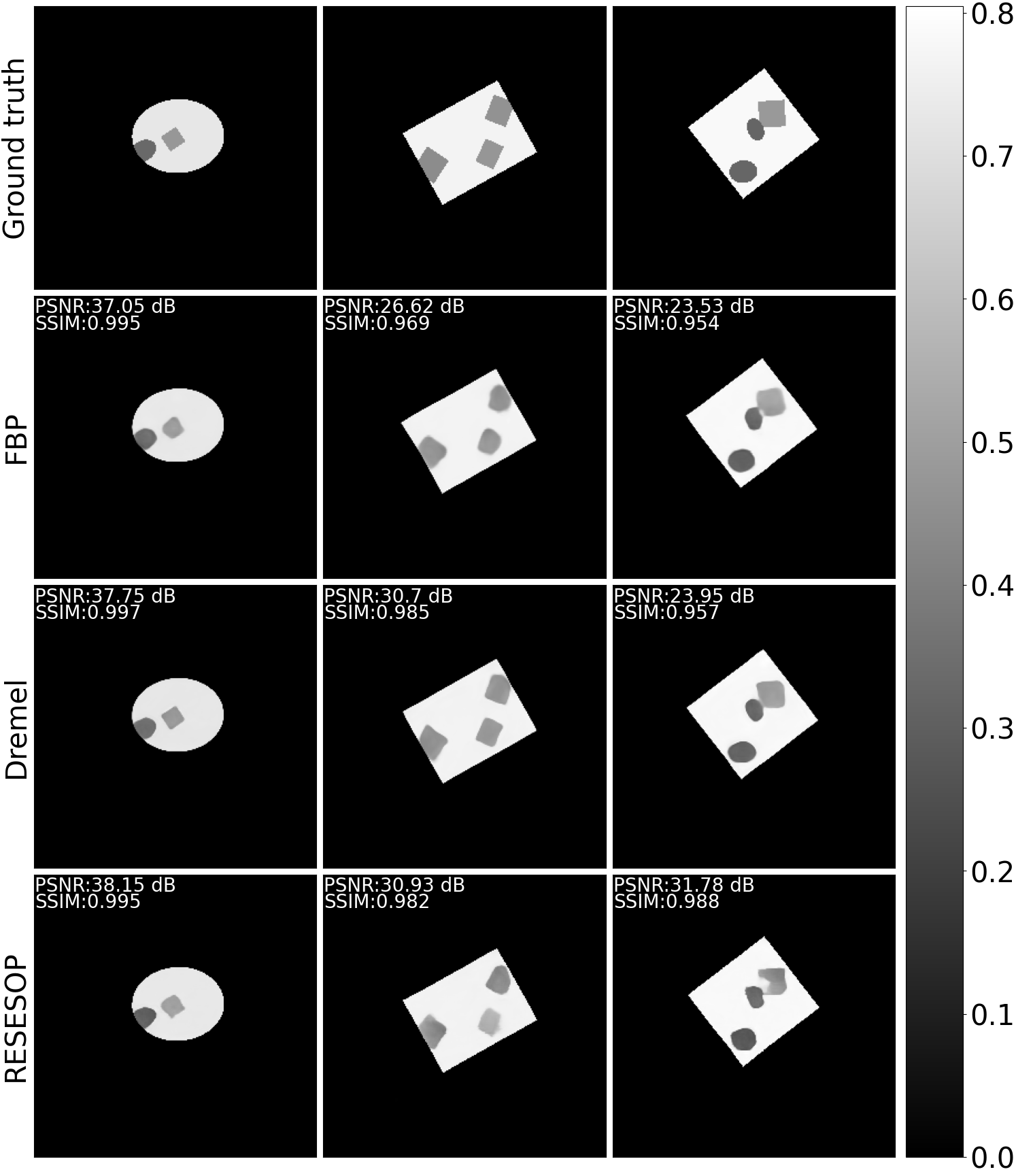}\\
    Differences\\[2mm] 
        \includegraphics[width=\columnwidth]{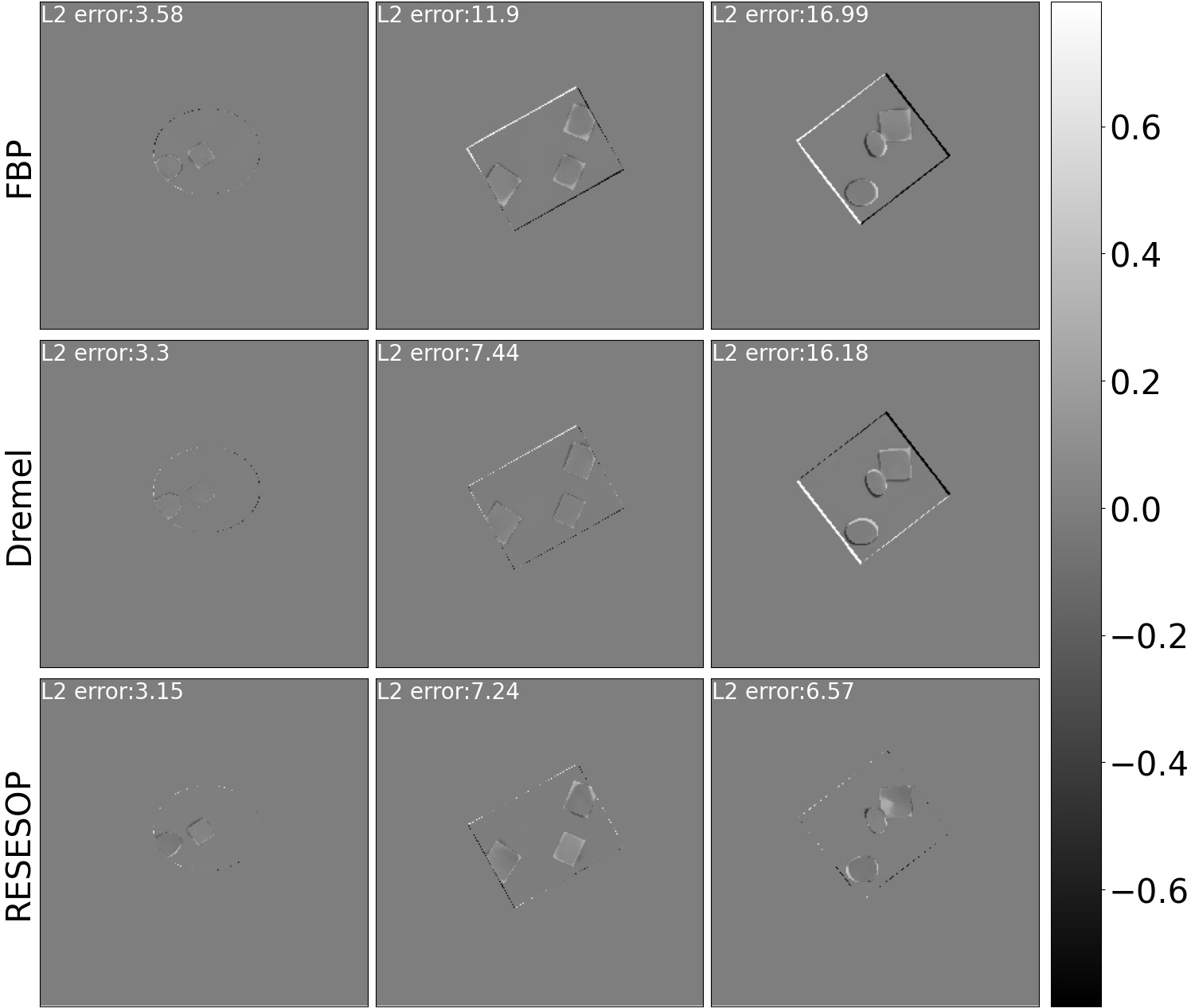}
    \caption{Image reconstructions $x_\mathrm{reco}$ and differences to ground truth using the non-trained reconstruction methods $\mathcal{T}$ combined with \textbf{UNet} post-processing on perturbed parallel beam data for the phantoms in \cref{phantoms}.}
    \label{fig:unet_perturbed}
\end{figure}

\underline{Qualitative results:} Qualitatively, we can also observe different characteristics in the remaining artifacts in the reconstruction. 
The outcome of the non-learned reconstruction methods $\mathcal{T}$ which also serve as the input to the learned post-processing schemes is illustrated in \cref{fig:reco_perturbed}. 
All methods suffer from streaking artifacts emerging mainly on edges and corners in the phantom. 
The intensity of distortions decreases from FBP over Dremel to RESESOP. 
For FBP and RESESOP we can also observe more severe distortions within the objects.
The post-processing via the UNet is illustrated in \cref{fig:unet_perturbed}. %
Qualitatively, the reconstructions in \cref{fig:unet_perturbed} are superior to those provided solely by the methods $\mathcal{T}$. In all cases the streaking artifacts are removed, also in the objects interior for FBP and RESESOP.
Qualitative differences which cause quantitative differences become apparent in the difference images in \cref{fig:unet_perturbed}.
FBP+UNet and Dremel+UNet have more severe deviations particularly at phantom edges. For RESESOP these distortions are less distinct but in some cases other artifacts appear in the reconstruction, e.g., in the reconstruction at the bottom right in \cref{fig:unet_perturbed}.

\begin{figure}[t]
    \centering
     Reconstructions \\[2mm]
    \includegraphics[width=\columnwidth]{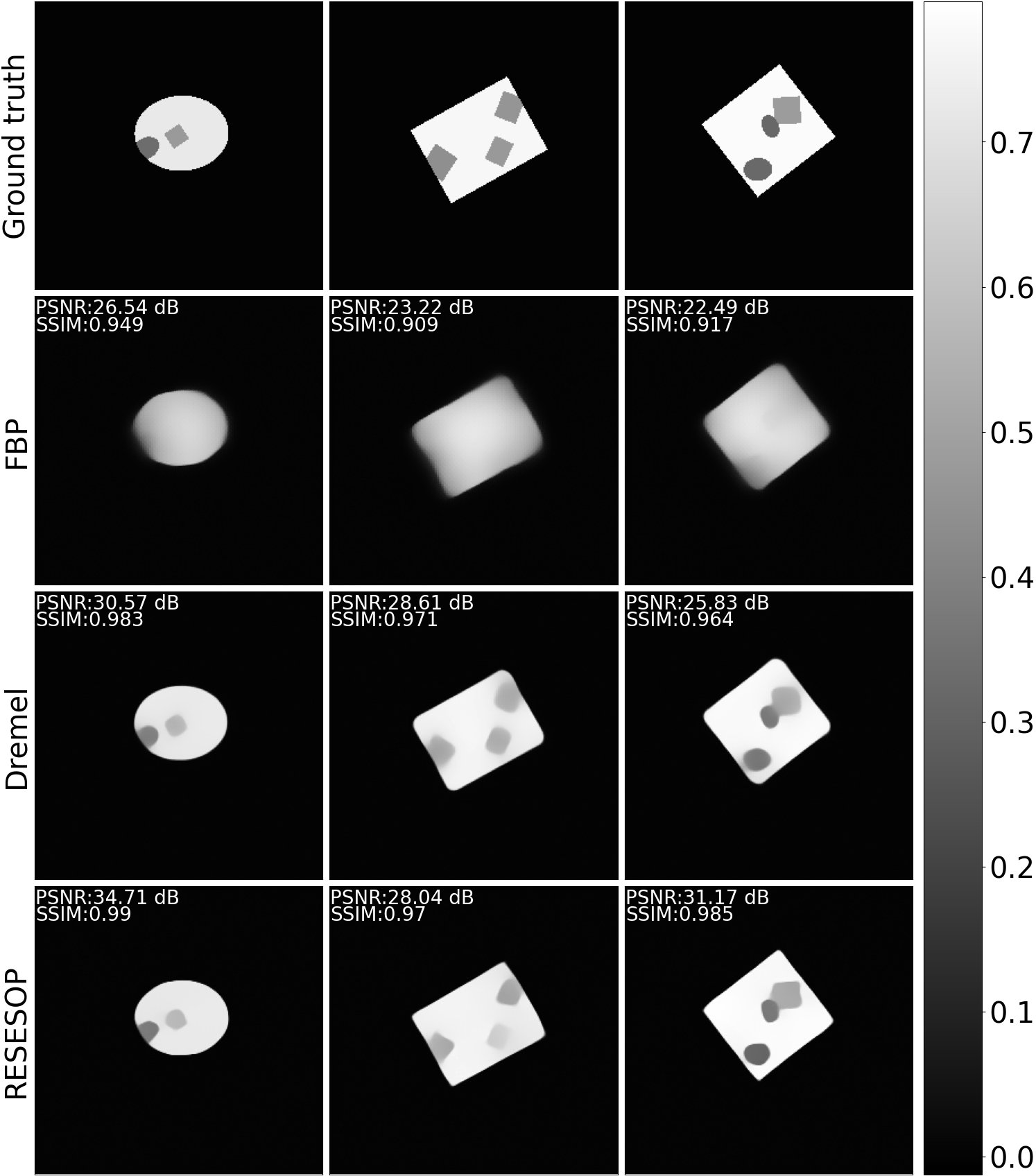}\\
    Differences\\[2mm]
        \includegraphics[width=\columnwidth]{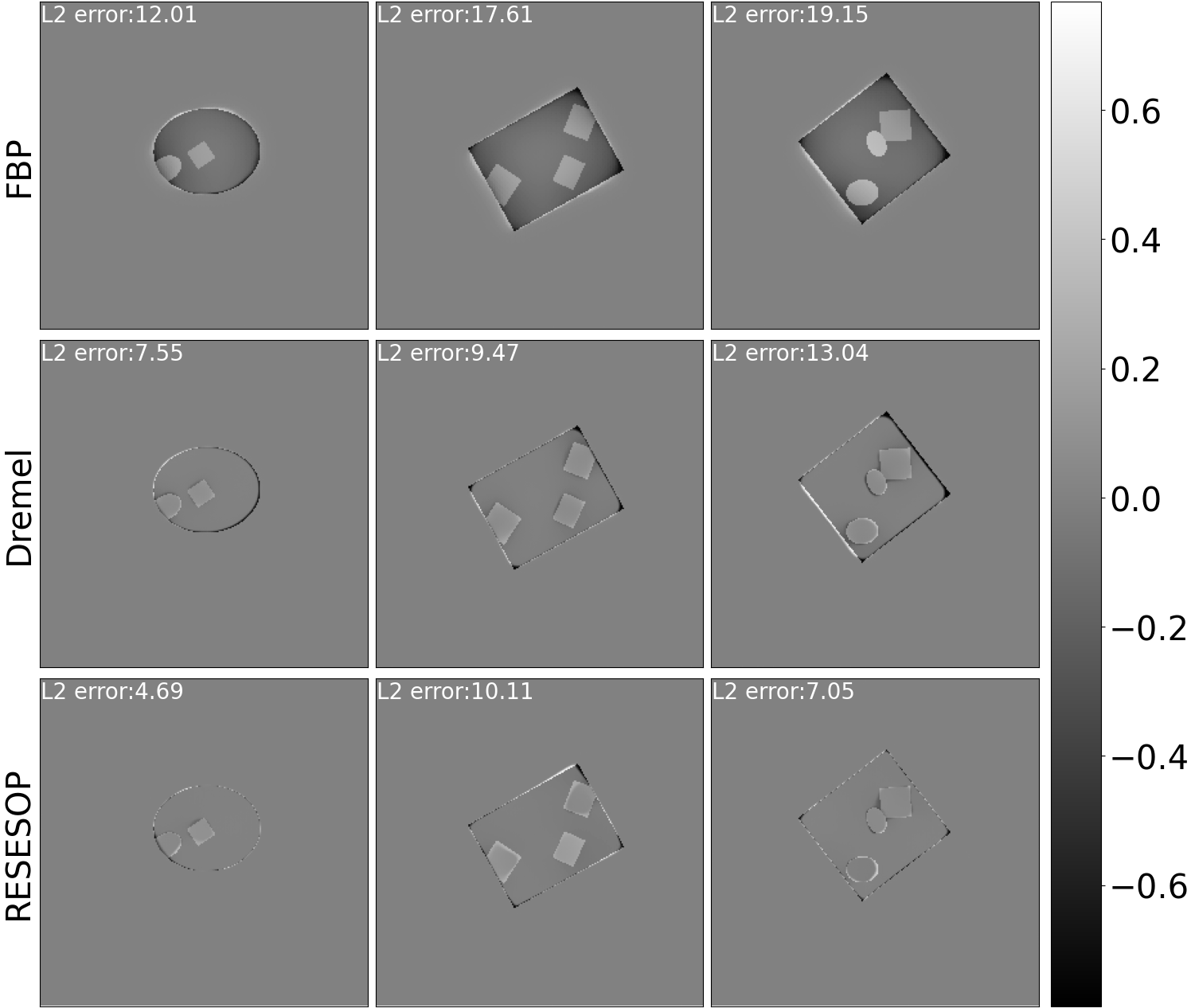}
    \caption{Image reconstructions $x_\mathrm{reco}$ differences to ground truth using the non-trained reconstruction methods $\mathcal{T}$ combined with \textbf{CiUNetRes} post-processing on perturbed parallel beam data for the phantoms in \cref{phantoms}.}
    \label{fig:iunet_perturbed}
\end{figure}

\underline{Qualitative results cond. iNNs:} The CiUNetRes showed the best performance among the iNNs such that we restrict the presented qualitative results to this particular network case. 
Further image reconstructions of other methods and also for the fan beam setting are provided in the Supplementary Materials~\ref{app:figs_perturbed_parallel} and~\ref{app:figs_perturbed_fan}.
The post-processing via the CiUNetRes is illustrated in \cref{fig:iunet_perturbed}. %
Here, we can observe that the CiUNetRes tends to provide smoother reconstructions. In the FBP case we observe severe smoothing which is most likely the reason for the quantitative performance drop in PSNR. For Dremel and RESESOP we observe much less smoothing and in the RESESOP case the artifact in the interior of the right phantom disappeared. 
The extent of smoothing becomes more apparent in the difference images in \cref{fig:iunet_perturbed}. 
Here, the rectangular shapes show more severe differences in the corners where these effects decrease from FBP to RESESOP again. 
In any $\mathcal{T}$ case the CiUNetRes is able to remove the streaking artifacts but the more severe these artifacts have been the more smoothing is imposed which results in oversmoothed image reconstructions with varying intensity. Phantoms with smooth contours are likely to be advantageous for the CiUNetRes due to the observed effect in the corners. In sum, the CiUNetRes is able improve the reconstructions from Dremel and RESESOP but the quantitative improvement is smaller when compared with the UNet approach which is likely to be caused by the imposed invertibility in the network and the sampling mean which can also contribute to the observed smoothing effect.

\begin{figure}
    \centering
    \vspace{0.2em}
    \includegraphics[width=\linewidth]{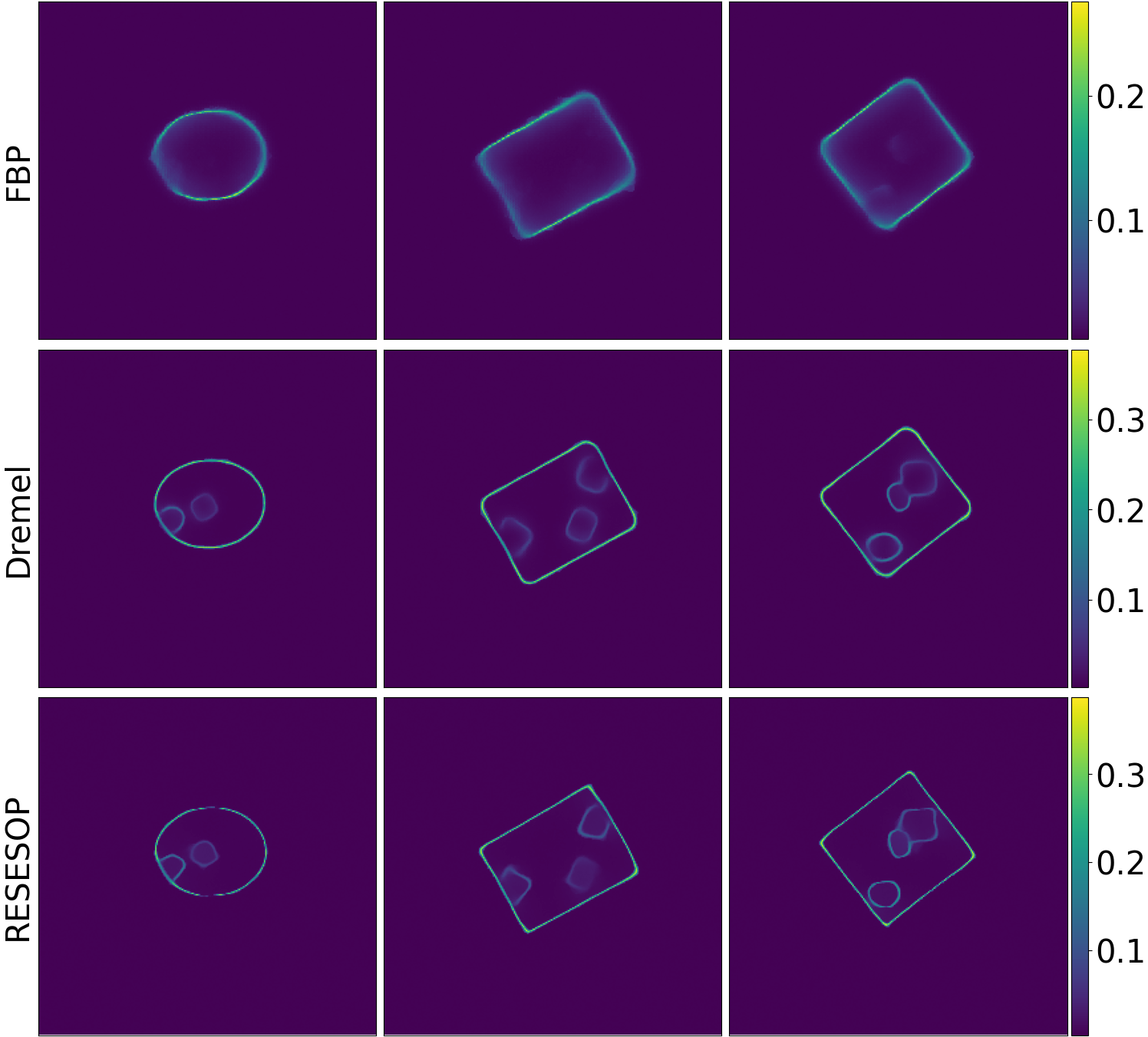}
    \caption{Standard deviation $\sigma_\mathrm{reco}$ according to \cref{sec:CINN} for the \textbf{CiUNetRes} reconstructions in \cref{fig:iunet_perturbed}.}
    \label{fig:std_iunet_perturbed}
\end{figure}

\underline{Uncertainty estimation cond. iNNs:}
But the invertibility and the sampling also opens the door for extracting additional features. For this we illustrate the pixel-wise standard deviation in \cref{fig:std_iunet_perturbed}. In the FBP reconstruction the outer contour is only extracted which gives limited information about the actual uncertainty, e.g., when compared with the differences to ground truth illustrated in \cref{fig:iunet_perturbed}. In case of Dremel and RESESOP we observe a larger similarity to the actual differences to ground truth. The standard deviations indicate an uncertainty in the reconstructed location of the edges. But similar to the actual reconstruction the uncertainty in non-smooth contours like the edges is not properly predicted. But on the other side the standard deviation can also indicate differences in the amount of uncertainty, e.g., being apparent in the middle phantom for RESESOP where the left and right inclusions in the bigger rectangle are more severely mis-located than the middle one.
This can also be observed in the standard deviation. On the circular phantom at the left the standard deviation also predicts less uncertainty in the horizontal edges compared to the vertical edges which is in line with the difference to ground truth.

\section{Discussion and Conclusion}
\label{sec:discussion}

\begin{table*}[b]
\caption{Quantitative evaluation (PSNR (top) and SSIM (bottom)) for mixed training and testing inputs to the post-processing \textbf{UNet} evaluated on the perturbed parallel beam test data for measurements from perturbed phantoms. The rows encode the reconstruction method used for training the columns encode the reconstruction method used for quantitative evaluation on the test set.}
\include{tables/input_mixing_psnr_perturbed_parallel_beam_test_data_unet.tex}
\include{tables/input_mixing_ssim_perturbed_parallel_beam_test_data_unet.tex}
\label{tab:parallel_perturbed_quantitative_mixing}
\end{table*}

In the present work we investigated different classical and learning-based computational strategies how to deal with inexactness in the underlying forward operator motivated by nanoCT applications where the inexactness has a severe probabilistic characteristic due to environmental influences. 
On the methodological side sequential algorithms already taking into account that the forward operator is inexact were combined with learning-based components such as classic U-Nets and more novel network architectures such as conditional invertible neural networks. The latter class also allowed for additional feature extraction for the purpose of uncertainty quantification.
Overall, the non-learned reconstructions can be further improved by a learned post-processing network in most cases. In the proposed computational pipeline the performance of the learned component relies on the quality of the preliminary reconstruction, i.e, the post-processing method only has access to the information from the output of the previous process and not from the actual measurement. Consequently, an initial reconstruction with too little or too distorted information can only be improved to a limited extent. Future works include an extension taking into account a conditioning depending on both actual measurement and preliminary reconstruction.  

One particular focus of the work was the investigation of conditional iNN approaches as learned post-processing tools for reconstructions from inexact operators as they also provide the opportunity to extract features to quantify the uncertainty. Here, the kind of training has a larger impact on the result, i.e., training with respect to the residual is superior. This might be related to the complexity of the learned feature distributions. In the residual case the error pattern of the preliminary reconstruction needs to be learned solely, while in the full reconstruction case the error pattern as well as the image class needs to be encoded in the network somehow. 

The overall reconstruction performance in particular for the CiUNetRes has been only slightly worse compared to the classic U-Net approach. 
We observed that this might be related to an inherent smoothing effect in the iNNs' reconstructions along contours.  Here, one needs to take into account that the underlying architectures have a different structure in particular within the invertible path of the iNN. 
The particular choice of invertible network architecture is likely to be one of the origins of the observed smoothing effect.
In the present work we use coupling layers which tend to result in universal approximators of diffeomorphisms \cite{teshima2020coupling, ishikawa2022universal} which are compared to classic U-Nets being composed of solely continuous and not necessarily differentiable layers. 
Here, future works may relax this property in the invertible path of the post-processing network using invertible residual networks \cite{behrmann_invertible_2019} which tend to approximate homeomorphisms only \cite{zhang2020approximation}.
In summary, the conditional invertible networks for post-processing can provide a powerful tool to further improve image reconstructions from inexact forward operators and to simultaneously provide insights into the remaining operator uncertainty.

Besides the previously mentioned possible extensions, future research directions include the investigation of estimated higher order statistics, the application to measured nanoCT data, and also the transfer to other imaging modalities suffering from operator inexactness, e.g., such as magnetic particle imaging \cite{Gleich2005}.

\section*{Acknowledgments}
J.L.\ and M.S.\ were funded by the German Research Foundation (DFG; GRK 2224/1). J.L. was additionally supported from the DELETO project funded by the German Federal Ministry of Education and Research (Bundesministerium für Bildung und Forschung, BMBF, project number 05M20LBB). T.K. acknowledges financial support from the KIWi project funded by the German Federal Ministry of Education and Research (BMBF, project number 05M22LBA). A.O.~and A.W.~acknowledge funding by the German Federal Ministry of Education and Research (BMBF) under 05M20TSA (DELETO) and by Hermann und Dr. Charlotte Deutsch–Stiftung, A.W.~additionally acknowledges support by the German Research Foundation (DFG; Project-ID 432680300 -- SFB 1456).

\bibliographystyle{IEEEtran}
\bibliography{references}

\vfill

\clearpage
{\appendices
\newpage 
\onecolumn
\section{Pseudocode Dremel}
\label{app:dremel}
\RestyleAlgo{ruled}

\begin{algorithm}[h]
\setstretch{1.2}
\SetKw{KwInit}{Init:}
\caption{Dremel Pseudocode (based on \cite{Dremel_2018})}
\label{alg:dremel}
\SetKwInOut{MyInput}{Input}
\MyInput{Models for each scanner angle with no shifts: $\mathcal{A}_{[1]}, \dots, \mathcal{A}_{[K]}$\\
Perturbed sinogram: $y^\delta$}
\KwResult{$x$}
\SetKwInOut{MyInit}{Init}
\MyInit{${\rm scannerAngles} \gets 1 \ldots K$\\
 ${\rm detectorPixels} \gets 1 \ldots 2 \cdot P + 1$\\
 $x \gets \textrm{empty phantom}$\\
 $s_{[1]}, \dots, s_{[K]} \gets 0, \dots, 0$}
\While{$\neg ${\rm stoppingCriterion}}{
    \For{angle index $k$ in {\rm scannerAngles}}{
        $y_\mathrm{FP} \gets \mathcal{A}_{[k]}\,x$ \Comment{compute forward projection}\\
        $x \gets x - \omega \mathcal{A}_{[k]}^* (y_\mathrm{FP} - y^\delta_{k,:})$ \Comment{Kaczmarz step}\\
        \Comment{Shift correction step}\\
        $s_{[k]} \gets s_{[k]} - \lambda \frac{\Delta d}{m} \mathop{\rm shift\_cross\_corr}(y_\mathrm{FP}, y^\delta_{k,:})$\\ %
        $\mathcal{A}_{[k]} \gets \text{model at angle}\ k\ \text{with object shift}\ s_{[k]}$\\
    }
}
return $x$
\end{algorithm}

In \cref{alg:dremel}, $\mathrm{shift\_cross\_corr}(y_\mathrm{FP}, y^\delta_{k,:})$ determines the shift on the detector axis that maximizes the cross-correlation between the single-angle projection vectors $y_\mathrm{FP},y^\delta_{k,:}\in\mathbb{R}^{2\cdot P + 1}$ (or transformed versions of these), which is multiplied with the detector pixel size $\Delta d$ divided by the magnification $m$ (which is given by the source-to-detector and source-to-object distances) to estimate the object shift orthogonal to the current projection direction, and additionally scaled by an optional relaxation parameter $\lambda\in (0,1]$.
To achieve sub-pixel resolution for the shifts, it is proposed in \cite{Dremel_2018} to upsample the projections before computing the cross-correlation, for which we choose the super-resolution factor $f_\text{SR} = 2$.
Specifically, the shift computation $\mathrm{shift\_cross\_corr}(\cdot,\cdot)$ is given by
\begin{align*}
 \mathrm{shift\_cross\_corr}(z, y) &= \frac{1}{f_\text{SR}} \cdot \left(-\left\lfloor\frac{n_\text{SR}}{2}\right\rfloor + \mathop{\rm mod}\left(\left\lfloor\frac{n_\text{SR}}{2}\right\rfloor + \argmax_l \,\mathrm{cross\_corr}\!\left(\mathrm{up}_{f_\text{SR}}(z), \mathrm{up}_{f_\text{SR}}(y)\right)_l, n_\text{SR}\right) \right)\\
 \mathrm{cross\_corr}(z_\text{SR}, y_\text{SR}) &= \mathop{\rm post}\!\left(\mathcal{F}^{-1}\left(\mathcal{F}(\mathop{\rm pre}(y_\text{SR}))^* \mathcal{F}(\mathop{\rm pre}(z_\text{SR}))\right)\right),
\end{align*}
where $\mathrm{up}_{f_\text{SR}}(\cdot)$ is the upsampling operator, $n_\text{SR} = f_\text{SR} (2 \cdot P + 1)$ is the dimension of the upsampled projection vectors $z_\text{SR},y_\text{SR}\in\mathbb{R}^{n_\text{SR}}$, and $\mathcal{F}$ is the discrete Fourier transform; $\mathop{\rm pre}(\cdot)$ and $\mathop{\rm post}(\cdot)$ are pre- and post-processing functions that we describe next.
In \cite{Dremel_2018}, an edge-filter has been suggested as a pre-processing when correlating low-contrast signals, which we did not find to be beneficial in our validation.
We instead first perform mean subtraction in $\mathop{\rm pre}$ followed by zero-padding, while in $\mathop{\rm post}$ the correlation signal is cropped to the original size of $z$ and $y$, and a weighting is applied that boosts the correlation values for higher shifts as they are biased to be smaller due to the zero-padding.
Applying the latter weighting lead to improved results on validation data.

The learning rate $\omega$ and the relaxation parameter $\lambda$ are also found on validation data, we use $\omega = 1$ and $\lambda = 1$ (no relaxation).
In \cite{Dremel_2018}, it is suggested that the model correction step can be omitted in some initial warm-up iterations; we tried this on validation data and found that starting the correction in the first iteration performed best.
The ${\rm stoppingCriterion}$ is fulfilled when the maximum number of iterations is reached, which is set to $32$ in all our experiments.
\clearpage

\section{Pseudocode RESESOP}
\label{app:resesop}
\vspace{-0.3cm}
\RestyleAlgo{ruled}
  
\begin{algorithm}[h]
\setstretch{1.2}
\newcommand{\prodUF}{\langle u_\text{\rm old}, \tilde{x}_\text{\rm new} \rangle}
\newcommand{\prodUU}{\langle u_\text{\rm new}, u_\text{\rm old} \rangle}
\SetKw{KwInit}{Init:}
\caption{RESESOP Pseudocode for two search directions from two consecutive iterations}
\label{alg:resesop}
\newlength\mylen
\SetKwInOut{MyInput}{Input}  
\MyInput{$\mathcal{A} \gets$ weight matrix for each detector point and scanner angle\\
$\eta,\delta \gets \textrm{global model inexactness/noise level}$\\
$y^\delta \gets \textrm{perturbed sinogram}$}

\KwResult{$x$}
\SetKwInOut{MyInit}{Init}  
\MyInit{${\rm scannerAngles} \gets 1 \ldots K$, \\
${\rm detectorPixels} \gets 1 \ldots 2 \cdot P + 1$,\\
$x \gets \textrm{empty phantom}$}

\While{$\neg {\rm stoppingCriterion}$}{
    \For{angle index $k$ in {\rm scannerAngles}}{
        \For{detector index $l$ in {\rm detectorPixels}}{
            $w_{k,l} \gets \mathcal{A}_{k,l} \cdot x - y^\delta_{k,l}$ \Comment{define search direction} 

            \uIf{$\norm{w_{k,l}} \leqslant \tau \cdot (\delta + \eta_{k,l})$ }{\vspace*{-1em}\Comment{the discrepancy principle is satisfied}\\
                ${\rm stoppingCondition}_{k,l} \gets {\rm True}$ \\
                continue
            }\vspace*{0.75em}
            
            $\alpha \gets w_{k,l} \cdot y_{k,l}$\\
            $\xi \gets (\delta + \eta_{k,l}) \cdot \norm{w_{k,l}}$\\
            $u_\text{new} \gets (\mathcal{A}_{k,l})^* \cdot w_{k,l}$ \Comment{compute first search direction $u_\text{new}$}\\
             $\tilde{x}_\text{new} \gets x - \frac{\norm{w_{k,l}} \cdot (\norm{w_{k,l}} - (\delta + \eta_{k,l}))}{\norm{u_\text{new}}^2 }\cdot u_\text{new}$ \Comment{projection w.r.t. first search direction}
            
            \uIf{is first iteration}{\vspace*{-1em}\Comment{solely one search direction in first iteration}\\
                $x \gets \tilde{x}_\text{new}$
            }
            \Else{
                
                \uIf{$\prodUF > \alpha_\text{old} + \xi_\text{old}$}{\vspace*{-1em}\Comment{consider second search direction $u_\text{old}$}\\
                    $t \gets \frac{\prodUF - (\alpha_\text{old} + \xi_\text{old})}{\norm{u_\text{new}}^2 \cdot \norm{u_\text{old}}^2 - \prodUU^2}$
                }
            
                \uElseIf{$\prodUF < \alpha_\text{old} - \xi_\text{old}$}{
                    $t \gets \frac{\prodUF - (\alpha_\text{old} - \xi_\text{old})}{\norm{u_\text{new}}^2 \cdot \norm{u_\text{old}}^2 - \prodUU^2}$
                }
                \Else{
                     $t \gets 0$
                }
                \vspace{5pt}
                $x \gets \tilde{x}_\text{new} + \prodUU \cdot t \cdot u_\text{new} - {} \norm{u_\text{new}}^2 \cdot t \cdot u_\text{old}$\\

            }
            $u_\text{old} \gets u_\text{new}$\Comment{save variables from this iteration}\\
            $\alpha_\text{old} \gets \alpha$\\
            $\xi_\text{old} \gets \xi$\\
        }
    }
}
return $x$
\end{algorithm}

In \cref{alg:resesop}, the ${\rm stoppingCriterion}$ is fulfilled when either (i) ${\rm stoppingCondition}_{k,l}$ has been set to ${\rm True}$ for all angle and detector indices $k,l$, or (ii) a maximum number of entire sweeps over all angle and detector indices is reached, which is set to $20$ in all our experiments. Norm and scalar product in $\mathcal{X}$ are discretizations of the $L^2([-1,1]^2)$ norm and scalar product using piecewise constant basis functions with pairwise disjoint rectangular supports.

\clearpage
\FloatBarrier
\section{Training details}
\label{app:network_training}
\subsection{UNet}
The details of the UNet architecture used in our experiments are depicted in Figure \ref{fig:unet_scheme}. We trained the network for 500 epochs and used the network parameters for testing that performed best during validation. As a loss function we chose the mean squared error. We used the Adam optimizer with an initial learning rate of 0.001 and a batch size of 16.
\begin{figure}
    \centering
    \includegraphics[width=\columnwidth]{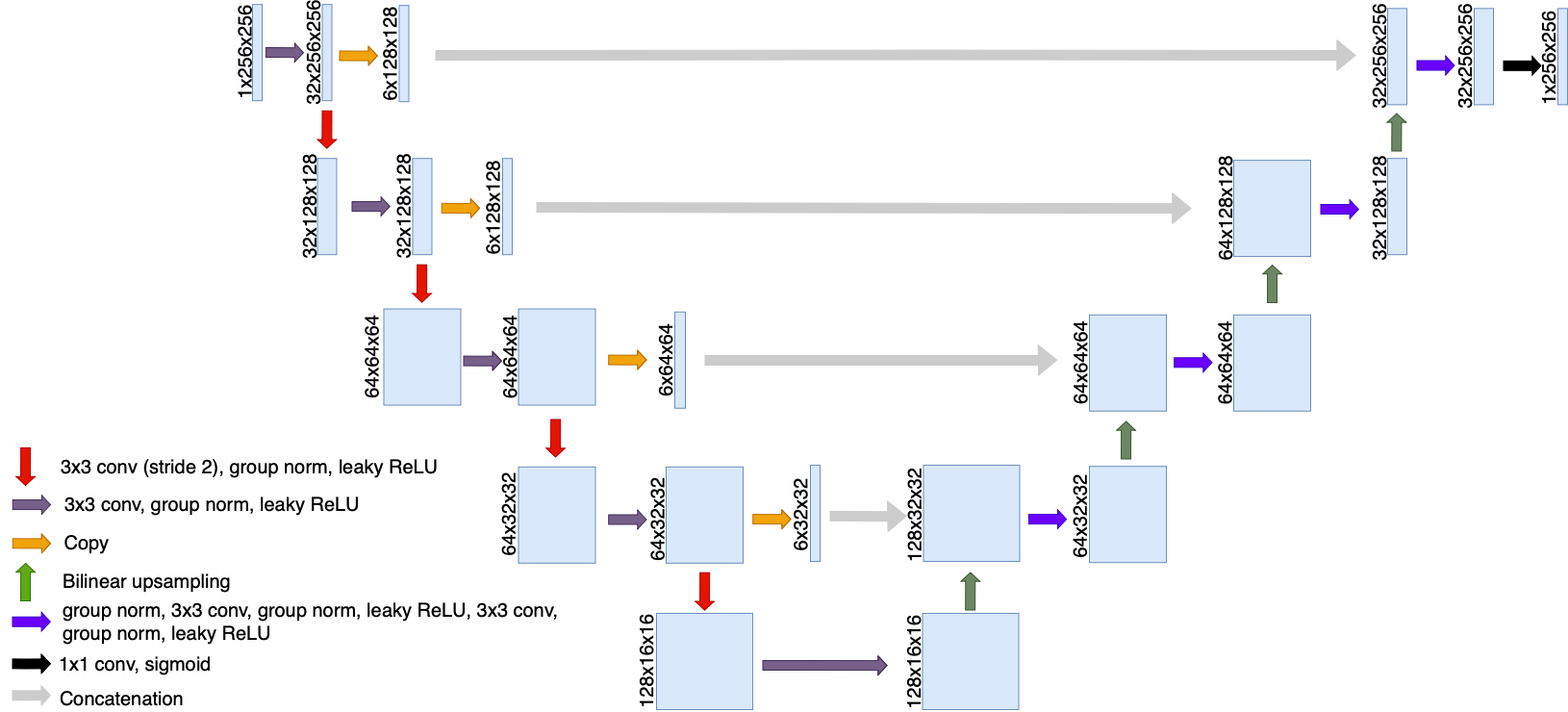}
    \caption{UNet architecture}
    \label{fig:unet_scheme}
\end{figure}
\subsection{CiNN / CiNNRes}
The details of the CiNN architecture used in our experiments are depicted in Figure \ref{fig:cinn_scheme}. We trained the network for 500 epochs and used the network parameters for testing that performed best during validation. As a loss function we chose the mean squared error. We used the Adam optimizer with an initial learning rate of 0.00001 and a batch size of 16. The coupling blocks follow the NICE design. The conditioning network is a ResNet, which is trained jointly with the invertible network, i.e.\ as part of the conditional normalizing flow.
\begin{figure}
    \centering
    \includegraphics[width=\columnwidth]{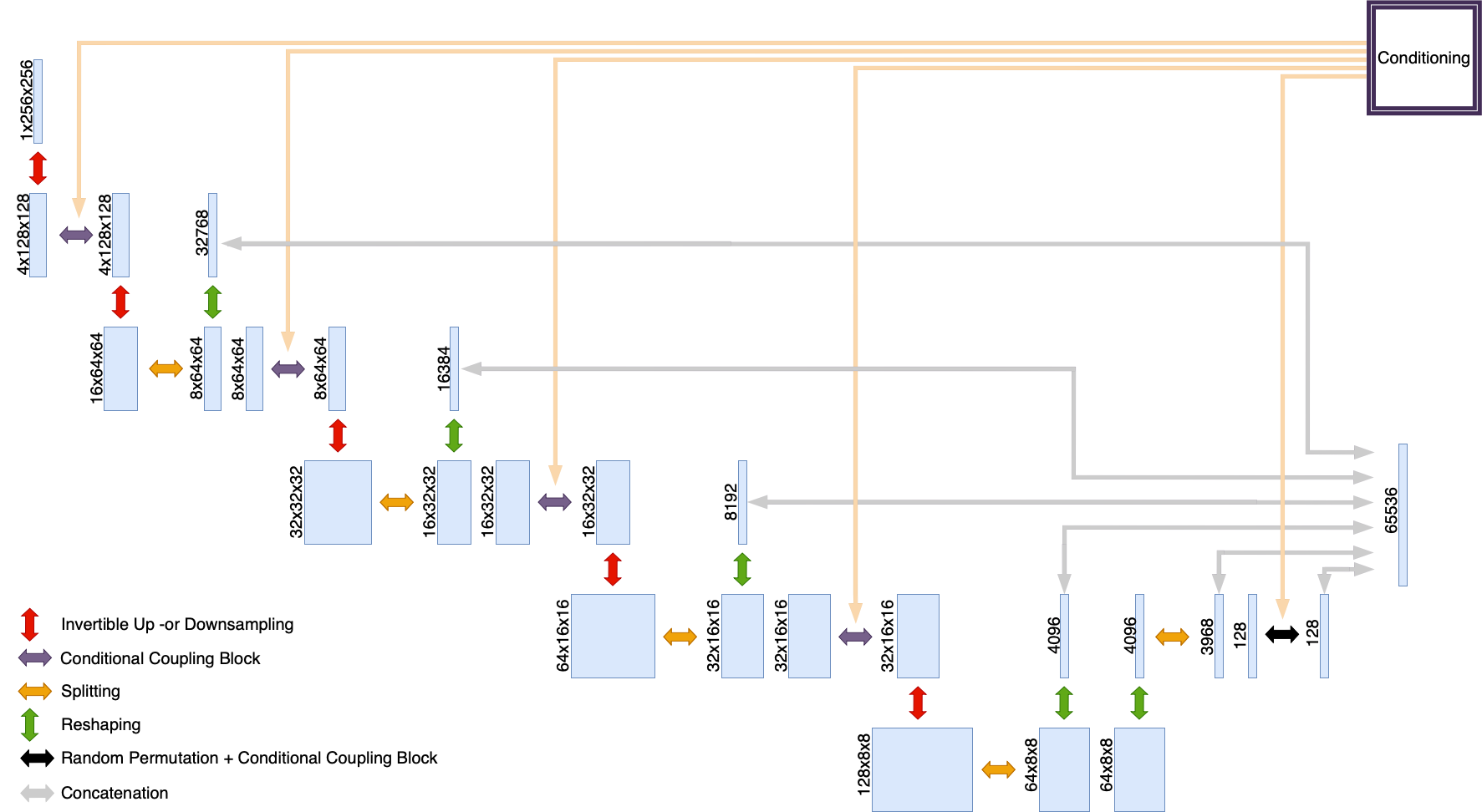}
    \caption{CiNN architecture. Part of the information is factored out at different scales, and concatenated to form the latent vector $z$. The other part is continued to be transformed, using additive coupling blocks, which receive conditioning inputs from a ResNet, and invertible down-sampling.}
    \label{fig:cinn_scheme}
\end{figure}

\subsection{CiUNet / CiUNetRes}
The details of the iUNet architecture used in our experiments are depicted in Figure \ref{fig:iunet_scheme}. We trained the network for 500 epochs and used the network parameters for testing that performed best during validation. As a loss function we chose the mean squared error. We used the Adam optimizer with an initial learning rate of 0.00001 and a batch size of 16. The coupling blocks follow the NICE design. The conditioning network is a UNet with the same architecture as for the UNet post-processing, shown in Figure \ref{fig:unet_scheme}, which is trained jointly with the invertible network, i.e.\ as part of the conditional normalizing flow.
\begin{figure}
    \centering
    \includegraphics[width=\columnwidth]{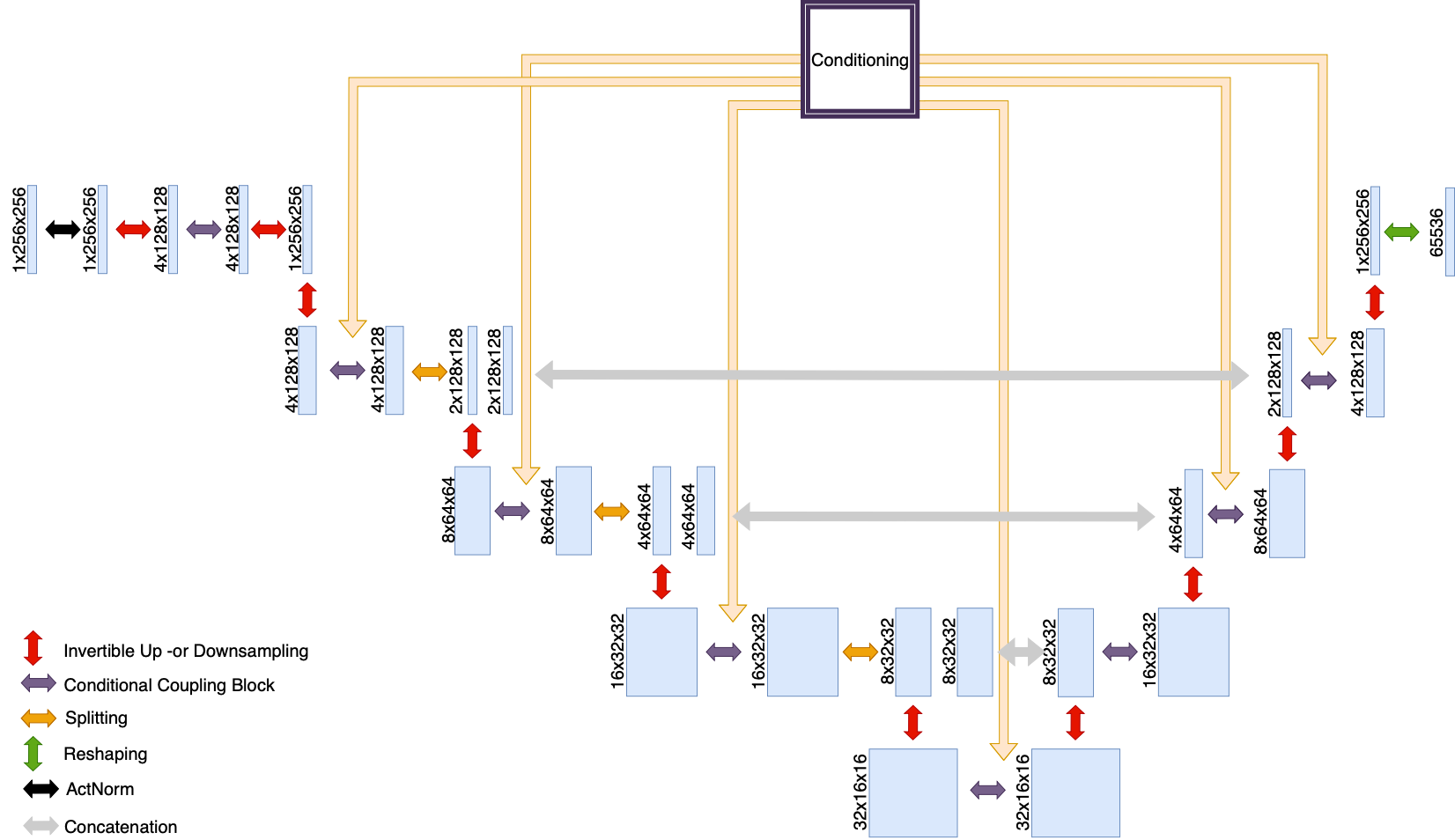}
    \caption{CiUNet architecture. The additive coupling blocks in the encoder and decoder receive conditioning inputs from the decoder and encoder of a non-invertible U-Net, respectively. Encoder and decoder are linked by concatenating skip connections.}
    \label{fig:iunet_scheme}
\end{figure}
\clearpage

\twocolumn
\section{Additional figures - non-perturbed parallel beam} 
\label{app:figs_non_perturbed_parallel}
\
\begin{figure}[h]
    \centering
    Reconstructions\\[2mm] 
    \includegraphics[width=\linewidth]{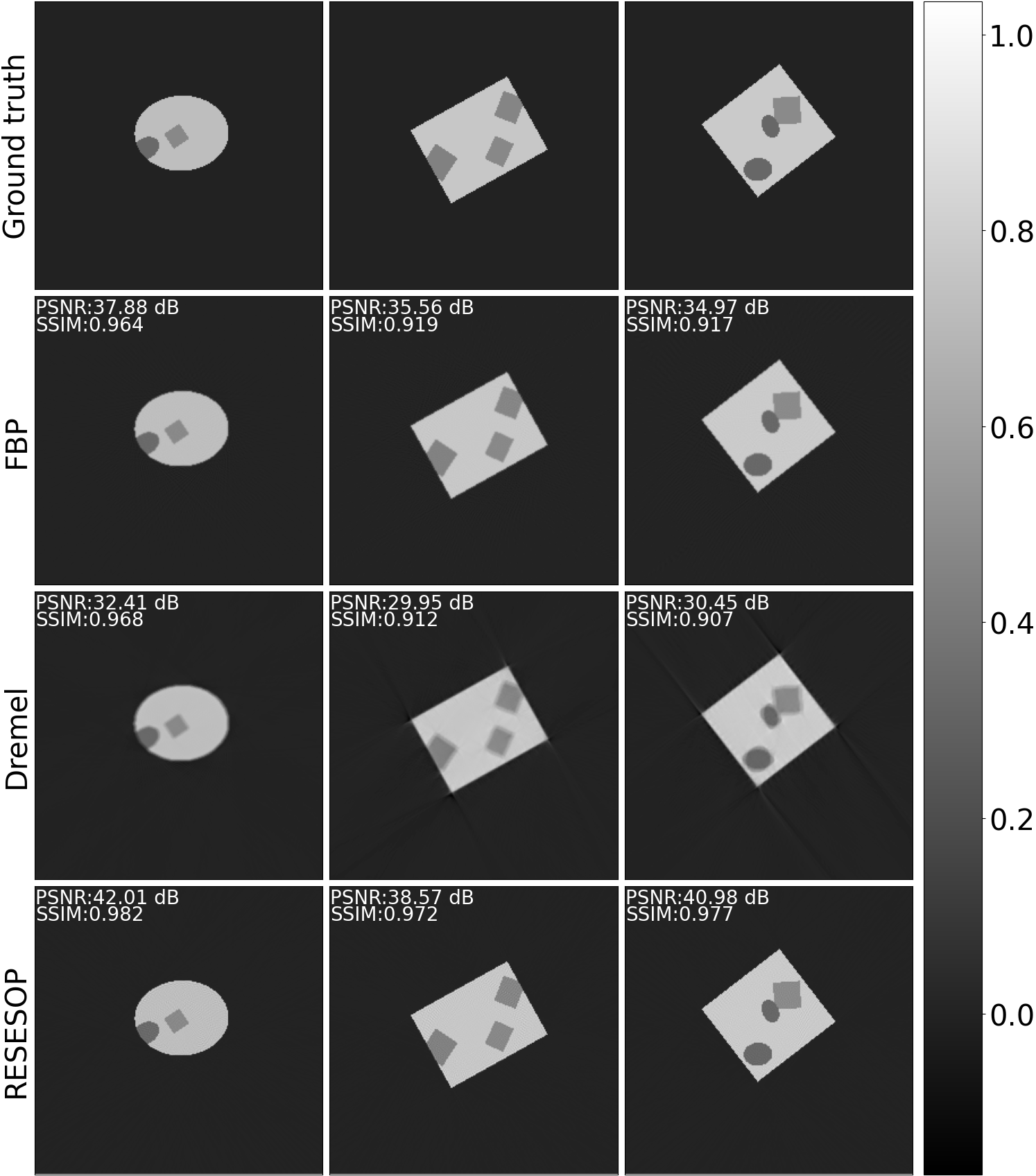}\\
    Differences\\[2mm] 
    \includegraphics[width=\linewidth]{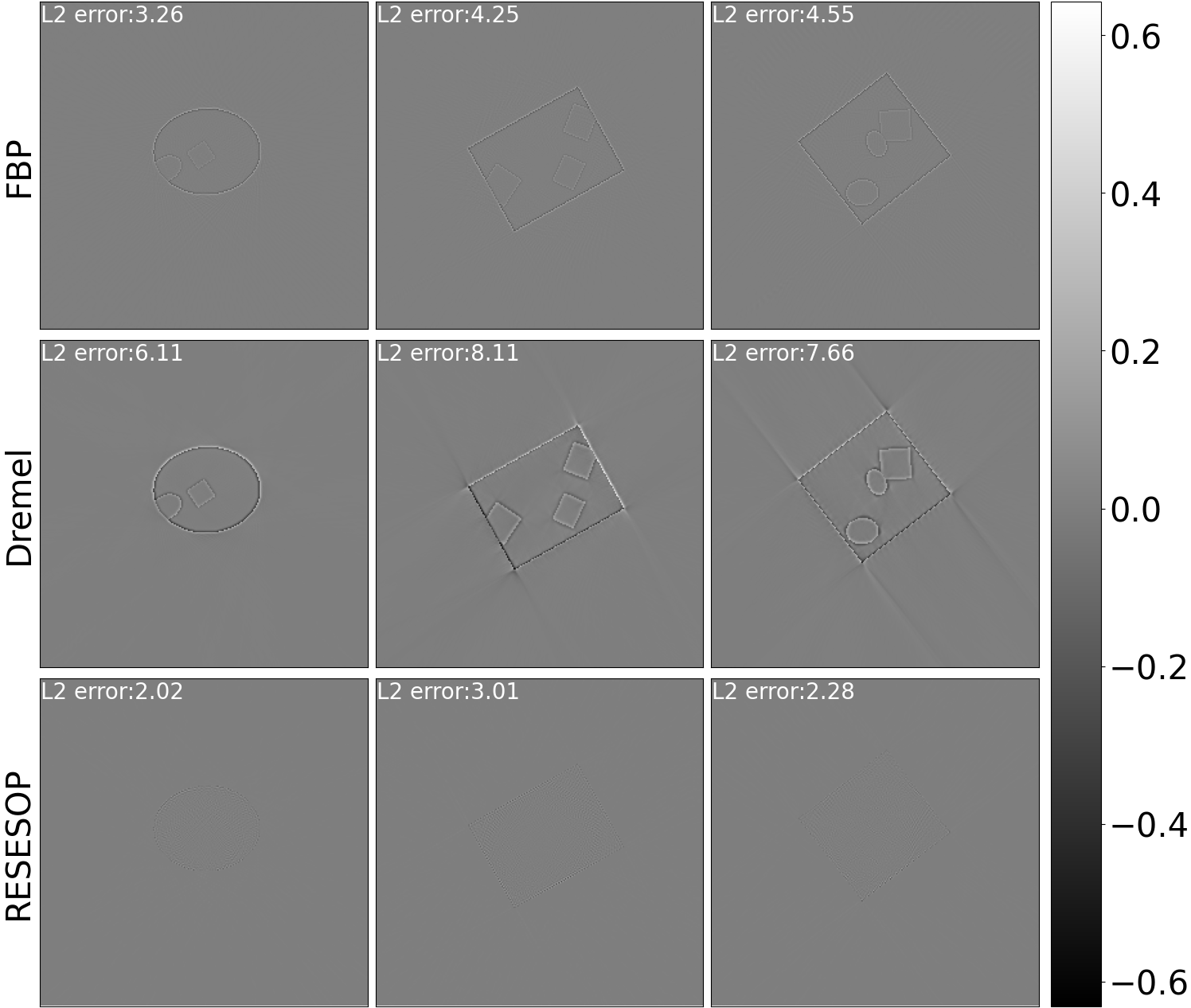}
    \caption{Image reconstructions and difference to ground truth using the non-trained reconstruction methods $\mathcal{T}$ on non-perturbed parallel beam data for the phantoms in \cref{phantoms}.}
    \label{fig:reco_unperturbed_parallel_std_algo}
\end{figure}

\begin{figure}[t]
    \centering
    Reconstructions\\[2mm]
    \includegraphics[width=\linewidth]{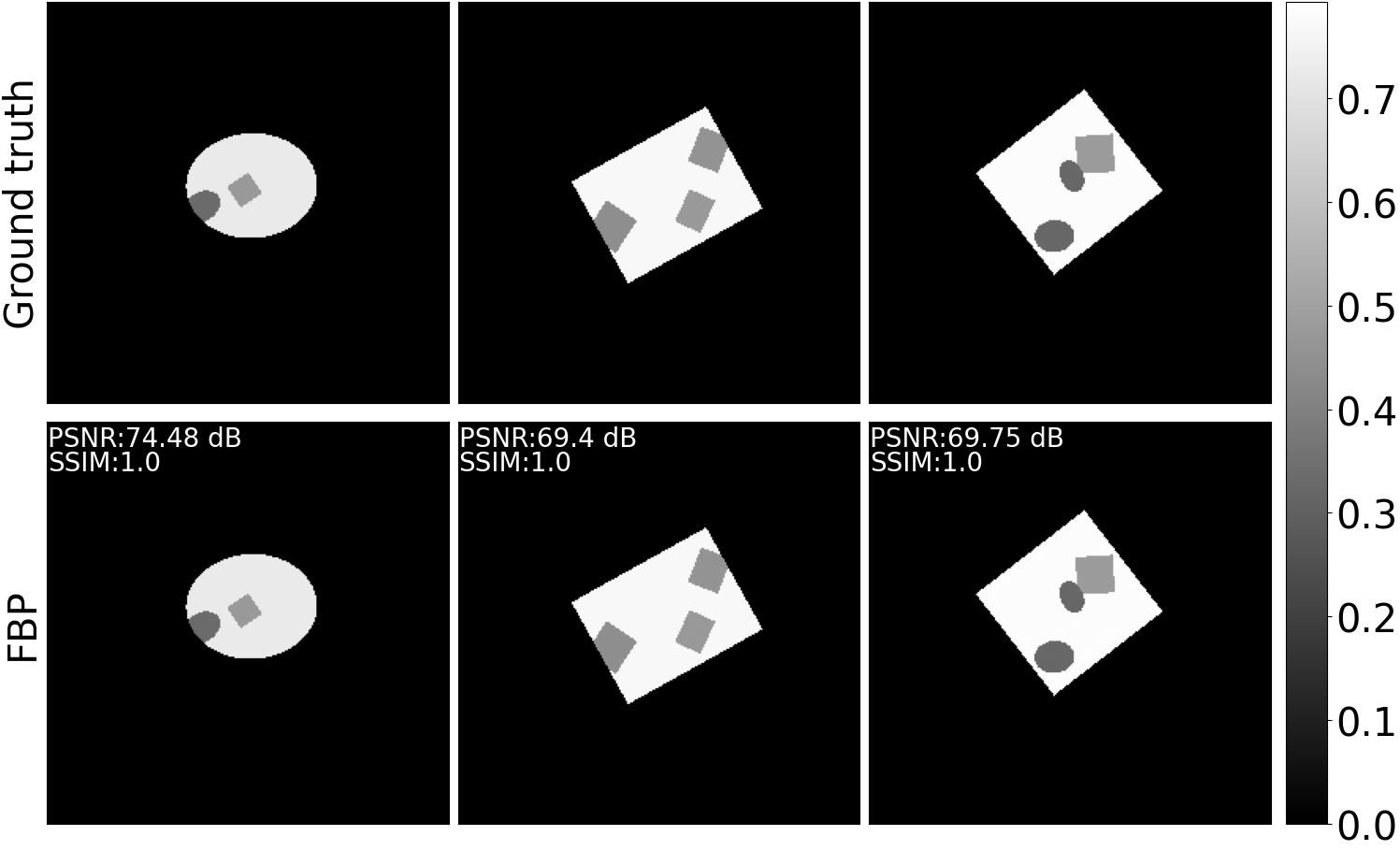}\\
    Differences\\[2mm] 
        \includegraphics[width=\linewidth]{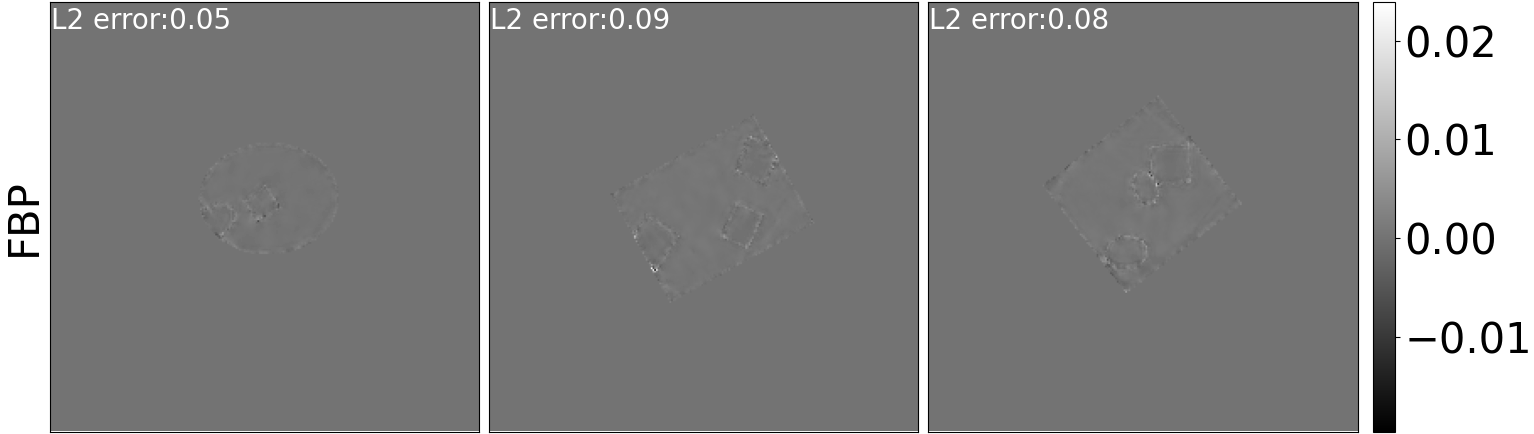}
    \caption{Image reconstructions $x_\mathrm{reco}$ and differences to ground truth using the non-trained reconstruction methods $\mathcal{T}$ combined with \textbf{UNet} post-processing on non-perturbed parallel beam data for the phantoms in \cref{phantoms}.}
    \label{fig:reco_unperturbed_parallel_U-Net}
\end{figure}

\begin{figure}[b]
    \centering
     Reconstructions\\[2mm]
    \includegraphics[width=\linewidth]{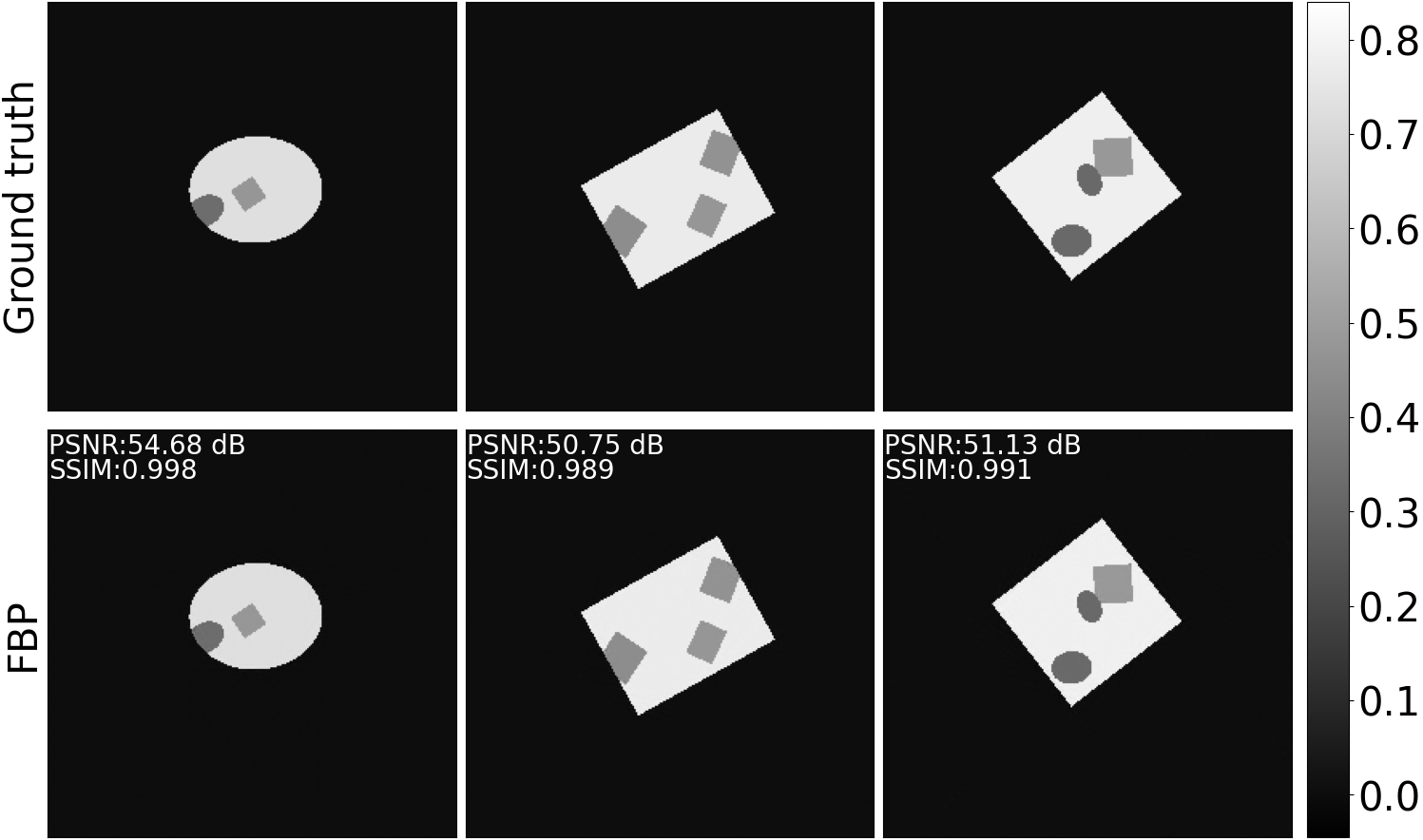}\\
    Differences\\[2mm] 
        \includegraphics[width=\linewidth]{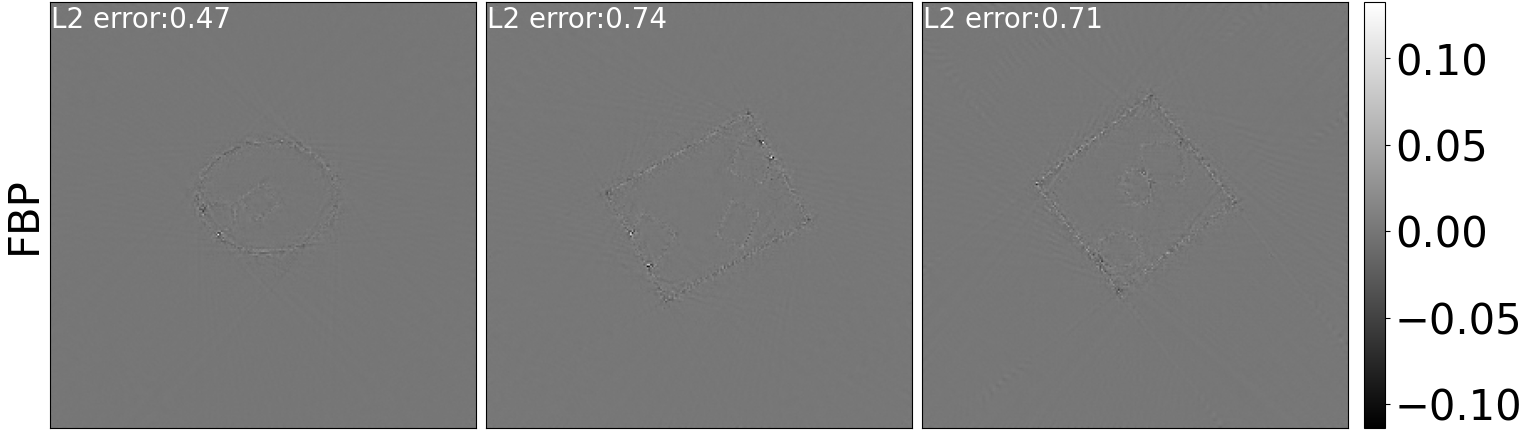}
    \caption{Image reconstructions $x_\mathrm{reco}$ and differences to ground truth using the non-trained reconstruction methods $\mathcal{T}$ combined with \textbf{CiNNRes} post-processing on non-perturbed parallel beam data for the phantoms in \cref{phantoms}.}
    \label{fig:reco_unperturbed_parallel_CINN}
\end{figure}

\begin{figure}[h]
    \centering
      Reconstructions\\[2mm]
    \includegraphics[width=\linewidth]{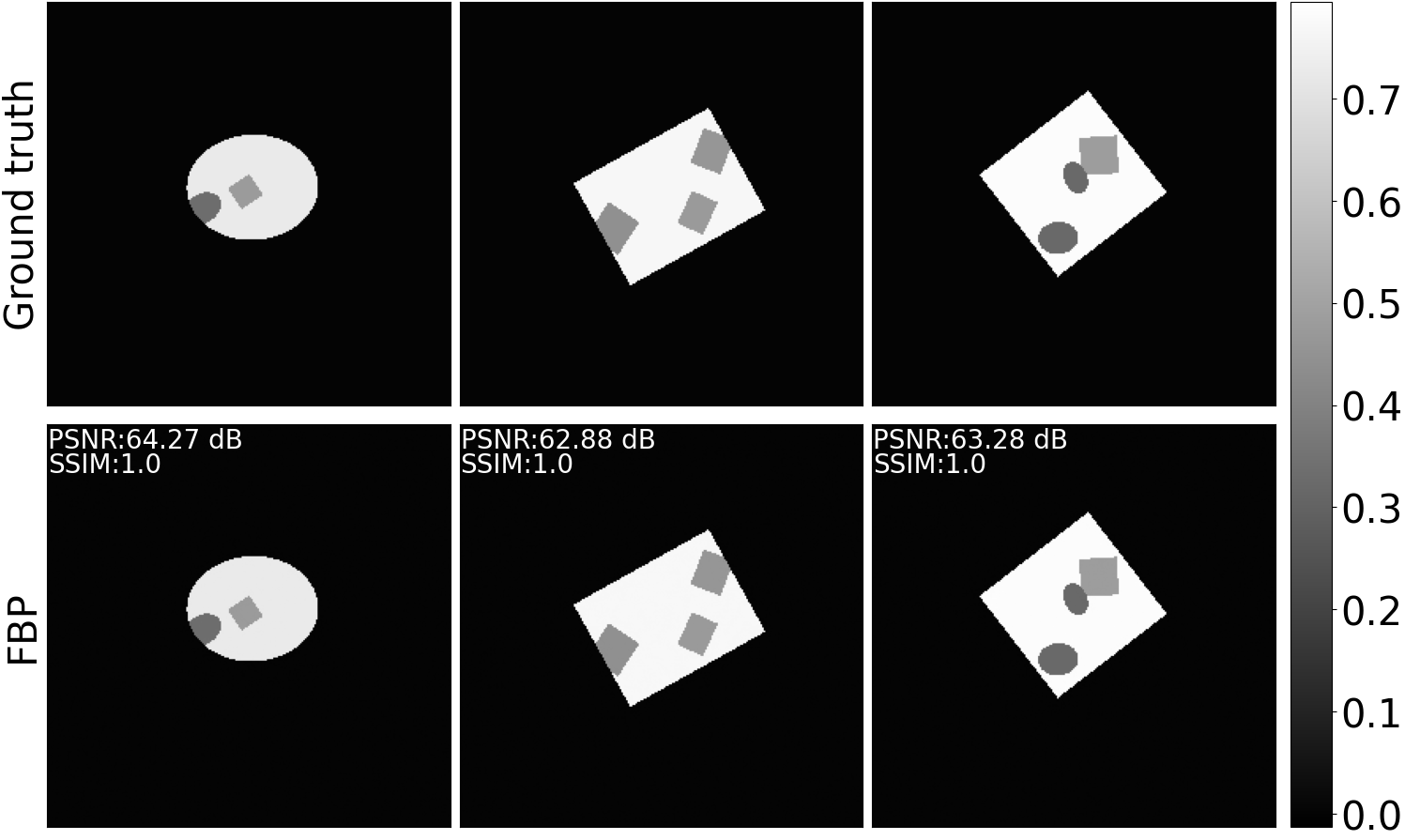}\\
    Differences\\[2mm] 
        \includegraphics[width=\linewidth]{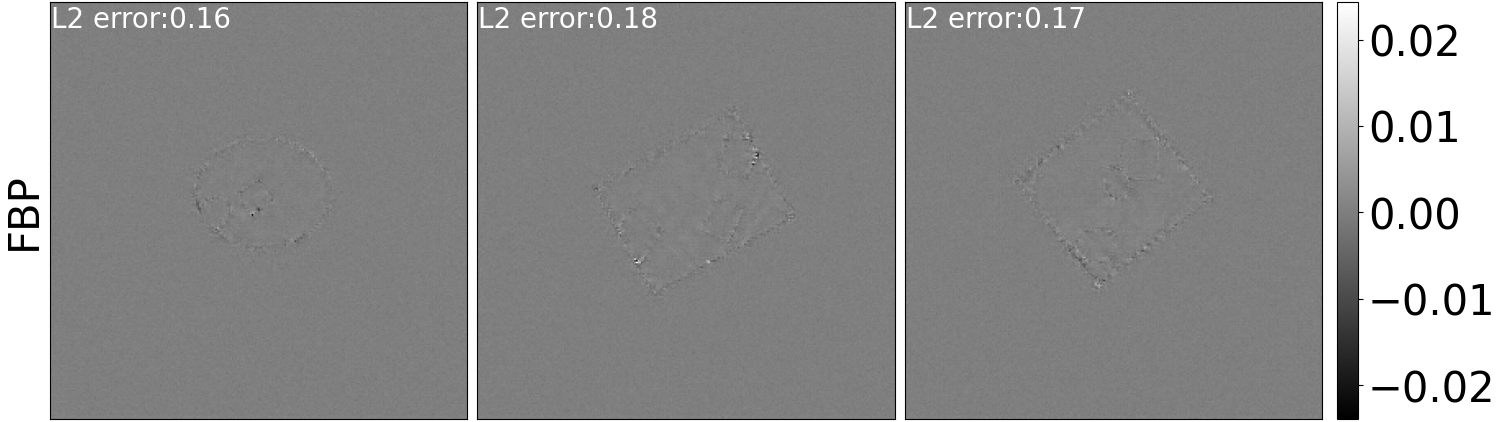}
    \caption{Image reconstructions $x_\mathrm{reco}$ and differences to ground truth using the non-trained reconstruction methods $\mathcal{T}$ combined with \textbf{CiUNetRes} post-processing on non-perturbed parallel beam data for the phantoms in \cref{phantoms}.}
    \label{fig:reco_unperturbed_parallel_CI-U-Net}
\end{figure}

\section{Additional figures - perturbed parallel beam} 
\label{app:figs_perturbed_parallel}
\
\begin{figure}[h]
    \centering
    \includegraphics[width=\linewidth]{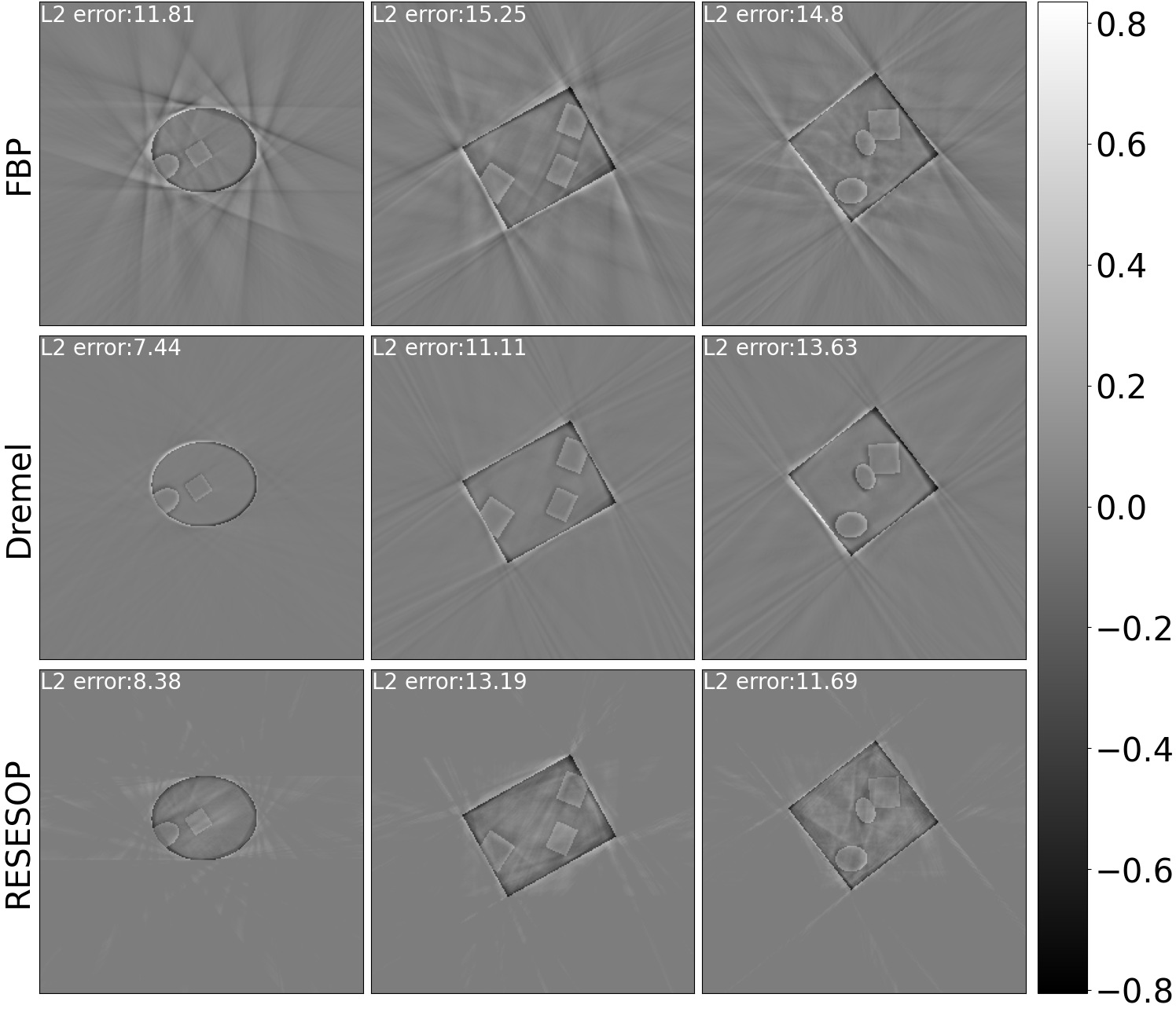}
    \caption{Difference of ground truth phantoms (\cref{phantoms}) and $\mathcal{T}$ reconstructions (\cref{fig:reco_perturbed}).}
    \label{fig:diff_reco_perturbed}
\end{figure}

\begin{figure}
    \centering
    \includegraphics[width=\linewidth]{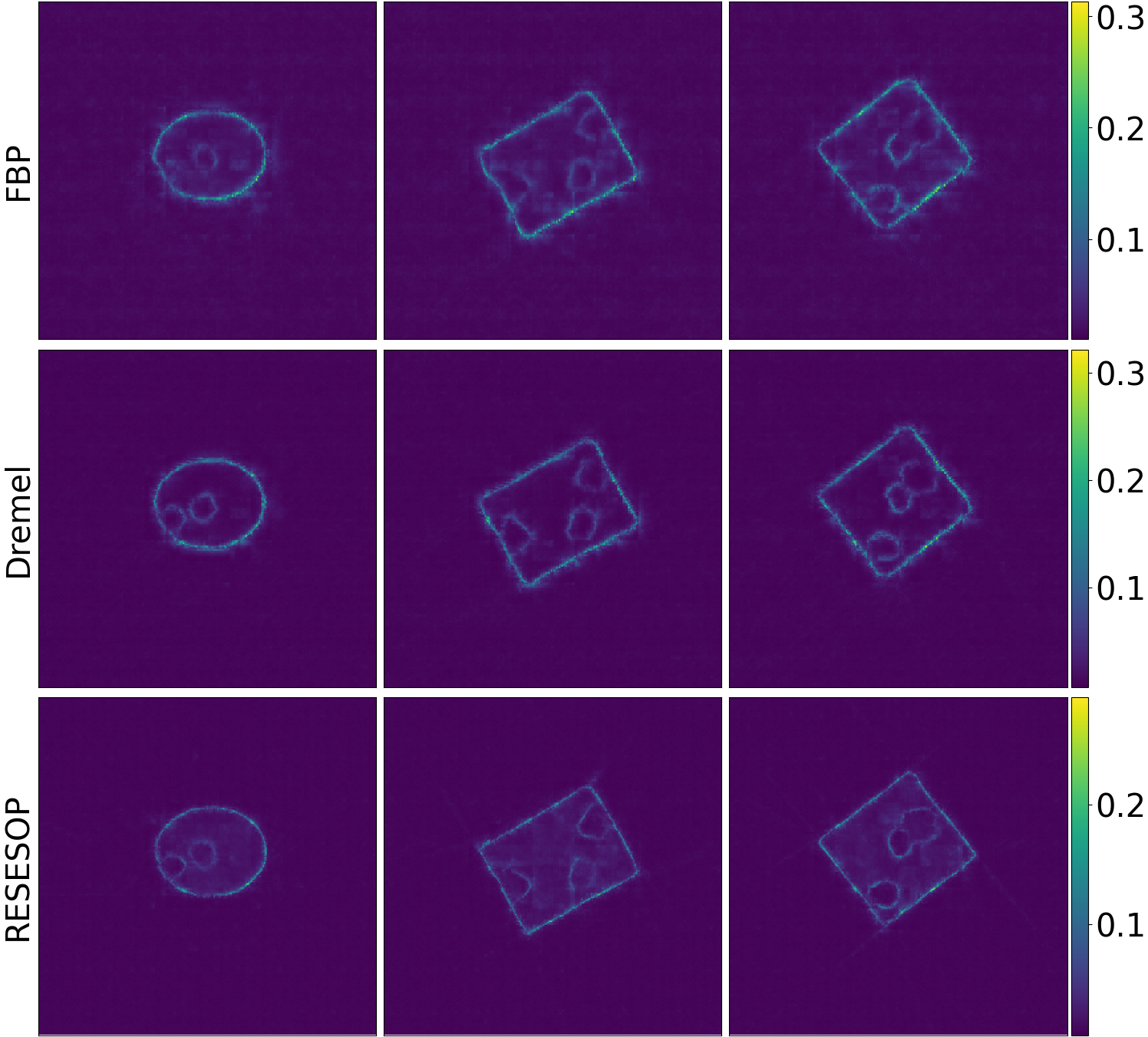}
    \caption{Standard deviation $\sigma_\mathrm{reco}$ according to \cref{sec:CINN} for the \textbf{CiNNRes} reconstructions in \cref{fig:cinn_perturbed}.}
    \label{fig:std_cinn_perturbed}
\end{figure}

\begin{figure}
    \centering
    Reconstructions\\[2mm] 
    \includegraphics[width=\linewidth]{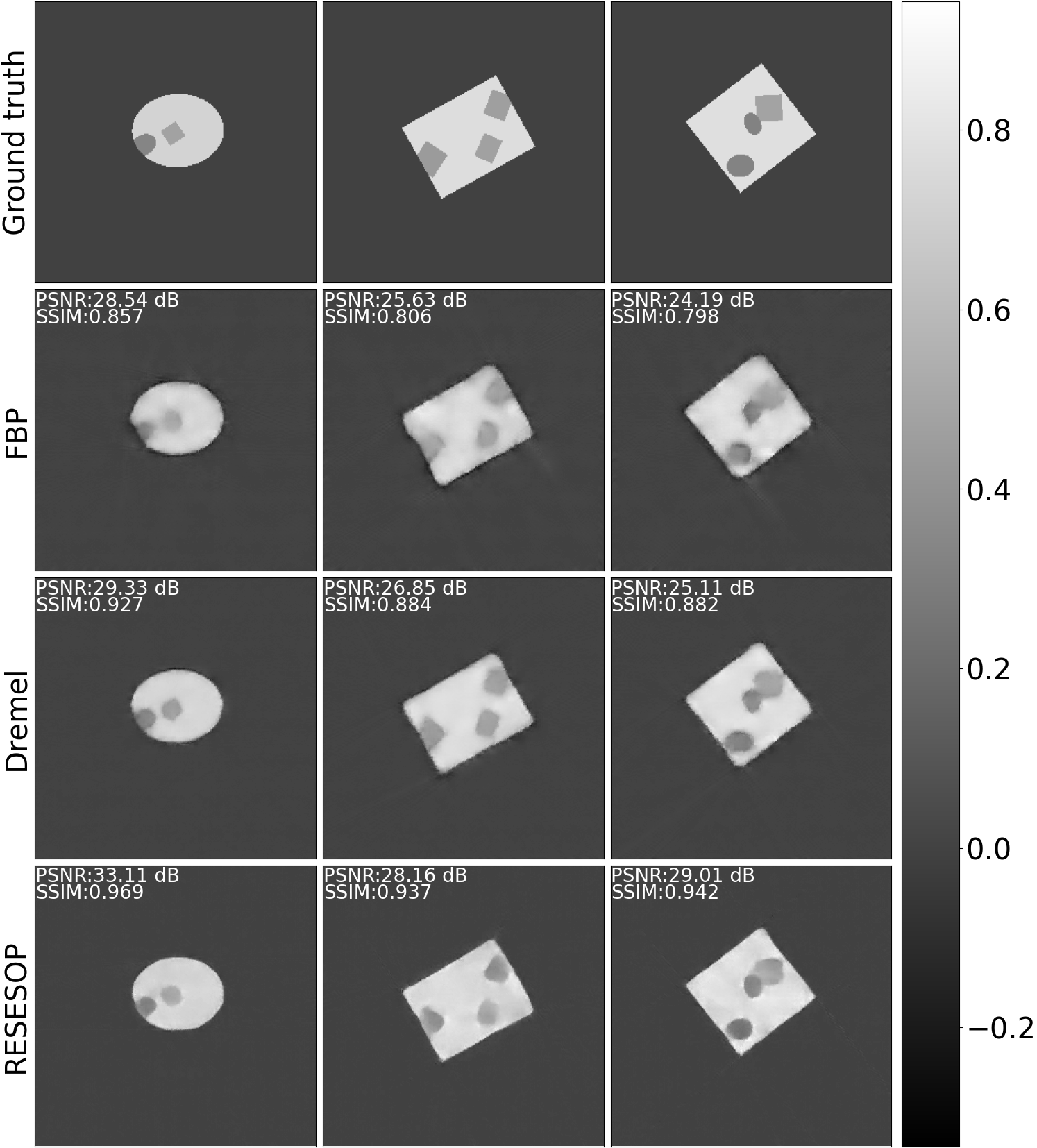}\\
    Differences\\[2mm] 
        \includegraphics[width=\linewidth]{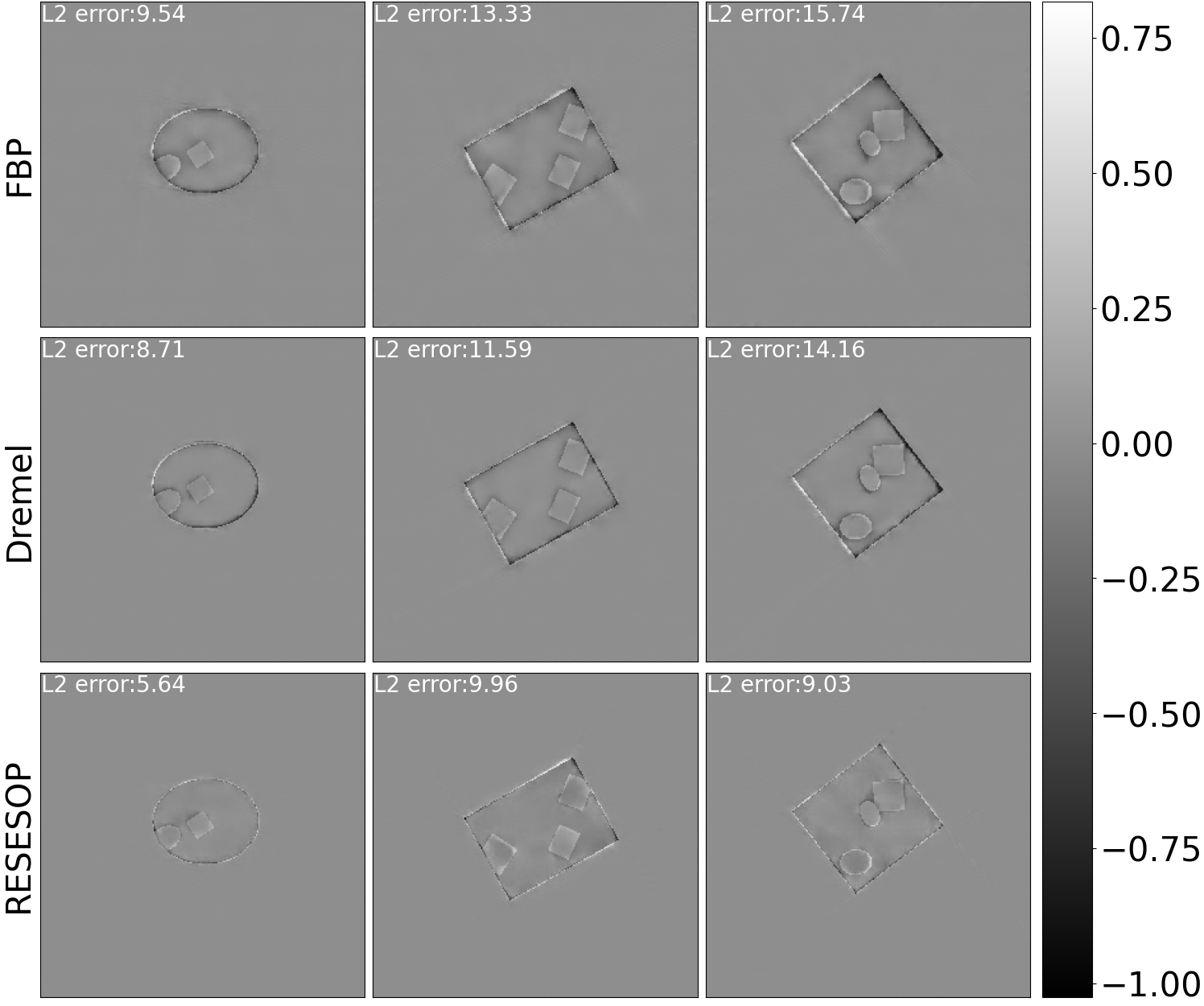}
    \caption{Image reconstructions $x_\mathrm{reco}$ and differences to ground truth using the non-trained reconstruction methods $\mathcal{T}$ combined with \textbf{CiNNRes} post-processing on perturbed parallel beam data for the phantoms in \cref{phantoms}.}
    \label{fig:cinn_perturbed}
\end{figure}

\clearpage
\section{Additional figures - perturbed fan beam} 
\
\label{app:figs_perturbed_fan}
\begin{figure}[h]
    \centering
    Reconstructions\\[2mm] 
    \includegraphics[width=\linewidth]{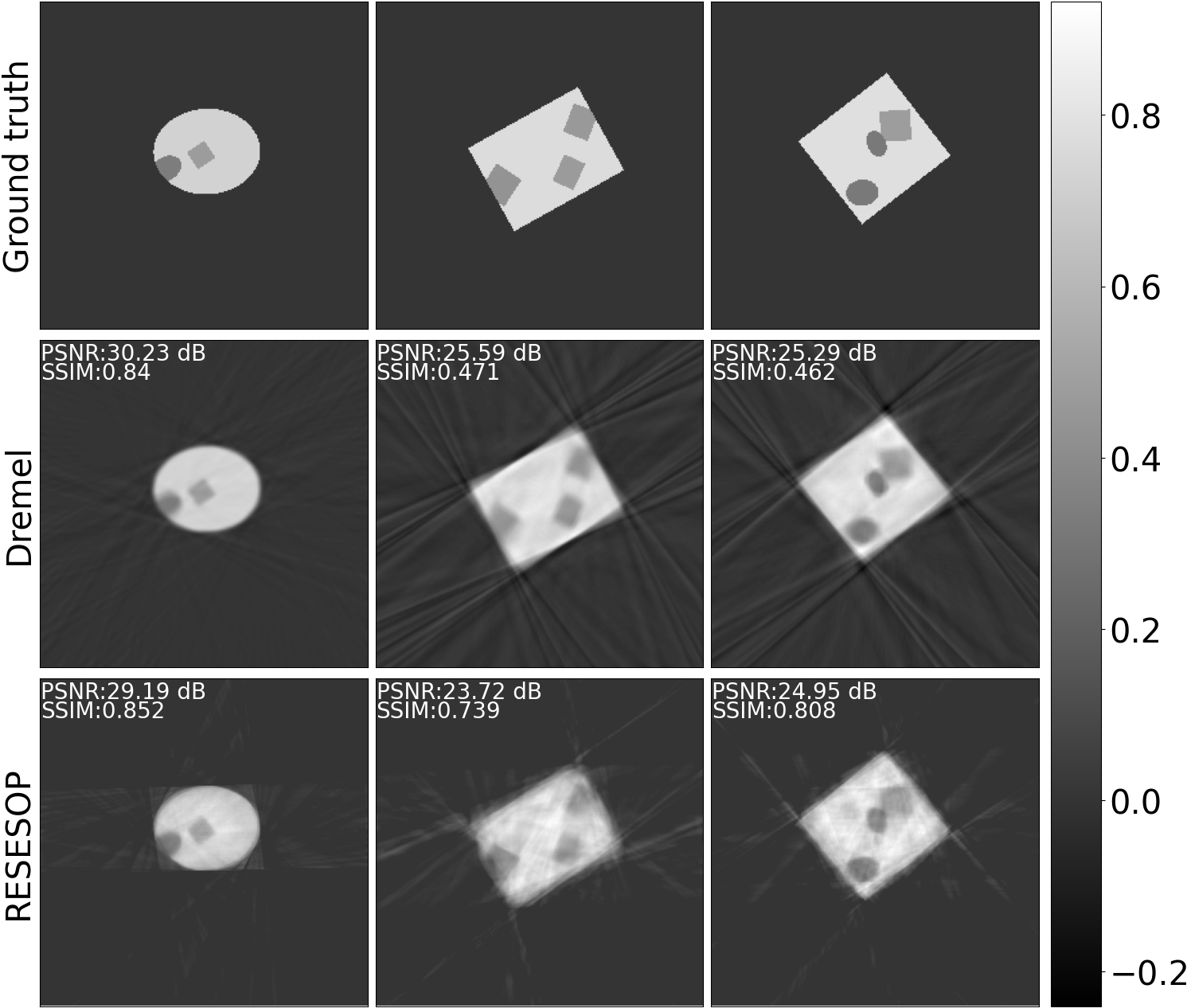}\\
    Differences\\[2mm] 
        \includegraphics[width=\linewidth]{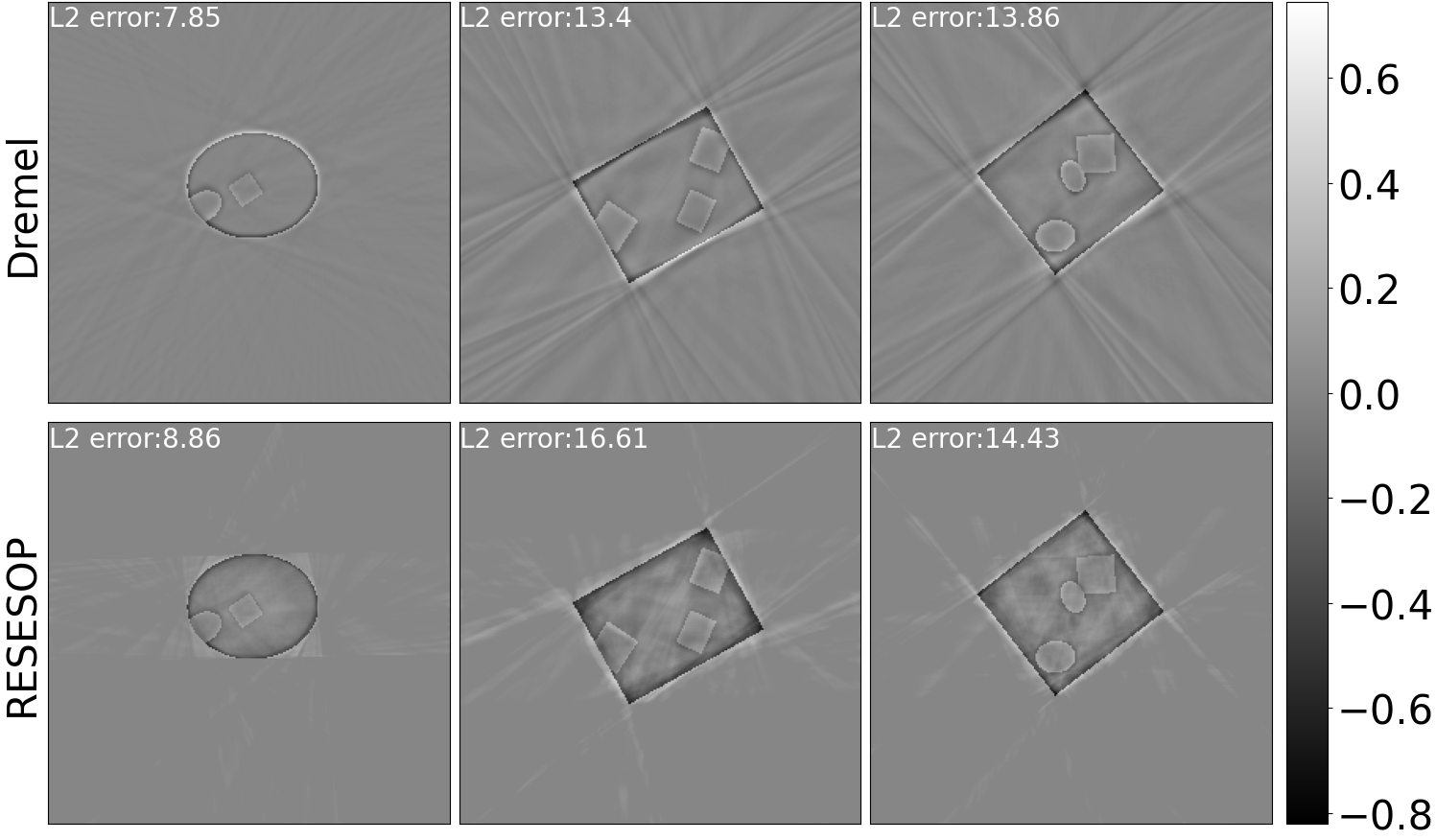}
    \caption{Image reconstructions and differences to ground truth using the non-trained reconstruction methods $\mathcal{T}$ on perturbed fan beam data for the phantoms in \cref{phantoms}.}
    \label{fig:reco_fan}
\end{figure}

\begin{figure}
    \centering
        Reconstructions\\[2mm] 
    \includegraphics[width=\linewidth]{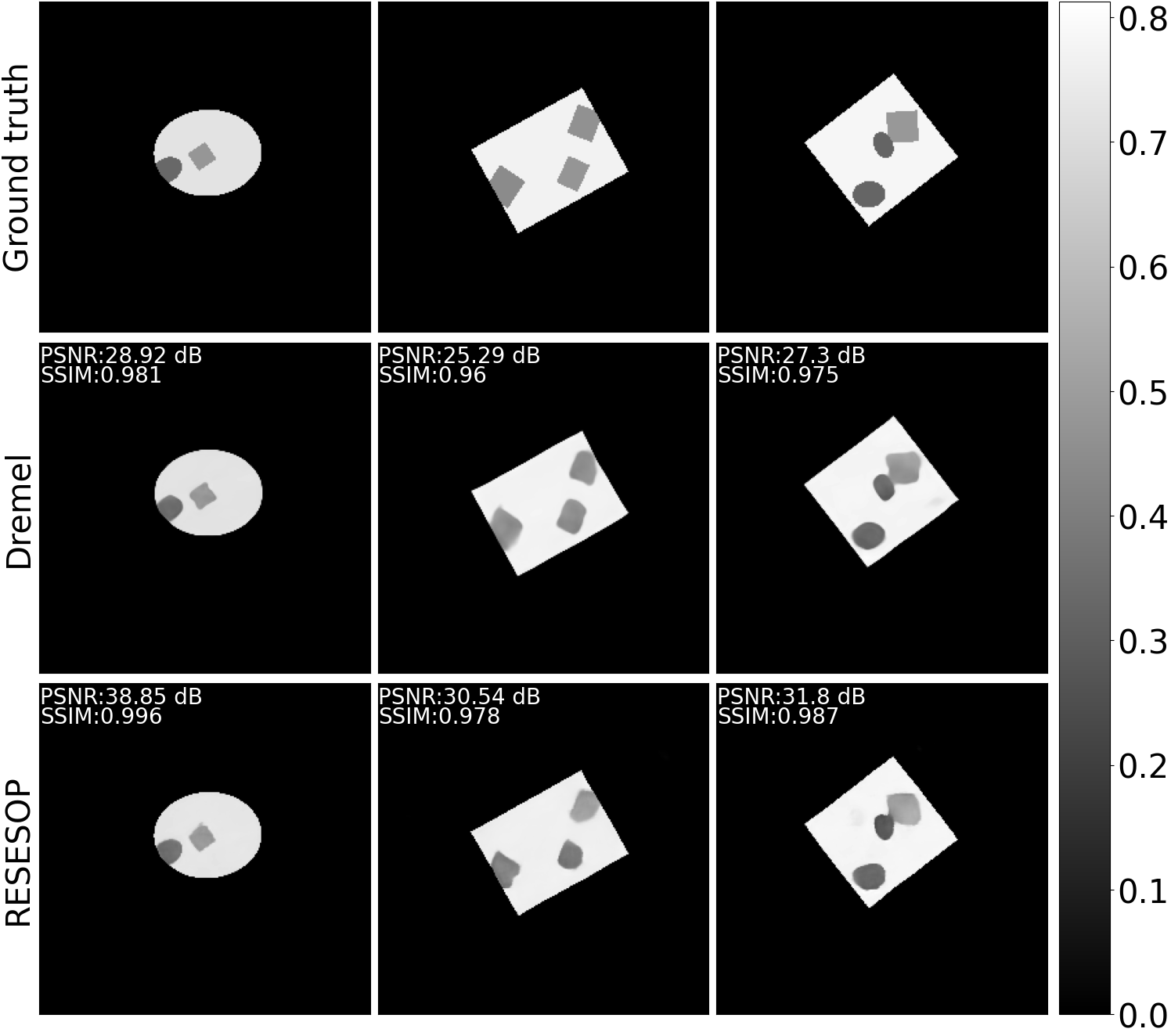}\\
    Differences\\[2mm] 
        \includegraphics[width=\linewidth]{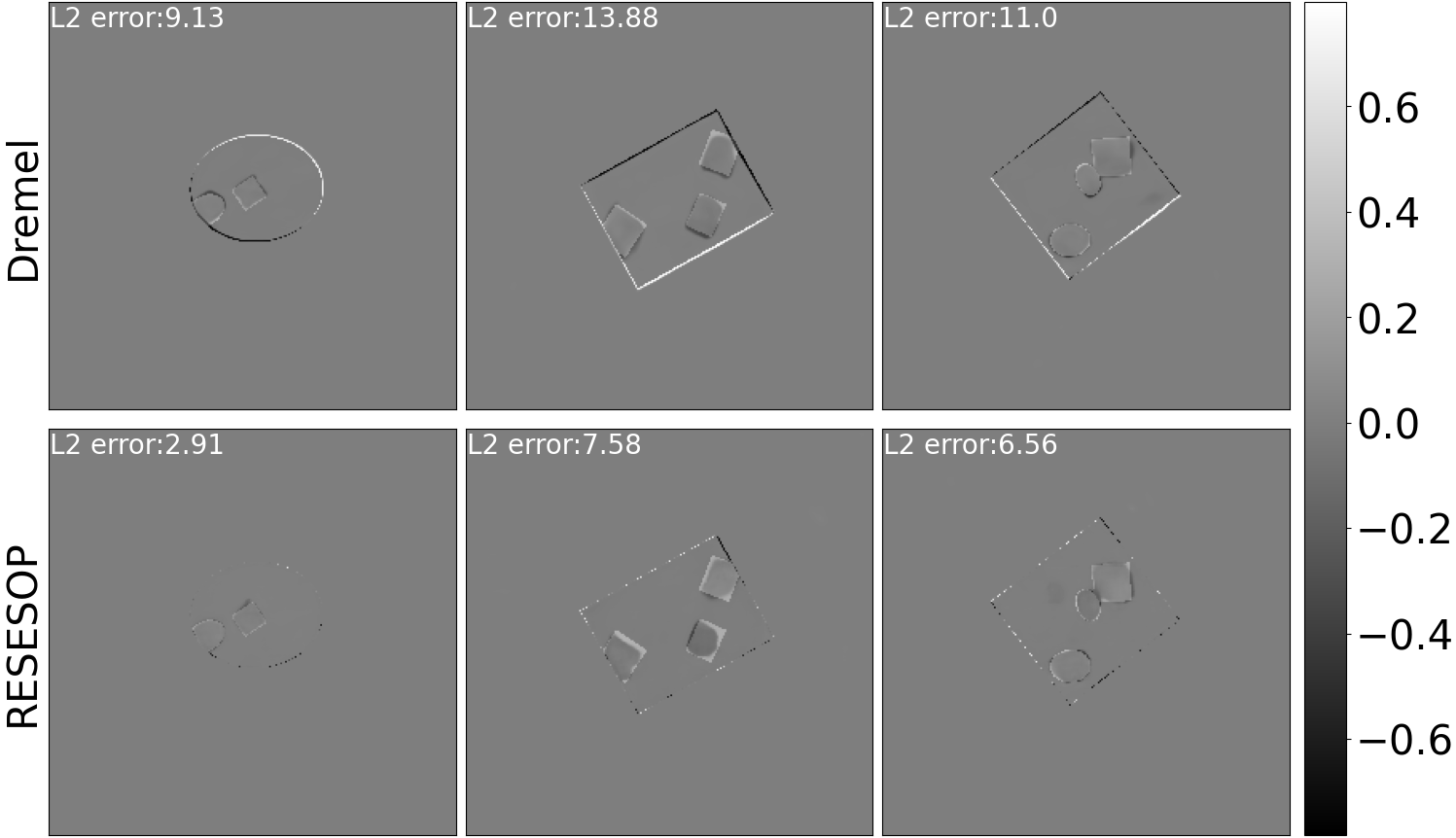}
    \caption{Image reconstructions $x_\mathrm{reco}$ and differences to ground truth using the non-trained reconstruction methods $\mathcal{T}$ combined with \textbf{UNet} post-processing on perturbed fan beam data for the phantoms in \cref{phantoms}.}
    \label{fig:unet_fan}
\end{figure}

\begin{figure}
    \centering
    Reconstructions\\[2mm] 
    \includegraphics[width=\linewidth]{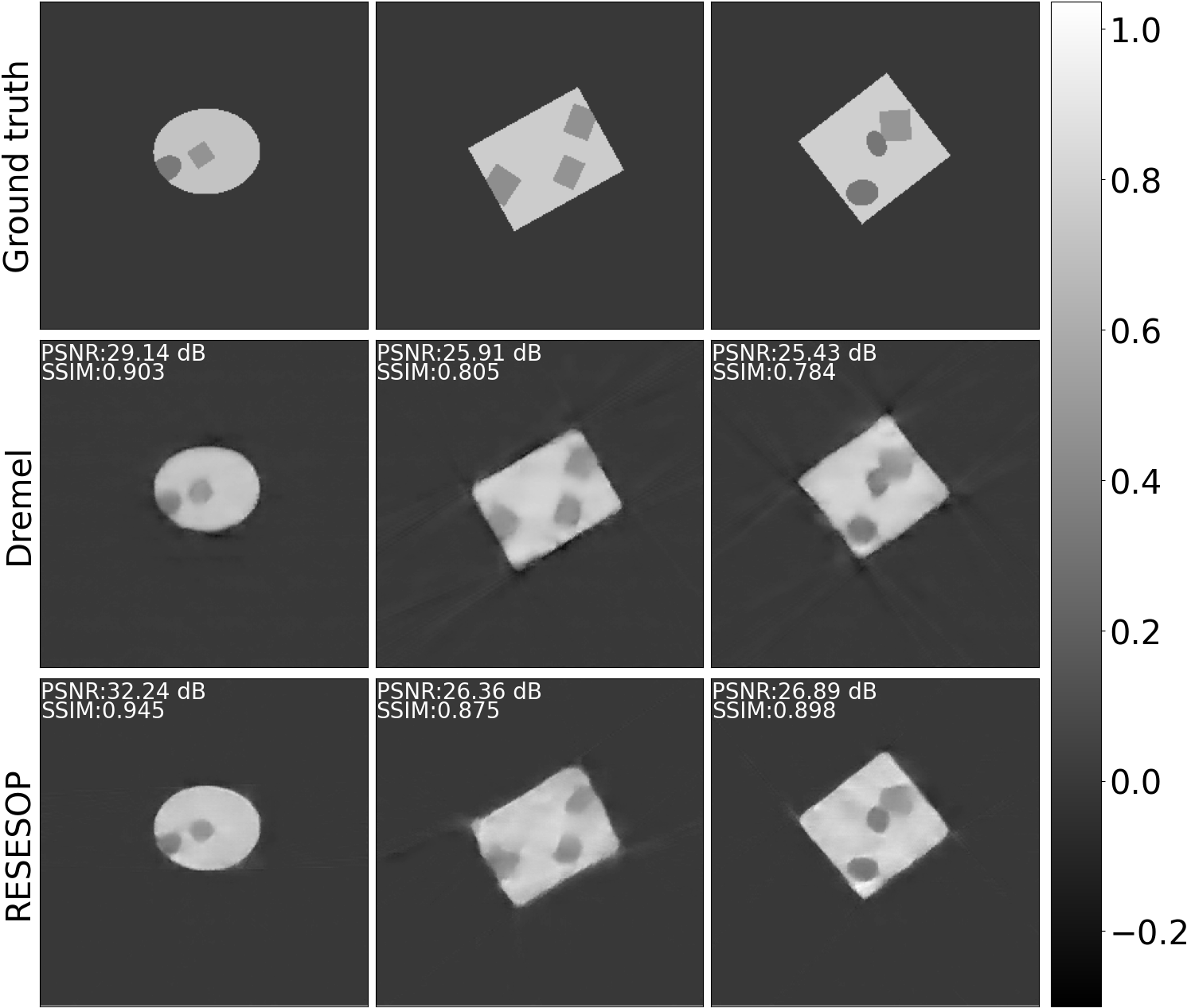}\\
    Differences\\[2mm] 
        \includegraphics[width=\linewidth]{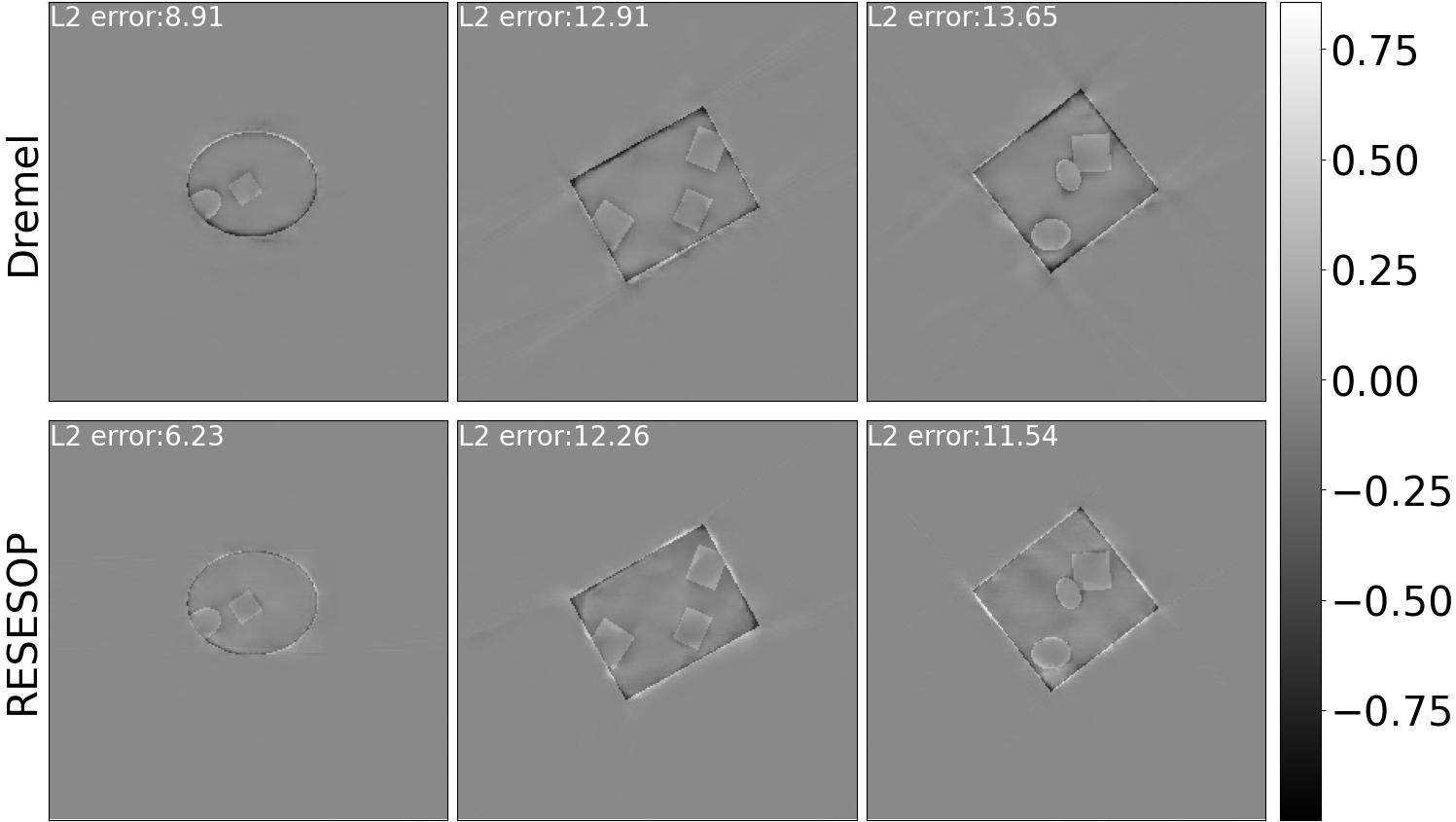}
    \caption{Image reconstructions $x_\mathrm{reco}$ and differences to ground truth using the non-trained reconstruction methods $\mathcal{T}$ combined with \textbf{CiNNRes} post-processing on perturbed fan beam data for the phantoms in \cref{phantoms}.}
    \label{fig:cinn_fan}
\end{figure}

\begin{figure}
    \centering
    Reconstructions\\[2mm] 
    \includegraphics[width=\linewidth]{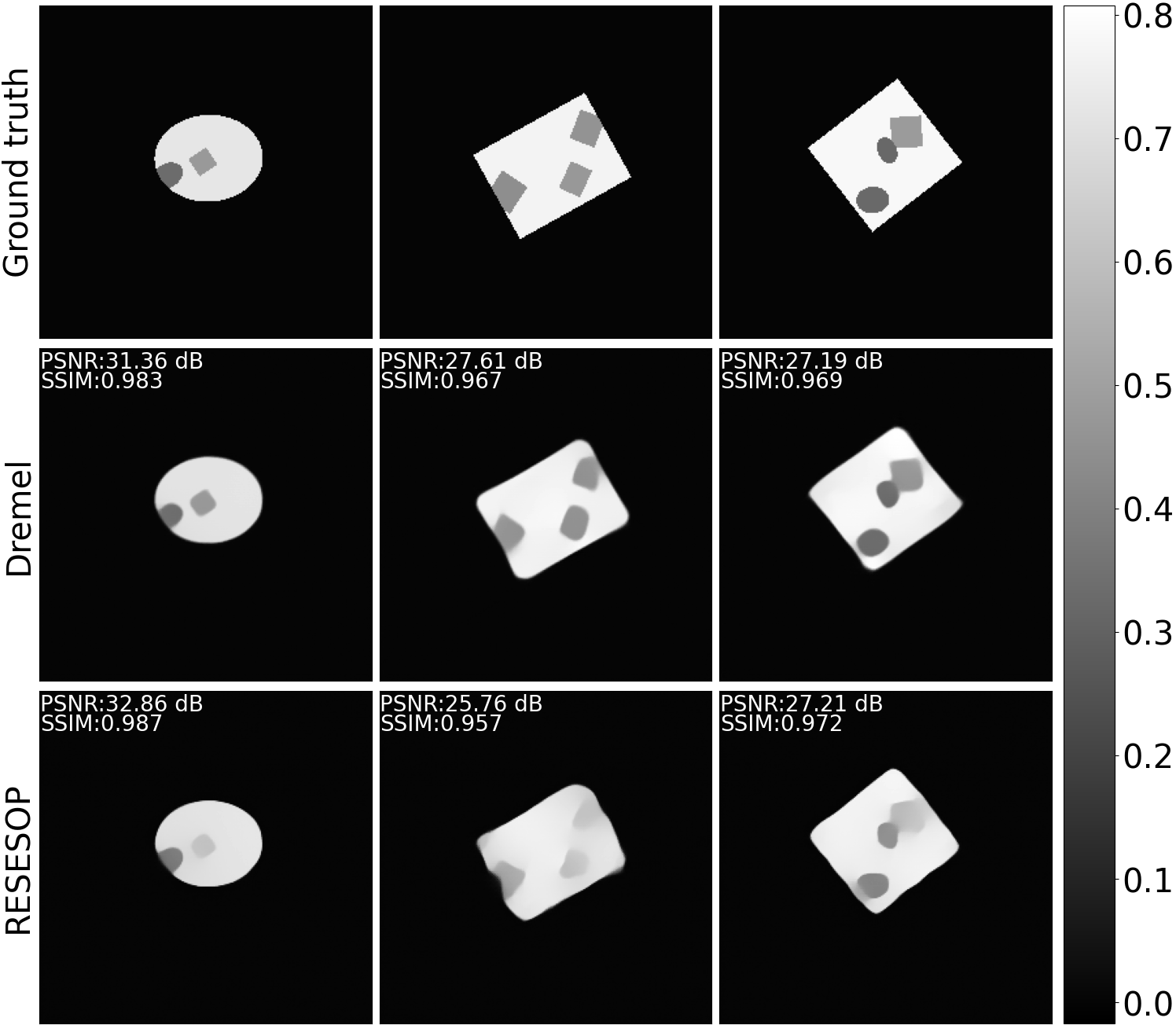}\\
    Differences\\[2mm]     
        \includegraphics[width=\linewidth]{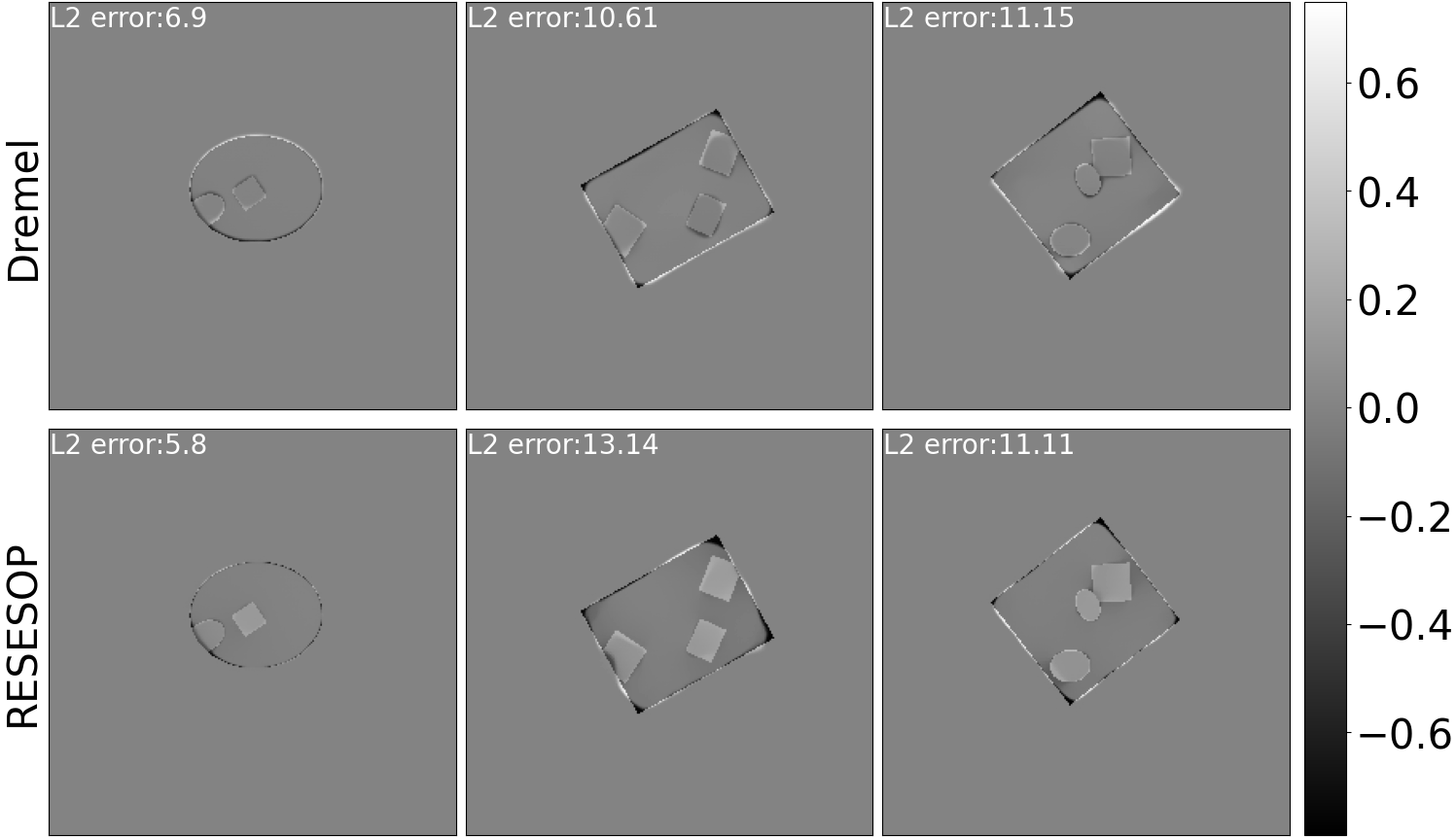}
    \caption{Image reconstructions $x_\mathrm{reco}$ and differences to ground truth using the non-trained reconstruction methods $\mathcal{T}$ combined with \textbf{CiUNetRes} post-processing on perturbed fan beam data for the phantoms in \cref{phantoms}.}
    \label{fig:iunet_fan}
\end{figure}

\clearpage
\onecolumn
\section{Additional quantitative results - perturbed parallel beam} 
\label{app:quant_perturbed_parallel}
\begin{table*}[ht]
\caption{Quantitative evaluation (PSNR (top) and SSIM (bottom)) for mixed training and testing inputs to the post-processing {\textbf{CiNN}} evaluated on the perturbed parallel beam test data for measurements from perturbed phantoms. The rows encode the reconstruction method used for training the columns encode the reconstruction method used for quantitative evaluation on the test set.}
\include{tables/input_mixing_psnr_perturbed_parallel_beam_test_data_cinn.tex}
\include{tables/input_mixing_ssim_perturbed_parallel_beam_test_data_cinn.tex}
\end{table*}

\begin{table*}
\caption{Quantitative evaluation (PSNR (top) and SSIM (bottom)) for mixed training and testing inputs to the post-processing {\textbf{CiNNRes}} evaluated on the perturbed parallel beam test data for measurements from perturbed phantoms. The rows encode the reconstruction method used for training the columns encode the reconstruction method used for quantitative evaluation on the test set.}
\include{tables/input_mixing_psnr_perturbed_parallel_beam_test_data_cinnres.tex}
\include{tables/input_mixing_ssim_perturbed_parallel_beam_test_data_cinnres.tex}
\end{table*}

\begin{table*}
\caption{Quantitative evaluation (PSNR (top) and SSIM (bottom)) for mixed training and testing inputs to the post-processing {\textbf{CiUnet}} evaluated on the perturbed parallel beam test data for measurements from perturbed phantoms. The rows encode the reconstruction method used for training the columns encode the reconstruction method used for quantitative evaluation on the test set.}
\include{tables/input_mixing_psnr_perturbed_parallel_beam_test_data_iunet.tex}
\include{tables/input_mixing_ssim_perturbed_parallel_beam_test_data_iunet.tex}
\end{table*}

\begin{table*}
\caption{Quantitative evaluation (PSNR (top) and SSIM (bottom)) for mixed training and testing inputs to the post-processing {\textbf{CiUnetRes}} evaluated on the perturbed parallel beam test data for measurements from perturbed phantoms. The rows encode the reconstruction method used for training the columns encode the reconstruction method used for quantitative evaluation on the test set.}
\include{tables/input_mixing_psnr_perturbed_parallel_beam_test_data_iunetres.tex}
\include{tables/input_mixing_ssim_perturbed_parallel_beam_test_data_iunetres.tex}
\end{table*}

}

\end{document}

%% file: tables/geometry.tex
\begin{tabular}{lcccc}
\rowcolor[HTML]{DAE8FC} 
\textbf{Geometry}  & \textbf{$S$} & \textbf{$K$} & \textbf{$L$} & \textbf{$D$} \\
Parallel   & $255^2$ & 567 & 363 & \\
Fan        & $255^2$ & 133 & 723 & $7773.4$
\end{tabular}

%% file: tables/unperturbed_parallel_beam_test_data.tex
\centering
\scalebox{0.9}{
\begin{tabular}{lll}
\rowcolor[HTML]{DAE8FC} 
\textbf{$\mathcal{T}$ + $F_\theta$}  & \textbf{PSNR [dB]} & \textbf{SSIM} \\
\hline
\textbf{ No $F_\theta$ } & & \\
\hline
FBP & $40.10 \pm 2.60$ & $0.967 \pm 0.022$ \\
Dremel & $35.29 \pm 2.81$ & $0.945 \pm 0.028$ \\
RESESOP ($\eta=0$) & $43.96 \pm 2.34$ & $0.986 \pm 0.006$ \\
\hline
\textbf{ + $F_\theta$ } & & \\
\hline
FBP + UNet & $71.76 \pm 3.18$ & $0.999 \pm 9.027 \cdot 10^{-6}$ \\
FBP + CiNN & $43.35 \pm 1.47$ & $0.869 \pm 0.004$ \\
FBP + CiNNRes & $55.07 \pm 2.06$ & $0.997 \pm 0.003$ \\
FBP + CiUNet & $46.21 \pm 0.25$ & $0.852 \pm 0.005$ \\
FBP + CiUNetRes & $63.53 \pm 1.27$ & $0.999 \pm 1.784 \cdot 10^{-5}$ \\
\end{tabular}
}

%% file: tables/perturbed_parallel_beam_test_data.tex
\centering
\scalebox{0.9}{
\begin{tabular}{lll}
\rowcolor[HTML]{DAE8FC} 
\textbf{$\mathcal{T}$ + $F_\theta$}  & \textbf{PSNR [dB]} & \textbf{SSIM} \\
\hline
\textbf{ $F_\theta = Id$ } & & \\
\hline
FBP & $27.94 \pm 2.81$ & $0.510 \pm 0.103$ \\
Dremel & $30.98 \pm 3.01$ & $0.829 \pm 0.076$ \\
RESESOP & $30.65 \pm 2.74$ & $0.868 \pm 0.032$ \\
\hline
\textbf{ $F_\theta$ = UNet } & & \\
\hline
FBP & $\mathbf{33.04 \pm 3.87}$ & $\mathbf{0.984 \pm 0.008}$ \\
Dremel & $\mathbf{31.67 \pm 3.79}$ & $\mathbf{0.980 \pm 0.010}$ \\
RESESOP & $\underline{\mathbf{38.17 \pm 3.36}}$ & $\underline{\mathbf{0.993 \pm 0.004}}$ \\
RESESOP 1.2$\eta$ & $\mathbf{37.75 \pm 3.40}$ & $\mathbf{0.992 \pm 0.004}$ \\
RESESOP 0.8$\eta$ & $\mathbf{37.93 \pm 3.44}$ & $\mathbf{0.993 \pm 0.004}$ \\
\hline
\textbf{ $F_\theta$ = CiNN } & & \\
\hline
FBP & {\color{teal}$30.60 \pm 2.31$} & $0.644 \pm 0.009$ \\
Dremel & $30.69 \pm 2.69$ & $0.647 \pm 0.009$ \\
RESESOP & $32.58 \pm 2.16$ & $0.608 \pm 0.010$ \\
RESESOP 1.2$\eta$ & $31.51 \pm 2.25$ & $0.586 \pm 0.011$ \\
RESESOP 0.8$\eta$ & $31.89 \pm 2.16$ & $0.562 \pm 0.012$ \\
\hline
\textbf{ $F_\theta$ = CiNNRes } & & \\
\hline
FBP & $28.82 \pm 2.39$ & $0.862 \pm 0.032$ \\
Dremel & $30.40 \pm 2.74$ & $0.919 \pm 0.018$ \\
RESESOP & $33.30 \pm 2.60$ & $0.964 \pm 0.011$ \\
RESESOP 1.2$\eta$ & {\color{teal}$33.11 \pm 2.65$} & $0.964 \pm 0.012$ \\
RESESOP 0.8$\eta$ & $32.93 \pm 2.59$ & $0.940 \pm 0.017$ \\
\hline
\textbf{ $F_\theta$ = CiUNet } & & \\
\hline
FBP & $22.26 \pm 1.42$ & $0.761 \pm 0.015$ \\
Dremel & $21.98 \pm 1.57$ & $0.686 \pm 0.017$ \\
RESESOP & $22.36 \pm 1.64$ & $0.737 \pm 0.016$ \\
RESESOP 1.2$\eta$ & $25.47 \pm 1.25$ & $0.755 \pm 0.014$ \\
RESESOP 0.8$\eta$ & $27.52 \pm 1.84$ & $0.768 \pm 0.009$ \\
\hline
\textbf{ $F_\theta$ = CiUNetRes } & & \\
\hline
FBP & $25.81 \pm 1.37$ & {\color{teal}$0.931 \pm 0.016$} \\
Dremel & {\color{teal}$30.94 \pm 2.80$} & {\color{teal}$0.975 \pm 0.009$} \\
RESESOP & {\color{teal}$33.87 \pm 2.86$} & {\color{teal}$0.984 \pm 0.006$} \\
RESESOP 1.2$\eta$ & $30.13 \pm 2.49$ & {\color{teal}$0.974 \pm 0.011$} \\
RESESOP 0.8$\eta$ & {\color{teal}$33.50 \pm 2.99$} & {\color{teal}$0.983 \pm 0.007$} \\
\end{tabular}
}

%% file: tables/perturbed_fan_beam_test_data.tex
\centering
\scalebox{0.9}{
\begin{tabular}{lll}
\rowcolor[HTML]{DAE8FC} 
\textbf{$\mathcal{T}$ + $F_\theta$}  & \textbf{PSNR [dB]} & \textbf{SSIM} \\
\hline
\textbf{ $F_\theta = Id$ } & & \\
\hline
Dremel & $30.41 \pm 3.01$ & $0.778 \pm 0.117$ \\
RESESOP & $30.13 \pm 2.87$ & $0.859 \pm 0.039$ \\
\hline
\textbf{ $F_\theta$ = UNet } & & \\
\hline
Dremel & $30.37 \pm 3.77$ & $\mathbf{0.975 \pm 0.012}$ \\
RESESOP & $\underline{\mathbf{37.33 \pm 3.61}}$ & $\underline{\mathbf{0.992 \pm 0.005}}$ \\
\hline
\textbf{ $F_\theta$ = CiNN } & & \\
\hline
Dremel & $29.83 \pm 2.52$ & $0.511 \pm 0.015$ \\
RESESOP & $31.19 \pm 2.19$ & $0.506 \pm 0.015$ \\
\hline
\textbf{ $F_\theta$ = CiNNRes } & & \\
\hline
Dremel & $29.55 \pm 2.69$ & $0.896 \pm 0.027$ \\
RESESOP & {\color{teal}$32.34 \pm 2.66$} & $0.939 \pm 0.019$ \\
\hline
\textbf{ $F_\theta$ = CiUNet } & & \\
\hline
Dremel & $26.21 \pm 1.77$ & $0.440 \pm 0.015$ \\
RESESOP & $25.62 \pm 1.23$ & $0.775 \pm 0.017$ \\
\hline
\textbf{ $F_\theta$ = CiUNetRes } & & \\
\hline
Dremel & {\color{teal}$\mathbf{31.01 \pm 3.26}$} & {\color{teal}$0.974 \pm 0.010$} \\
RESESOP & $32.03 \pm 2.77$ & {\color{teal}$0.978 \pm 0.010$} \\
\end{tabular}
}

%% file: tables/input_mixing_psnr_perturbed_parallel_beam_test_data_unet.tex
\centering
\scalebox{0.9}{
\begin{tabular}{llllll}
\rowcolor[HTML]{DAE8FC} 
\textbf{UNet/Reconstruction} & FBP & Dremel & RESESOP & RESESOP 1.2$\eta$ & RESESOP 0.8$\eta$ \\
FBP & $33.04 \pm 3.87$ & $31.16 \pm 3.64$ & $32.06 \pm 2.68$ & $31.37 \pm 2.79$ & $32.68 \pm 2.68$ \\
Dremel & $29.53 \pm 3.00$ & $31.67 \pm 3.79$ & $31.12 \pm 2.76$ & $30.63 \pm 2.80$ & $31.59 \pm 2.75$ \\
RESESOP & $26.41 \pm 2.21$ & $29.90 \pm 2.99$ & $38.17 \pm 3.36$ & $37.45 \pm 3.55$ & $36.43 \pm 3.71$ \\
RESESOP 1.2$\eta$ & $28.00 \pm 2.97$ & $30.17 \pm 3.08$ & $37.74 \pm 3.39$ & $37.75 \pm 3.40$ & $36.03 \pm 3.73$ \\
RESESOP 0.8$\eta$ & $27.64 \pm 3.07$ & $30.44 \pm 3.13$ & $37.03 \pm 3.75$ & $35.91 \pm 3.79$ & $37.93 \pm 3.44$ \\
\end{tabular}
}

%% file: tables/input_mixing_ssim_perturbed_parallel_beam_test_data_unet.tex
\centering
\scalebox{0.9}{
\begin{tabular}{llllll}
\rowcolor[HTML]{DAE8FC} 
\textbf{UNet/Reconstruction} & FBP & Dremel & RESESOP & RESESOP 1.2$\eta$ & RESESOP 0.8$\eta$ \\
FBP & $0.984 \pm 0.008$ & $0.979 \pm 0.010$ & $0.978 \pm 0.008$ & $0.975 \pm 0.009$ & $0.981 \pm 0.006$ \\
Dremel & $0.969 \pm 0.012$ & $0.980 \pm 0.010$ & $0.971 \pm 0.012$ & $0.968 \pm 0.013$ & $0.973 \pm 0.011$ \\
RESESOP & $0.958 \pm 0.011$ & $0.975 \pm 0.009$ & $0.993 \pm 0.004$ & $0.992 \pm 0.004$ & $0.990 \pm 0.005$ \\
RESESOP 1.2$\eta$ & $0.932 \pm 0.038$ & $0.976 \pm 0.009$ & $0.992 \pm 0.004$ & $0.992 \pm 0.004$ & $0.990 \pm 0.005$ \\
RESESOP 0.8$\eta$ & $0.958 \pm 0.019$ & $0.978 \pm 0.009$ & $0.991 \pm 0.005$ & $0.989 \pm 0.006$ & $0.993 \pm 0.004$ \\
\end{tabular}
}

%% file: tables/input_mixing_psnr_perturbed_parallel_beam_test_data_cinn.tex
\centering
\scalebox{0.9}{
\begin{tabular}{llllll}
\rowcolor[HTML]{DAE8FC}
\textbf{CiNN/Reconstruction} & FBP & Dremel & RESESOP & RESESOP 1.2$\eta$ & RESESOP 0.8$\eta$ \\
FBP & $30.60 \pm 2.31$ & $30.54 \pm 2.39$ & $31.20 \pm 2.10$ & $30.86 \pm 2.22$ & $31.47 \pm 2.02$ \\
Dremel & $29.75 \pm 2.74$ & $30.69 \pm 2.69$ & $30.66 \pm 2.34$ & $30.24 \pm 2.43$ & $31.02 \pm 2.29$ \\
RESESOP & $29.02 \pm 2.69$ & $30.07 \pm 2.53$ & $32.58 \pm 2.16$ & $32.35 \pm 2.26$ & $32.57 \pm 2.16$ \\
RESESOP 1.2$\eta$ & $28.94 \pm 2.72$ & $29.75 \pm 2.52$ & $31.67 \pm 2.21$ & $31.51 \pm 2.25$ & $31.72 \pm 2.21$ \\
RESESOP 0.8$\eta$ & $29.28 \pm 2.81$ & $30.14 \pm 2.63$ & $31.64 \pm 2.19$ & $31.33 \pm 2.28$ & $31.89 \pm 2.16$ \\
\end{tabular}
}

%% file: tables/input_mixing_ssim_perturbed_parallel_beam_test_data_cinn.tex
\centering
\scalebox{0.9}{
\begin{tabular}{llllll}
\rowcolor[HTML]{DAE8FC}
\textbf{CiNN/Reconstruction} & FBP & Dremel & RESESOP & RESESOP 1.2$\eta$ & RESESOP 0.8$\eta$ \\
FBP & $0.644 \pm 0.009$ & $0.644 \pm 0.010$ & $0.643 \pm 0.010$ & $0.641 \pm 0.010$ & $0.645 \pm 0.010$ \\
Dremel & $0.637 \pm 0.011$ & $0.647 \pm 0.009$ & $0.643 \pm 0.008$ & $0.641 \pm 0.008$ & $0.645 \pm 0.008$ \\
RESESOP & $0.596 \pm 0.011$ & $0.601 \pm 0.011$ & $0.608 \pm 0.010$ & $0.607 \pm 0.010$ & $0.608 \pm 0.011$ \\
RESESOP 1.2$\eta$ & $0.576 \pm 0.011$ & $0.581 \pm 0.012$ & $0.587 \pm 0.011$ & $0.586 \pm 0.011$ & $0.587 \pm 0.011$ \\
RESESOP 0.8$\eta$ & $0.550 \pm 0.013$ & $0.556 \pm 0.013$ & $0.561 \pm 0.012$ & $0.560 \pm 0.012$ & $0.562 \pm 0.012$ \\
\end{tabular}
}

%% file: tables/input_mixing_psnr_perturbed_parallel_beam_test_data_cinnres.tex
\centering
\scalebox{0.9}{
\begin{tabular}{llllll}
\rowcolor[HTML]{DAE8FC}
\textbf{CiNNRes/Reconstruction} & FBP & Dremel & RESESOP & RESESOP 1.2$\eta$ & RESESOP 0.8$\eta$ \\
FBP & $28.82 \pm 2.39$ & $28.76 \pm 2.37$ & $28.77 \pm 2.11$ & $28.52 \pm 2.19$ & $29.01 \pm 2.07$ \\
Dremel & $29.09 \pm 2.89$ & $30.40 \pm 2.74$ & $29.92 \pm 2.45$ & $29.51 \pm 2.53$ & $30.32 \pm 2.39$ \\
RESESOP & $30.05 \pm 2.96$ & $31.07 \pm 2.95$ & $33.30 \pm 2.60$ & $32.98 \pm 2.68$ & $33.37 \pm 2.59$ \\
RESESOP 1.2$\eta$ & $29.26 \pm 2.93$ & $30.68 \pm 2.94$ & $33.30 \pm 2.61$ & $33.11 \pm 2.65$ & $33.29 \pm 2.62$ \\
RESESOP 0.8$\eta$ & $28.22 \pm 2.91$ & $30.21 \pm 2.89$ & $32.57 \pm 2.64$ & $32.17 \pm 2.71$ & $32.93 \pm 2.59$ \\
\end{tabular}
}

%% file: tables/input_mixing_ssim_perturbed_parallel_beam_test_data_cinnres.tex
\centering
\scalebox{0.9}{
\begin{tabular}{llllll}
\rowcolor[HTML]{DAE8FC}
\textbf{CiNNRes/Reconstruction} & FBP & Dremel & RESESOP & RESESOP 1.2$\eta$ & RESESOP 0.8$\eta$ \\
FBP & $0.862 \pm 0.032$ & $0.878 \pm 0.025$ & $0.880 \pm 0.023$ & $0.880 \pm 0.024$ & $0.872 \pm 0.029$ \\
Dremel & $0.793 \pm 0.082$ & $0.919 \pm 0.018$ & $0.893 \pm 0.027$ & $0.886 \pm 0.030$ & $0.887 \pm 0.034$ \\
RESESOP & $0.874 \pm 0.051$ & $0.955 \pm 0.015$ & $0.964 \pm 0.011$ & $0.962 \pm 0.012$ & $0.961 \pm 0.013$ \\
RESESOP 1.2$\eta$ & $0.795 \pm 0.067$ & $0.933 \pm 0.027$ & $0.965 \pm 0.011$ & $0.964 \pm 0.012$ & $0.961 \pm 0.014$ \\
RESESOP 0.8$\eta$ & $0.584 \pm 0.097$ & $0.822 \pm 0.070$ & $0.935 \pm 0.017$ & $0.929 \pm 0.020$ & $0.940 \pm 0.017$ \\
\end{tabular}
}

%% file: tables/input_mixing_psnr_perturbed_parallel_beam_test_data_iunet.tex
\centering
\scalebox{0.9}{
\begin{tabular}{llllll}
\rowcolor[HTML]{DAE8FC} 
\textbf{CiUNet/Reconstruction} & FBP & Dremel & RESESOP & RESESOP 1.2$\eta$ & RESESOP 0.8$\eta$ \\
FBP & $22.26 \pm 1.42$ & $22.38 \pm 1.42$ & $21.82 \pm 1.49$ & $21.68 \pm 1.52$ & $21.92 \pm 1.47$ \\
Dremel & $21.80 \pm 1.60$ & $21.98 \pm 1.57$ & $21.50 \pm 1.66$ & $21.38 \pm 1.69$ & $21.57 \pm 1.64$ \\
RESESOP & $22.55 \pm 1.63$ & $22.78 \pm 1.66$ & $22.36 \pm 1.64$ & $22.21 \pm 1.66$ & $22.46 \pm 1.62$ \\
RESESOP 1.2$\eta$ & $25.95 \pm 1.43$ & $26.39 \pm 1.43$ & $25.74 \pm 1.22$ & $25.47 \pm 1.25$ & $25.97 \pm 1.20$ \\
RESESOP 0.8$\eta$ & $27.50 \pm 2.22$ & $27.78 \pm 2.20$ & $27.28 \pm 1.83$ & $27.03 \pm 1.84$ & $27.52 \pm 1.84$ \\
\end{tabular}
}

%% file: tables/input_mixing_ssim_perturbed_parallel_beam_test_data_iunet.tex
\centering
\scalebox{0.9}{
\begin{tabular}{llllll}
\rowcolor[HTML]{DAE8FC} 
\textbf{CiUNet/Reconstruction} & FBP & Dremel & RESESOP & RESESOP 1.2$\eta$ & RESESOP 0.8$\eta$ \\
FBP & $0.761 \pm 0.015$ & $0.762 \pm 0.015$ & $0.757 \pm 0.014$ & $0.756 \pm 0.014$ & $0.758 \pm 0.014$ \\
Dremel & $0.685 \pm 0.017$ & $0.686 \pm 0.017$ & $0.682 \pm 0.016$ & $0.681 \pm 0.017$ & $0.683 \pm 0.016$ \\
RESESOP & $0.736 \pm 0.017$ & $0.740 \pm 0.017$ & $0.737 \pm 0.016$ & $0.736 \pm 0.017$ & $0.738 \pm 0.016$ \\
RESESOP 1.2$\eta$ & $0.757 \pm 0.012$ & $0.760 \pm 0.013$ & $0.756 \pm 0.014$ & $0.755 \pm 0.014$ & $0.758 \pm 0.014$ \\
RESESOP 0.8$\eta$ & $0.766 \pm 0.009$ & $0.768 \pm 0.009$ & $0.766 \pm 0.009$ & $0.765 \pm 0.008$ & $0.768 \pm 0.009$ \\
\end{tabular}
}

%% file: tables/input_mixing_psnr_perturbed_parallel_beam_test_data_iunetres.tex
\centering
\scalebox{0.9}{
\begin{tabular}{llllll}
\rowcolor[HTML]{DAE8FC} 
\textbf{CiUNetRes/Reconstruction} & FBP & Dremel & RESESOP & RESESOP 1.2$\eta$ & RESESOP 0.8$\eta$ \\
FBP & $25.81 \pm 1.37$ & $26.15 \pm 1.32$ & $25.57 \pm 1.28$ & $25.34 \pm 1.34$ & $25.71 \pm 1.23$ \\
Dremel & $30.00 \pm 2.89$ & $30.94 \pm 2.80$ & $30.95 \pm 1.91$ & $30.64 \pm 1.98$ & $31.14 \pm 1.91$ \\
RESESOP & $29.08 \pm 3.09$ & $30.76 \pm 3.18$ & $33.87 \pm 2.86$ & $33.45 \pm 2.90$ & $34.13 \pm 2.88$ \\
RESESOP 1.2$\eta$ & $28.44 \pm 2.77$ & $29.78 \pm 2.84$ & $30.56 \pm 2.51$ & $30.13 \pm 2.49$ & $31.08 \pm 2.53$ \\
RESESOP 0.8$\eta$ & $29.16 \pm 3.14$ & $30.72 \pm 3.27$ & $32.95 \pm 2.95$ & $32.52 \pm 2.96$ & $33.50 \pm 2.99$ \\
\end{tabular}
}

%% file: tables/input_mixing_ssim_perturbed_parallel_beam_test_data_iunetres.tex
\centering
\scalebox{0.9}{
\begin{tabular}{llllll}
\rowcolor[HTML]{DAE8FC} 
\textbf{CiUNetRes/Reconstruction} & FBP & Dremel & RESESOP & RESESOP 1.2$\eta$ & RESESOP 0.8$\eta$ \\
FBP & $0.931 \pm 0.016$ & $0.936 \pm 0.016$ & $0.933 \pm 0.016$ & $0.931 \pm 0.017$ & $0.935 \pm 0.016$ \\
Dremel & $0.965 \pm 0.013$ & $0.975 \pm 0.009$ & $0.972 \pm 0.008$ & $0.971 \pm 0.010$ & $0.971 \pm 0.009$ \\
RESESOP & $0.891 \pm 0.038$ & $0.965 \pm 0.015$ & $0.984 \pm 0.006$ & $0.982 \pm 0.007$ & $0.984 \pm 0.006$ \\
RESESOP 1.2$\eta$ & $0.884 \pm 0.041$ & $0.959 \pm 0.016$ & $0.976 \pm 0.010$ & $0.974 \pm 0.011$ & $0.978 \pm 0.009$ \\
RESESOP 0.8$\eta$ & $0.937 \pm 0.029$ & $0.975 \pm 0.010$ & $0.981 \pm 0.007$ & $0.980 \pm 0.008$ & $0.983 \pm 0.007$ \\
\end{tabular}
}

%% file: main_arxiv.bbl
\begin{thebibliography}{10}
\providecommand{\url}[1]{#1}
\csname url@samestyle\endcsname
\providecommand{\newblock}{\relax}
\providecommand{\bibinfo}[2]{#2}
\providecommand{\BIBentrySTDinterwordspacing}{\spaceskip=0pt\relax}
\providecommand{\BIBentryALTinterwordstretchfactor}{4}
\providecommand{\BIBentryALTinterwordspacing}{\spaceskip=\fontdimen2\font plus
\BIBentryALTinterwordstretchfactor\fontdimen3\font minus
  \fontdimen4\font\relax}
\providecommand{\BIBforeignlanguage}[2]{{%
\expandafter\ifx\csname l@#1\endcsname\relax
\typeout{** WARNING: IEEEtran.bst: No hyphenation pattern has been}%
\typeout{** loaded for the language `#1'. Using the pattern for}%
\typeout{** the default language instead.}%
\else
\language=\csname l@#1\endcsname
\fi
#2}}
\providecommand{\BIBdecl}{\relax}
\BIBdecl

\bibitem{Gleich2005}
B.~Gleich and J.~Weizenecker, ``Tomographic imaging using the nonlinear
  response of magnetic particles,'' \emph{Nature}, vol. 435, no. 7046, pp.
  1214--1217, 2005.

\bibitem{rtkms17}
A.~Ruhlandt, M.~T\"opperwien, M.~Krenkel, R.~Mokso, and T.~Salditt, ``{Four
  dimensional material movies: High speed phase-contrast tomography by
  backprojection along dynamically curved paths},'' \emph{Scientific Reports},
  vol.~7, p. 6487, 2017.

\bibitem{youssef_2013}
\BIBentryALTinterwordspacing
S.~Youssef, H.~Deschamps, J.~Dautriat, E.~Rosenberg, R.~Oughanem, E.~Maire, and
  R.~Mokso, ``{4{D} imaging of fluid flow dynamics in natural porous media with
  ultra-fast {X}-ray microtomography},'' in \emph{{International {S}ymposium of
  the {S}{C}{A}, {N}apa {V}alley, {C}alifornia}}, vol. 176, 2013. [Online].
  Available: \url{https://hal.archives-ouvertes.fr/hal-01538218}
\BIBentrySTDinterwordspacing

\bibitem{golub1999tikhonov}
G.~H. Golub, P.~C. Hansen, and D.~P. O'Leary, ``Tikhonov regularization and
  total least squares,'' \emph{SIAM journal on matrix analysis and
  applications}, vol.~21, no.~1, pp. 185--194, 1999.

\bibitem{Bleyer:2013cw}
I.~R. Bleyer and R.~Ramlau, ``{A double regularization approach for inverse
  problems with noisy data and inexact operator},'' \emph{Inverse Problems},
  vol.~29, no.~2, pp. 025\,004--17.

\bibitem{KluthBathkeJiangMaass2020}
T.~Kluth, C.~Bathke, M.~Jiang, and P.~Maass, ``Joint super-resolution image
  reconstruction and parameter identification in imaging operator: analysis of
  bilinear operator equations, numerical solution, and application to magnetic
  particle imaging,'' \emph{Inverse Problems}, vol.~36, no.~12, p. 124006,
  2020.

\bibitem{Blanke_2020}
S.~E. Blanke, B.~N. Hahn, and A.~Wald, ``Inverse problems with inexact forward
  operator: iterative regularization and application in dynamic imaging,''
  \emph{Inverse Problems}, vol.~36, no.~12, p. 124001, 10 2020.

\bibitem{lunz_learned_2021}
S.~Lunz, A.~Hauptmann, T.~Tarvainen, C.-B. Schönlieb, and S.~Arridge, ``On
  learned operator correction in inverse problems,'' \emph{SIAM Journal on
  Imaging Sciences}, vol.~14, no.~1, pp. 92--127, 2021, publisher: SIAM.

\bibitem{gregor_learning_2010}
K.~Gregor and Y.~LeCun, ``Learning fast approximations of sparse coding,'' in
  \emph{Proceedings of the 27th international conference on international
  conference on machine learning}, 2010, pp. 399--406.

\bibitem{adler2017solving}
J.~Adler and O.~{\"O}ktem, ``Solving ill-posed inverse problems using iterative
  deep neural networks,'' \emph{Inverse Problems}, vol.~33, no.~12, p. 124007,
  2017.

\bibitem{adler_learned_2018}
J.~Adler and O.~Öktem, ``Learned primal-dual reconstruction,'' \emph{IEEE
  Transactions on Medical Imaging}, vol.~37, no.~6, pp. 1322--1332, 2018.

\bibitem{Dremel_2018}
K.~Dremel, ``Modellbildung des messprozesses und umsetzung eines
  modellbasierten iterativen l{\"o}sungsverfahrens der
  schnittbild-rekonstruktion f{\"u}r die r{\"o}ntgen-computertomographie,''
  doctoralthesis, Universit{\"a}t W{\"u}rzburg, 2018.

\bibitem{Wang_2019}
S.~Wang, J.~Liu, Y.~Li, J.~Chen, Y.~Guan, and L.~Zhu, ``{Jitter correction for
  transmission X-ray microscopy via measurement of geometric moments},''
  \emph{Journal of Synchrotron Radiation}, vol.~26, no.~5, pp. 1808--1814, Sep
  2019.

\bibitem{natterer_mathematical_2001}
F.~Natterer and F.~Wübbeling, \emph{Mathematical {Methods} in {Image}
  {Reconstruction}}.\hskip 1em plus 0.5em minus 0.4em\relax Society for
  Industrial and Applied Mathematics, 2001.

\bibitem{kaczmarz1937angenaherte}
S.~Kaczmarz, ``Angenäherte auflösung von systemen linearer gleichungen,''
  \emph{Bull. Int. Acad. Sci. Pologne, A}, vol.~35, pp. 355--357, 1937.

\bibitem{Gordon1970art}
R.~Gordon, R.~Bender, and G.~T. Herman, ``Algebraic reconstruction techniques
  (art) for three-dimensional electron microscopy and x-ray photography,''
  \emph{Journal of Theoretical Biology}, vol.~29, no.~3, pp. 471--481, 1970.

\bibitem{narkiss2005sequential}
G.~Narkiss and M.~Zibulevsky, \emph{Sequential subspace optimization method for
  large-scale unconstrained problems}.\hskip 1em plus 0.5em minus 0.4em\relax
  Technion-IIT, Department of Electrical Engineering, 2005.

\bibitem{Wang2020dl_tomo_reco}
G.~Wang, J.~C. Ye, and B.~De~Man, ``Deep learning for tomographic image
  reconstruction,'' \emph{Nature Machine Intelligence}, vol.~2, no.~12, pp.
  737--748, 2020.

\bibitem{VillarragaGomez2022microscope_ct_dl}
H.~Villarraga-G{\'o}mez, M.~N. Rad, M.~Andrew, A.~Andreyev, R.~Sanapala,
  L.~Omlor, and C.~G. vom Hagen, ``Improving throughput and image quality of
  high-resolution 3d x-ray microscopes using deep learning reconstruction
  techniques,'' in \emph{11th Conference on Industrial Computed Tomography
  (iCT)}, 2022, pp. 8--11.

\bibitem{AnayaIsaza2021overview_dl_medical}
A.~Anaya-Isaza, L.~Mera-Jiménez, and M.~Zequera-Diaz, ``An overview of deep
  learning in medical imaging,'' \emph{Informatics in Medicine Unlocked},
  vol.~26, p. 100723, 2021.

\bibitem{Rafiei2023cv_industrial_security}
M.~Rafiei, J.~Raitoharju, and A.~Iosifidis, ``Computer vision on x-ray data in
  industrial production and security applications: A comprehensive survey,''
  vol.~11, pp. 2445--2477.

\bibitem{Szczykutowicz2022review_ct_reco}
T.~P. Szczykutowicz, G.~V. Toia, A.~Dhanantwari, and B.~Nett, ``A review of
  deep learning ct reconstruction: Concepts, limitations, and promise in
  clinical practice,'' \emph{Current Radiology Reports}, vol.~10, no.~9, pp.
  101--115, 2022.

\bibitem{jin_deep_2017}
K.~H. Jin, M.~T. McCann, E.~Froustey, and M.~Unser, ``Deep convolutional neural
  network for inverse problems in imaging,'' \emph{IEEE Transactions on Image
  Processing}, vol.~26, no.~9, pp. 4509--4522, 2017, publisher: IEEE.

\bibitem{yang_low-dose_2018}
Q.~Yang, P.~Yan, Y.~Zhang, H.~Yu, Y.~Shi, X.~Mou, M.~K. Kalra, Y.~Zhang,
  L.~Sun, and G.~Wang, ``Low-{Dose} {CT} {Image} {Denoising} {Using} a
  {Generative} {Adversarial} {Network} {With} {Wasserstein} {Distance} and
  {Perceptual} {Loss},'' \emph{IEEE Transactions on Medical Imaging}, vol.~37,
  no.~6, pp. 1348--1357, Jun. 2018.

\bibitem{yuan_sipid_2018}
H.~Yuan, J.~Jia, and Z.~Zhu, ``{SIPID}: {A} deep learning framework for
  sinogram interpolation and image denoising in low-dose {CT} reconstruction.''
  in \emph{15th {IEEE} {International} {Symposium} on {Biomedical} {Imaging},
  {ISBI} 2018, {Washington}, {DC}, {USA}, {April} 4-7, 2018}, 2018, pp.
  1521--1524.

\bibitem{liu_interpreting_2020}
T.~Liu, A.~Chaman, D.~Belius, and I.~Dokmanic, ``Learning {Multiscale}
  {Convolutional} {Dictionaries} for {Image} {Reconstruction}.'' \emph{IEEE
  Trans. Computational Imaging}, vol.~8, pp. 425--437, 2022.

\bibitem{you_ct_2020}
C.~You, W.~Cong, M.~W. Vannier, P.~K. Saha, E.~A. Hoffman, G.~Wang, G.~Li,
  Y.~Zhang, X.~Zhang, H.~Shan, M.~Li, S.~Ju, Z.~Zhao, and Z.~Zhang, ``{CT}
  {Super}-{Resolution} {GAN} {Constrained} by the {Identical}, {Residual}, and
  {Cycle} {Learning} {Ensemble} ({GAN}-{CIRCLE}).'' \emph{IEEE Trans. Medical
  Imaging}, vol.~39, no.~1, pp. 188--203, 2020.

\bibitem{ronneberger_u-net_2015}
O.~Ronneberger, P.~Fischer, and T.~Brox, ``U-{Net}: {Convolutional} {Networks}
  for {Biomedical} {Image} {Segmentation},'' in \emph{Medical {Image}
  {Computing} and {Computer}-{Assisted} {Intervention} – {MICCAI}
  2015}.\hskip 1em plus 0.5em minus 0.4em\relax Cham: Springer International
  Publishing, 2015, pp. 234--241.

\bibitem{kaipio_statistical_2005}
J.~Kaipio and E.~Somersalo, \emph{Statistical and {Computational} {Inverse}
  {Problems}}, 1st~ed., ser. Applied {Mathematical} {Sciences}.\hskip 1em plus
  0.5em minus 0.4em\relax New York, NY: Springer, 2005, no. 160.

\bibitem{tarantola_inverse_1982}
A.~Tarantola, B.~Valette, and {others}, ``Inverse problems = quest for
  information,'' \emph{Journal of Geophysics}, vol.~50, no.~1, pp. 159--170,
  1982.

\bibitem{Ardizzone2019cinn}
L.~Ardizzone, C.~Lüth, J.~Kruse, C.~Rother, and U.~Köthe, ``Guided {Image}
  {Generation} with {Conditional} {Invertible} {Neural} {Networks},''
  \emph{arXiv preprint arXiv:1907.02392}, 2019.

\bibitem{Winkler2019cnf}
\BIBentryALTinterwordspacing
C.~Winkler, D.~Worrall, E.~Hoogeboom, and M.~Welling, ``Learning likelihoods
  with conditional normalizing flows.'' [Online]. Available:
  \url{https://arxiv.org/abs/1912.00042}
\BIBentrySTDinterwordspacing

\bibitem{anantha_padmanabha_solving_2021}
G.~Anantha~Padmanabha and N.~Zabaras, ``Solving inverse problems using
  conditional invertible neural networks,'' \emph{Journal of Computational
  Physics}, vol. 433, p. 110194, May 2021.

\bibitem{denker_conditional_2020}
\BIBentryALTinterwordspacing
A.~Denker, M.~Schmidt, J.~Leuschner, P.~Maass, and J.~Behrmann, ``Conditional
  normalizing flows for low-dose computed tomography image reconstruction,'' in
  \emph{{ICML} workshop on invertible neural networks, normalizing flows, and
  explicit likelihood models, {Vienna}, {Austria}, 18 {July}}, 2020. [Online].
  Available:
  \url{https://invertibleworkshop.github.io/INNF_2020/accepted_papers/index.html}
\BIBentrySTDinterwordspacing

\bibitem{denker_conditional_2021}
A.~Denker, M.~Schmidt, J.~Leuschner, and P.~Maass, ``Conditional {Invertible}
  {Neural} {Networks} for {Medical} {Imaging},'' \emph{Journal of Imaging},
  vol.~7, no.~11, 2021.

\bibitem{dinh_nice_2015}
L.~Dinh, D.~Krueger, and Y.~Bengio, ``{NICE}: {Non}-linear {Independent}
  {Components} {Estimation},'' in \emph{3rd {International} {Conference} on
  {Learning} {Representations}, {ICLR} 2015, {San} {Diego}, {CA}, {USA}, {May}
  7-9, 2015, {Workshop} {Track} {Proceedings}}, 2015.

\bibitem{dinh_density_2017}
L.~Dinh, J.~Sohl-Dickstein, and S.~Bengio, ``Density estimation using {Real}
  {NVP},'' in \emph{5th {International} {Conference} on {Learning}
  {Representations}, {ICLR} 2017, {Toulon}, {France}, {April} 24-26, 2017,
  {Conference} {Track} {Proceedings}}, 2017.

\bibitem{etmann_iunets_2020}
C.~Etmann, R.~Ke, and C.-B. Schönlieb, ``{iUNets}: {Learnable} {Invertible}
  {Up}- and {Downsampling} for {Large}-{Scale} {Inverse} {Problems},'' in
  \emph{30th {IEEE} {International} {Workshop} on {Machine} {Learning} for
  {Signal} {Processing}, {MLSP} 2020, {Espoo}, {Finland}, {September} 21-24,
  2020}.\hskip 1em plus 0.5em minus 0.4em\relax IEEE, 2020, pp. 1--6.

\bibitem{Blaek2019}
P.~Bla{\v{z}}ek, J.~{\v{S}}r{\'a}mek, T.~Zikmund, D.~Kalasov{\'a},
  V.~Hortv{\'\i}k, P.~Klapetek, and J.~Kaiser, ``Voxel size and calibration for
  ct measurements with a small field of view,'' in \emph{Proceedings of the 9th
  Conference on Industrial Computed Tomography (iCT 2019), Padova, Italy},
  2019, pp. 13--15.

\bibitem{Cressa2022fib_sem_correlative}
L.~Cressa, J.~Fell, C.~Pauly, Q.~H. Hoang, F.~M{\"u}cklich, H.-G. Herrmann,
  T.~Wirtz, and S.~Eswara, ``A fib-sem based correlative methodology for x-ray
  nanotomography and secondary ion mass spectrometry: An application example in
  lithium batteries research,'' \emph{Microscopy and Microanalysis}, vol.~28,
  no.~6, pp. 1890--1895, 2022.

\bibitem{adler_operator_2018}
\BIBentryALTinterwordspacing
J.~Adler, H.~Kohr, A.~Ringh, J.~Moosmann, Sbanert, M.~J. Ehrhardt, G.~R. Lee,
  Niinimaki, Bgris, O.~Verdier, J.~Karlsson, Zickert, W.~J. Palenstijn,
  O.~Öktem, C.~Chen, H.~A. Loarca, and M.~Lohmann, ``Odlgroup/{Odl}: {Odl}
  0.7.0,'' Sep. 2018. [Online]. Available:
  \url{https://zenodo.org/record/1442734}
\BIBentrySTDinterwordspacing

\bibitem{aarle_astra_2015}
W.~v. Aarle, W.~J. Palenstijn, J.~D. Beenhouwer, T.~Altantzis, S.~Bals, K.~J.
  Batenburg, and J.~Sijbers, ``The {ASTRA} {Toolbox}: {A} platform for advanced
  algorithm development in electron tomography,'' \emph{Ultramicroscopy}, vol.
  157, pp. 35--47, 2015.

\bibitem{teshima2020coupling}
T.~Teshima, I.~Ishikawa, K.~Tojo, K.~Oono, M.~Ikeda, and M.~Sugiyama,
  ``Coupling-based invertible neural networks are universal diffeomorphism
  approximators,'' \emph{Advances in Neural Information Processing Systems},
  vol.~33, pp. 3362--3373, 2020.

\bibitem{ishikawa2022universal}
I.~Ishikawa, T.~Teshima, K.~Tojo, K.~Oono, M.~Ikeda, and M.~Sugiyama,
  ``Universal approximation property of invertible neural networks,''
  \emph{arXiv preprint arXiv:2204.07415}, 2022.

\bibitem{behrmann_invertible_2019}
J.~Behrmann, W.~Grathwohl, R.~T.~Q. Chen, D.~Duvenaud, and J.-H. Jacobsen,
  ``Invertible {Residual} {Networks},'' in \emph{Proceedings of the 36th
  {International} {Conference} on {Machine} {Learning}}, vol.~97, 2019, pp.
  573--582.

\bibitem{zhang2020approximation}
H.~Zhang, X.~Gao, J.~Unterman, and T.~Arodz, ``Approximation capabilities of
  neural odes and invertible residual networks,'' in \emph{International
  Conference on Machine Learning}.\hskip 1em plus 0.5em minus 0.4em\relax PMLR,
  2020, pp. 11\,086--11\,095.

\end{thebibliography}
